\documentclass[onefignum,onetabnum]{siamart190516}

\usepackage{multirow}
\usepackage{amsmath,amscd,amssymb}
\usepackage{xspace}
\usepackage{cite,alltt}
\usepackage{enumitem}
\setlist[itemize]{topsep=2pt,parsep=0pt,partopsep=2pt}
\setlist[enumerate]{topsep=2pt,parsep=0pt,partopsep=2pt}

\usepackage{mathtools, color}
\usepackage{multirow}
\usepackage{rotating}
\usepackage{subfigure}
\usepackage{epsfig}
\usepackage{url}
\usepackage{graphicx}
\usepackage{stackengine}
\usepackage{mathdots}
\usepackage[english]{babel}
\usepackage{csquotes}
\usepackage{mathrsfs}
\usepackage{overpic}
\usepackage{algorithmic} %algorithm package needed for example 1
\usepackage{algorithm}
\usepackage{blindtext}
\usepackage{mathrsfs}
\usepackage{hyperref}
\usepackage{tikz}
\usetikzlibrary{positioning}
\usepackage[export]{adjustbox}

\graphicspath{{figures/}}

\newsiamremark{remark}{Remark}
\newsiamremark{hypothesis}{Hypothesis}
\crefname{hypothesis}{Hypothesis}{Hypotheses}
\newsiamthm{claim}{Claim}

\newcommand{\ep}{\varepsilon}
\newcommand{\vp}{\mathbf{p}}

\newcommand{\vc}{\mathbf{c}}

\newcommand{\vn}{\mathbf{n}}

\newcommand{\vs}{\mathbf{s}}

\newcommand{\vu}{\mathbf{u}}

\newcommand{\vx}{\mathbf{x}}
\newcommand{\pts}{\mathbf{x}}
\newcommand{\Pts}{X}
\newcommand{\vy}{\mathbf{y}}

\newcommand{\vom}{\boldsymbol{\xi}}

\newcommand{\prad}{\rho}

\newcommand{\Phicurl}{\Phi}

\newcommand{\R}{\mathbb{R}}

\definecolor{matred}{rgb}{0.6350, 0.0780, 0.1840}
\definecolor{matblue}{rgb}{0, 0.4470, 0.7410}
\definecolor{matpurp}{rgb}{0.4940, 0.1840, 0.5560}
\definecolor{matgreen}{rgb}{0, 0.5, 0}

\newcommand{\cfpot}{s}  % cf potential for the vector approximant
\newcommand{\tcfpot}{f}  % target cf potential
\newcommand{\mcfpot}{\overline{\cfpot}}  % cf potential for the vector approximant with mean zero
\newcommand{\icfpot}{\widetilde{\cfpot}}  % cf potential for the vector approximant that exactly interpolates

\definecolor{mygreen}{rgb}{0.13,0.54,0.13}
\definecolor{purple}{rgb}{0.5412,0.1686,0.8863}
\newcommand{\bigO}{\mathcal{O}}
\newcommand{\ds}{\displaystyle}
\newcommand{\revision}[1]{{\color{black}#1}}
\newcommand{\refone}[1]{{\color{black}#1}}
\newcommand{\reftwo}[1]{{\color{black}#1}}

\newcommand{\vertiii}[1]{{\left\vert\kern-0.25ex\left\vert\kern-0.25ex\left\vert #1 \right\vert\kern-0.25ex\right\vert\kern-0.25ex\right\vert}}

\vspace{-3cm}
\headers{Implicit surface reconstructions with CFPU }{K.~P. Drake, E.~J. Fuselier, and G.~B. Wright}
\title{Implicit Surface Reconstruction with a Curl-free Radial Basis Function Partition of Unity Method}
\author{Kathryn P. Drake\thanks{Department of Mathematics, Boise State University, Boise, ID (\email{KathrynDrake@u.boisestate.edu}, \\ \email{gradywright@boisestate.edu})}
\and Edward J. Fuselier\thanks{Department of Mathematics, High Point University, High Point, NC (\email{efuselie@highpoint.edu})}
\and Grady B. Wright\footnotemark[1]}

\begin{document}
\maketitle

\begin{abstract}
Surface reconstruction from a set of scattered points, or a point cloud, has many applications ranging from computer graphics to remote sensing.  We present a new method for this task that produces an implicit surface (zero-level set) approximation for an oriented point cloud using only information about (approximate) normals to the surface.  The technique exploits the fundamental result from vector calculus that the normals to an implicit surface are curl-free.  By using curl-free radial basis function (RBF) interpolation of the normals, we can extract a potential for the vector field whose zero-level surface approximates the point cloud.  We use curl-free RBFs based on polyharmonic splines for this task, since they are free of any shape or support parameters.  To make this technique efficient and able to better represent local sharp features, we combine it with a partition of unity (PU) method. The result is the curl-free partition of unity (CFPU) method.  We show how CFPU can be  adapted to enforce exact interpolation of a point cloud and can be regularized to handle noise in both the normals and the point positions.  Numerical results are presented that demonstrate how the method converges for a known surface as the sampling density increases, how regularization handles noisy data, and how the method performs on various problems found in the literature. 
\end{abstract}

\begin{keywords}
  Implicit surface reconstruction; Meshfree; Point clouds; Radial basis functions; Partition of unity; Curl-free; Potential; Polyharmonic splines; Smoothing splines
\end{keywords}

\begin{AMS}
  41-04, 41A15, 41A29, 65D07, 65D10, 65D15, 65D17, 68U07
\end{AMS}

\section{Introduction}
The problem of reconstructing a surface from a set of unorganized points, or a point cloud, has been used in a variety of applications, including computer graphics, computer-aided design, medical imaging, image processing, manufacturing, and remote sensing. Many common methods developed to address this problem produce an implicit surface (zero-level set) using ``oriented" point clouds, which involve the unstructured points as well as their corresponding normal vectors. In this paper, we present the Curl-free Radial Basis Function Partition of Unity (CFPU) method for this problem.  The method is based on the following result from vector calculus:  If $\tcfpot:\R^{d} \longrightarrow \R$ defines a zero-level set $\mathcal{P}$ (i.e.\ implicit curve for $d=2$ or surface for $d=3$) and $\mathbf{n}$ is a normal vector to $\mathcal{P}$, then $\mathbf{n}$ is curl-free.  This follows since $\mathbf{n}$ is proportional to $\nabla \tcfpot$ and the curl of a gradient field is zero.  Given a set of points $\{\pts_1,\ldots,\pts_N\}\subset\mathcal{P}$ equipped with oriented normal vectors $\{\mathbf{n}_1,\mathbf{n}_2,\dots,\mathbf{n}_N\}$, we seek to recover the potential $\tcfpot$, such that $\nabla \tcfpot \approx \mathbf{n}$ at every point $\pts_j$.   The method we use to recover $\tcfpot$ comes, in part, from a curl-free RBF approximant to the normal vectors. These approximants were introduced in~\cite{AmodeiBenbourhim_1991,DoduRabut_2002_VectorialInterp,Fuselier08} for generating an analytically curl-free approximant of curl-free field from scattered samples of the field. The key to our method lies in the feature that a scalar potential can be extracted from these vector approximants that can be used to approximate the implicit surface $\mathcal{P}$.

Implementing curl-free approximants globally is too computationally expensive, requiring the solution of a $dN$-by-$dN$ system. To bypass this issue, we follow a similar approach to~\cite{DFW2020a} and combine the technique with a partition of unity (PU) method.  This allows the potential $\tcfpot$ to be solved for locally on patches involving \refone{$n \ll N$} points.  These local potentials can then be blended together into a global potential whose zero-level set approximates $\mathcal{P}$. An added benefit of this approach is that it \revision{appears to be better equipped at recovering sharp features}, which many global methods lack.  The CFPU method can also be adapted to enforce exact interpolation of the surface and can be regularized to handle noisy data.  Finally, by using curl-free polyharmonic splines (PHS) and curl-free vector polynomials to construct the local approximants, the method is free from shape and scaling parameters, which are common to other RBF methods and for which good values can be  computationally expensive to determine.

The rest of the paper is structured as follows. For the remainder of Section 1, we briefly discuss relevant surface reconstruction methods and compare them to CFPU. In Section 2, we provide background on curl-free RBF approximation and how it can be used for surface reconstruction. We then introduce the CFPU method in Section 3, along with modifications for exact surface interpolation and regularization. We discuss computational considerations and results of the CFPU method in Sections 4 and 5, respectively. Finally, we provide concluding remarks in Section 6.

\subsection{Relationship to previous work}
Reconstructing a surface from an unorganized point cloud has been extensively studied in literature since the 1990s. Some of the approaches involve signed distance methods~\cite{SR_Unorg,SSR_2011}, RBF-based methods~\cite{Morse_Chen, Ohtake_3D,Suss_Hier,Walder_2006,Closed_Form_HRBF,Macedo_HRBF,Carr_3D,Turk_2002,Wendland:2004,RBF_Med,MeshFree_Mat}, 
%~\cite{Morse_Chen, Ohtake_3D,Pan_HVI,Sam_VCRBF,Suss_Hier,Walder_2006,Closed_Form_HRBF,Macedo_HRBF,Carr_3D,Turk_2002,Wendland:2004,RBF_Med}
partition of unity methods~\cite{Adapt_PUM,Multi_PUM,Tobor_PUM,Ohtake_Sparse}, and methods which turn the reconstruction problem into a Poisson or screened Poisson problem~\cite{Pois_SR,Screen_Pois}. While a comprehensive review of the aforementioned methods is beyond the scope of this paper, we will briefly compare and contrast the methods related to CFPU; the interested reader is directed to the survey papers~\cite{Survey1,Survey3} for a more comprehensive review.

The first closely related methods express the surface reconstruction problem as the solution to a (screened) Poisson equation~\cite{Pois_SR,Screen_Pois} (the so called indicator function approach~\cite{Survey3}).  Similar to CFPU, these methods rely on the fact that the normals of an oriented point cloud are the gradient of a potential. They use the divergence of the normals to get a (screened) Poisson equation that can be solved directly for the potential $\tcfpot$.  A weak formulation of the equation is used which replaces the divergence of the normals with integrals in the surfaces ambient space that involve an extension of the normals from the surface.  One issue with these methods \refone{is} that they require solving a global set of equations to recover $\tcfpot$. In contrast, the CFPU method does not require extending the normals from the surface or computing their divergence since it reconstructs the potentials locally using curl-free RBF approximants of the normals directly from the surface.  Furthermore, the method produces the global potential $\tcfpot$ by blending the local potentials together, which makes the method more amenable to parallel implementations.  

A second closely related method to CFPU is HRBF Implicits~\cite{Macedo_HRBF}. The idea behind this method is to interpolate the potential (which is taken to be zero) at the point cloud samples  \emph{and} the normal vectors at these samples using a Hermite RBF interpolant. The surface $\mathcal{P}$ is then determined directly from this interpolant. Our method, by comparison, only interpolates the normal vectors using a specially constructed matrix-valued kernel that allows us to extract a potential for the field directly, which can be used to approximate $\mathcal{P}$. This allows us to immediately reduce the size of the linear systems that need to be solved by 33\% for 2D problems and 25\% for 3D problems. 

There are other surface reconstruction methods from the literature based on RBFs that are less closely related to the CFPU method, but are still relevant.  Global RBF methods were initially used for modeling surfaces due to their ability to handle sparse point clouds; however, their global nature restricted their applications to small problems~\cite{RBF_Med,Turk_99}. 
%~\cite{RBF_Med,Savchenko_95,Turk_99}. 
Subsequently, RBF-based methods have been developed which address this issue. Carr et al.\ introduced a reconstruction method which requires the addition of ``auxiliary'' points to the data in an $\ep$-width narrow band around the surface determined by (possibly approximated) normal vectors to the surface~\cite{Carr_3D}. This method can be sensitive to the selection of the $\ep$ parameter, introducing numerical instabilities into the reconstruction, especially around thin features~\cite{Closed_Form_HRBF}, and optimal choices for $\ep$ are lacking. While direct computation of this method can be expensive and requires solving a $3N$-by-$3N$ linear system, using fast summation algorithms~\cite{Carr_3D}, partition of unity~\cite{MeshFree_Mat}, and compactly supported RBFs~\cite{Ohtake_3D} have been shown to overcome this issue. The use of compactly supported RBFs has especially gained much attention, due to the resulting sparse linear systems; however, one must still choose a support radius for the compactly supported kernels, and if this value is too small, then the approximation power of the method can be impacted~\cite{Morse_Chen,Walder_2006,Ohtake_3D,Closed_Form_HRBF}. CFPU addresses the issue of computational complexity with RBF systems while still remaining numerically stable and accurate.  Additionally, the method does not require the choice of shape parameter or support radii for the RBFs.

\subsection{Contributions}
In this paper, we present a novel method for reconstructing curves and surfaces with curl-free, vector-valued RBFs. CFPU is fast and requires only points on the surface and their corresponding normals. The RBFs we use are free of shape and support parameters, and the resulting linear systems are well-conditioned. Additionally, our method can be made interpolatory or can be regularized to handle noisy data. Since the implementation of the method involves local partition of unity patches, it is also scalable and highly parallelizable.

%%%%%%%%%%%%%%%%%%%%%%%%%%%%%%%%%%%%%%%%%%%%%%%%%%%%%%%

%%%%%%%%%%%%%%%%%%%%%%%%%%%%%%%%%%%%%%%%%%%%%%%%%%%%%%%%%%%%
%%%%%%%%%%%%%%%%%%%%%%%%%%%%%%%%%%%%%%%%%%%%%%%%%%%%%%%%%%%%%
\section{Curl-free RBFs}
Curl-free RBFs were developed for the interpolation of curl-free vector fields that are given from scattered measurements as occurs, for example, in the areas of electrostatics and geodesy~\cite{AmodeiBenbourhim_1991,DoduRabut_2004_DFCFInterp,Fuselier08}.  This technique has the important features that the vector interpolants are analytically curl-free and are well-posed for scattered data.  Additionally, a scalar potential can be extracted directly from the interpolants that approximates the underlying potential of the field (up to an additive constant)~\cite{FuselierWright2017}.  Curl-free RBF interpolation is similar to scalar RBF interpolation in the sense that one constructs the interpolants from linear combinations of shifts of a kernel at each of the given data sites.  The difference between the approaches is that in the curl-free case one uses shifts of a \emph{matrix-valued kernel} whose columns are curl-free. For the sake of brevity, we do not review scalar RBF approximations but refer the reader to any of the books~\cite{MeshFree_Mat,Wendland:2004,FFBook}.

Let $\phi:\R^d\times\R^d\longrightarrow\R$ be a radial kernel in the sense that $\phi(\vx,\vy)=\eta(\|\vx-\vy\|)$, for some $\eta:[0,\infty)\longrightarrow\R$, where $\|\cdot\|$ is the vector two-norm, and $\vx,\vy\in\R^d$. It is common in this case to simply write $\phi(\vx,\vy) = \phi(\|\vx-\vy\|)$ and refer to $\phi$ as an RBF.  A matrix-valued curl-free kernel $\Phicurl$ is given as~\cite{Fuselier08}
\begin{align}
\Phicurl(\vx,\vy) = -\nabla\nabla^T \phi(\|\vx-\vy\|),
\label{eq:cf_rbf}
\end{align}
where $\nabla$ is the gradient in $\R^d$ applied to $\vx$, and $\phi$ is assumed to have two continuous derivatives.  Since $\Phicurl$ is built from an RBF, these kernels are simply called curl-free RBFs.  For any $\vc\in\R^d$ and fixed $\vy$, the vector field $\Phicurl(\vx,\vy)\vc$ is curl-free in $\vx$.  This follows since
\begin{align}
\Phicurl(\vx,\vy)\vc = \nabla\underbrace{\left(-\nabla^T\phi\left(\|\vx-\vy\|\right)\vc\right)}_{\ds g(\vx)} = \nabla(g(\vx)),
\label{eq:PhiCurlPot}
\end{align}
i.e.\  $\Phicurl(\vx,\vy)\vc$ is the gradient of a scalar function $g$.  Note that the second argument of $\Phicurl$ acts as a shift of the kernel and indicates where the field \eqref{eq:PhiCurlPot} is ``centered.''

An interpolant to a curl-free vector field $\vu\in\R^d$ sampled at distinct points $\Pts=\{\pts_j\}_{j=1}^N$ can be constructed from $\Phicurl$ as follows:
\begin{equation}
	\vs(\vx)=\sum_{j=1}^N\Phicurl(\vx,\pts_j)\vc_j,
	\label{eq:CurlInterp}
\end{equation}
where the expansion coefficients $\vc_j \in \R^d$ are found by enforcing the interpolation conditions $\vs\bigr|_\Pts = \vu\bigr|_\Pts$. This results in the linear system
\begin{align}
\sum_{j=1}^N\Phicurl(\pts_i,\pts_j)\vc_j = \vu_i,\; i=1,2,\ldots,N, 
\label{eq:lin_sys}
\end{align}
which is commonly written as $A\vc=\vu$, where $A$ is the block $dN$-by-$dN$ interpolation matrix
\begin{equation}
    A = \begin{bmatrix}
        \Phicurl (\pts_1,\pts_1)  &\Phicurl (\pts_1,\pts_2) &\cdots  &\Phicurl (\pts_1,\pts_N)\\
        \Phicurl (\pts_2,\pts_1)  &\Phicurl(\pts_2,\pts_2) &\cdots &\Phicurl (\pts_2,\pts_N)\\
        \vdots  &\vdots  &\ddots &\vdots\\
        \Phicurl (\pts_N,\pts_1)  &\Phicurl (\pts_N,\pts_2) &\cdots  &\Phicurl (\pts_N,\pts_N) 
    \end{bmatrix}.
    \label{AY}
\end{equation}
One can show that $A$ is positive definite if $\Phicurl$ is constructed from an appropriately chosen scalar-valued $\phi$~\cite{Fuselier08}; see Table~\ref{tbl:rbfs} for some examples.
\begin{table}[h!]
	\centering
	\begin{tabular}{|c|c|}
		\hline
		\textbf{Radial kernel}          & \textbf{Expression} \\ \hline
		Gaussian (GA)                   & $\phi(r) = \exp(-(\ep r)^2)$ \\ \hline
		Inverse multiquadric (IMQ) & $\phi(r) = (1 + (\ep r)^2)^{-\frac12}$ \\ \hline
		Multiquadric (MQ)               & $\phi(r) = -(1 + (\ep r)^2)^{\frac12}$ \\ \hline
		%Wendland (CSP)		 & $\phi(r) = (1-r/\ep)_{+}^{2\ell}q(r/\ep)$ \\ \hline
	\end{tabular}
	\caption{\label{tbl:rbfs} Examples of radial kernels that result in positive definite matrices $A$ \eqref{AY} for curl-free RBF interpolation.  Here $\ep>0$ is the shape parameter.}
\end{table}
An important feature from the construction of the curl-free RBF interpolant is that we can extract a scalar potential for the interpolated field by exploiting \eqref{eq:PhiCurlPot}: 
\begin{align}
\vs(\vx)= \nabla\biggl(\underbrace{-\sum_{j=1}^N\nabla^T\phi\left(\|\vx-\pts_j\|\right)\vc_j}_{\ds \cfpot(\vx)}\biggr).
\label{eq:CurlInterp_Rd}
\end{align}

While the method described above will ensure a curl-free interpolant of the sampled field, some issues do arise.  First, the size of the linear system \eqref{eq:lin_sys} grows rapidly with $N$, and, for a globally supported kernel, will be dense and computationally expensive to solve---requiring $\bigO((dN)^3)$ operations if a direct method is used.  Second, each evaluation of the interpolant (or potential \eqref{eq:CurlInterp_Rd}) involves $dN$ terms, which can become computationally expensive when many evaluations are necessary (such as occurs in the present application).  Other issues involve the shape parameter $\ep$ used in the radial kernels from Table \ref{tbl:rbfs}. This parameter controls how flat or peaked the radial kernels are and has a dramatic effect on both the accuracy of the interpolant as well as the conditioning of the interpolation matrix $A$. If $\ep$ is fixed and the total number of interpolation points $N$ grows, then the $A$ matrix becomes exponentially ill-conditioned with $N$.  Additionally, while extensive literature dedicated to finding the ``good'' values of $\ep$ to use exists for scalar RBF interpolation~\cite{Opti1,Opti2}, these approaches are computationally expensive, and this will only be exacerbated by the larger sizes of the linear systems for curl-free RBFs. To bypass these issues, we next discuss curl-free RBFs that do not feature a shape parameter.  We address the issues with the computational cost in Section \ref{sec:PUM}.
 
%%%%%%%%%%%%%%%%%%%%%%%%%%%%%%%%%%%%%%%%%%%%%%%%%%%%%%%%%%%%
%%%%%%%%%%%%%%%%%%%%%%%%%%%%%%%%%%%%%%%%%%%%%%%%%%%%%%%%%%%%%
\subsection{Curl-free polyharmonic splines}\label{sec:CFPHS}
Interpolants based on polyharmonic splines (PHS) were introduced by  Duchon as a generalization of univariate spline interpolants to higher dimensions~\cite{Du77}.  PHS interpolants have the property that they minimize an energy functional that can be interpreted as a type of ``bending energy'' for the surfaces they produce, similar to univariate splines~\cite[ch.\ 13]{Wendland:2004}.  PHS are radial kernels and come in the following two types:
\begin{align}
\phi_{\ell}(r) = 
(-1)^{\ell+1}
\begin{cases}
r^{2\ell} \log r,  & \ell\; \text{positive integer}, \\
r^{2\ell+1}, & \ell\; \text{non-negative integer}.
\end{cases}
\label{eq:phs_rbf}
\end{align}
For interpolation in $\mathbb{R}^d$, the first option is typically used for $d$ even and the second option for $d$ odd.  The choice for the order parameter $\ell$ is often made based on smoothness assumptions of the data, with larger $\ell$ for smoother data.  However, larger $\ell$ also negatively effects the numerical stability of the interpolants~\cite[ch.\ 12]{Wendland:2004}. The combination of $d$ and $\ell$ determines the minimization properties of the interpolants~\cite[ch.\ 13]{Wendland:2004}; the choice of $\ell=1$ for $d=2$ leads to the classical thin-plate spline.  PHS do not feature a shape parameter like other RBFs, as any scaling of $r$ just factors out of the kernels.  While $\ell$ is a free parameter, one does not need to continually search for a good value to use when the interpolation problem is changed, as is  typically the case for RBFs with shape parameters.  

Curl-free PHS were introduced in~\cite{AmodeiBenbourhim_1991} and have further been studied in~\cite{DoduRabut_2002_VectorialInterp,Poly_Vec}.  These matrix-valued kernels can be produced by using \eqref{eq:cf_rbf} with $\phi_{\ell}$ given by either of the choices in \eqref{eq:phs_rbf} and $\ell$ chosen large enough to ensure the derivatives make sense.  We will denote these kernels by $\Phicurl_{\ell}$. As with scalar PHS, it is necessary to modify the curl-free RBF interpolant \eqref{eq:CurlInterp} to ensure a well-posed problem.  In the curl-free PHS case, the interpolant \eqref{eq:CurlInterp} must be augmented with curl-free (vector) polynomials in $\mathbb{R}^d$ of degree $\ell-1$~\cite{DoduRabut_2002_VectorialInterp}, where degree refers to the total degree of any of the components of the (vector) polynomial.  A basis for curl-free polynomials up to degree $\ell-1$ can be generated as follows. Let $\{p_0,\ldots,p_L\}$ be a monomial basis for scalar polynomials up to degree $\ell$ in $\R^d$, where $L=\binom{\ell+d}{d}-1$ and $p_0 = 1$. Then a basis for curl-free polynomials up to degree $\ell-1$, $\{\vp_1,\ldots,\vp_L\}$, is given by applying the gradient to each $p_i$, i.e., $\vp_i = \nabla p_i$, $i=1,\ldots,L$.  \refone{Note that the gradient reduces the degree of the scalar polynomials by one, which is why we end up with curl-free polynomials of degree $\ell-1$.} As an example, we give a basis for degree 1 curl-free polynomials in $\mathbb{R}^2$ and $\mathbb{R}^3$ \refone{obtained from a (scaled) basis of monomials of degree $\ell=2$}:
\begin{align*}
\text{Poly.\ basis in $\mathbb{R}^2$}:& \quad
\left\{\begin{bmatrix}1\\
0
\end{bmatrix},\begin{bmatrix}0\\
1
\end{bmatrix},\begin{bmatrix}y\\
x
\end{bmatrix},\begin{bmatrix}x\\
0
\end{bmatrix},\begin{bmatrix}0\\
y
\end{bmatrix}\right\} \\
\text{Poly.\ basis in $\mathbb{R}^3$}:& \quad
\left\{\begin{bmatrix}1\\
0\\
0
\end{bmatrix},\begin{bmatrix}0\\
1\\
0
\end{bmatrix},\begin{bmatrix}0\\
0\\
1
\end{bmatrix},\begin{bmatrix}y\\
x\\
0
\end{bmatrix},\begin{bmatrix}z\\
0\\
x
\end{bmatrix},\begin{bmatrix}0\\
z\\
y
\end{bmatrix},\begin{bmatrix}x\\
0\\
0
\end{bmatrix},\begin{bmatrix}0\\
y\\
0
\end{bmatrix},\begin{bmatrix}0\\
0\\
z
\end{bmatrix}\right\}.
\end{align*}

%While the curl-free polyharmonic spline kernel does not lead to a positive definite interpolation matrix, a minor adjustment to the interpolant can ensure a conditionally positive definite system. To demonstrate this, we introduce the concept of a curl-free polynomial basis. Let $\{p_k(\mathbf{x})\}_{k=1}^{L+1}$ be a basis for scalar polynomials up to degree $\ell$ in $\R^d$ at $\mathbf{x}$, where $L=\binom{\ell+d}{d}-1$. One can generate the curl-free polynomial basis vectors, $\vp_k(\mathbf{x})$ by applying the gradient to $p_k(\mathbf{x})$. 

A curl-free PHS interpolant of order $\ell$ to a curl-free vector field $\vu\in\R^d$ sampled at distinct points $\Pts=\{\pts_j\}_{j=1}^N$ is given as follows:
\begin{equation}
\vs(\vx) = \sum_{j=1}^N  \Phicurl_{\ell} (\vx,\pts_j)\vc_j+\sum_{k=1}^L  b_k \vp_k \left(\vx\right),
\label{CF_PHS_interpolant}
\end{equation}
where the interpolation coefficients $\vc_j\in\mathbb{R}^d$ and $b_k\in\mathbb{R}$ are determined by the conditions $\vs\bigr|_\Pts = \vu\bigr|_\Pts$ as well as the constraints 
 \begin{equation*}
\sum_{j=1}^N  \vc_j^T \vp_k(\pts_j)=0, \qquad k=1,2,\dots,L.
\end{equation*}
These constraints are necessary for the interpolant to minimize a certain energy norm and also limit its far-field growth~\cite{DoduRabut_2002_VectorialInterp}. We can state both constraints in terms of the following linear system of equations:
\begin{equation}
\label{rbf_phs_ls}
    \begin{bmatrix}
        A &P\\
        P^T &\mathbf{0}
    \end{bmatrix}
    \begin{bmatrix}
        \mathbf{c}\\
        \mathbf{b}
    \end{bmatrix}
    =
    \begin{bmatrix}
        \mathbf{u}\\
        \mathbf{0}
    \end{bmatrix},
\end{equation}
where $A$ is defined in~\eqref{AY} and
\begin{align*}
    P&=
    \begin{bmatrix}
        \vp_1(\pts_1)  &\vp_2(\pts_1)  &\cdots  &\vp_L(\pts_1) \\
        \vp_1(\pts_2)  &\vp_2(\pts_2)  &\cdots &\vp_L(\pts_2) \\
        \vdots  &\vdots  &\ddots &\vdots\\
        \vp_1(\pts_N)  &\vp_2(\pts_N)  &\cdots &\vp_L(\pts_N) \\
    \end{bmatrix}.
\end{align*}
Provided the set of points $\Pts$ is unisolvent with respect to the curl-free polynomial basis (i.e.\ $P$ is full rank), this linear system is non-singular and thus the interpolation problem is well-posed~\cite{DoduRabut_2002_VectorialInterp}.  However, as with interpolation matrices based on curl-free RBFs with shape parameters, this interpolation matrix also becomes ill-conditioned as $N$ increases, but the growth is algebraic as opposed to exponential~\cite{Fuselier08}.

We note that a scalar potential $\cfpot$ can also be recovered similarly to~\eqref{eq:CurlInterp_Rd} as
\begin{align}
\vs(\vx)= \nabla\biggl(\underbrace{-\sum_{j=1}^N\nabla^T\phi_{\ell}\left(\|\vx-\pts_j\|\right)\vc_j+\sum_{k=1}^L  b_k p_k \left(\mathbf{x}\right)}_{\ds \cfpot(\vx)}\biggr),
\label{eq:PHSCurlInterp_Rd}
\end{align}
where $p_k$ are the scalar polynomial basis used to generate the curl-free polynomial basis. This potential plays a key role in the CFPU method.

\subsection{An example} \label{sec:globalCF}
\label{R2recon}
\begin{figure}[h!]
\centering
\begin{tabular}{cc}
\includegraphics[width=0.45\textwidth]{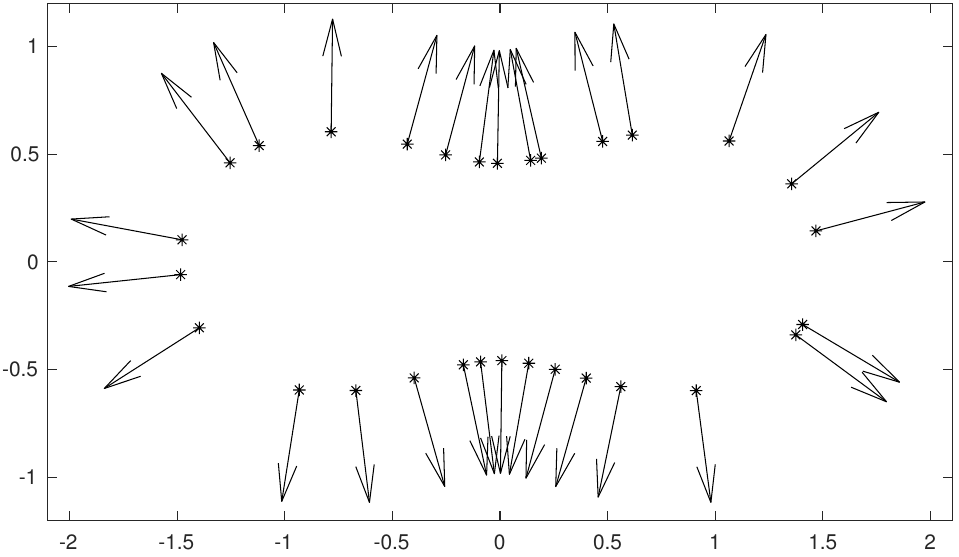} & 
\includegraphics[width=0.45\textwidth]{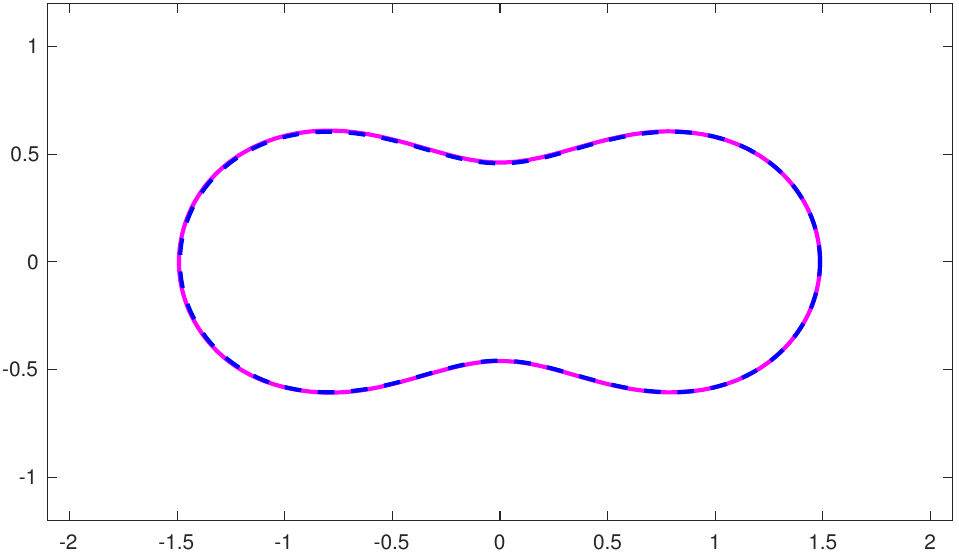} \\
(a) & (b) 
\end{tabular}
\caption{(a) $N=30$ points sampled from a Cassini oval \eqref{eq:cassini} with $a=1$ and $b=1.1$, together with the corresponding normal vectors to the curve. (b) Zero-level curve from the global curl-free PHS interpolation method with $\ell=2$ (magenta) together with the exact curve $\mathcal{P}$ (blue dashed line).}
		\label{cassini}
\end{figure}

At this point, it is illustrative to see how curl-free RBFs can be used to recover a level set from an oriented point cloud.  We focus here on the case of a level curve $\mathcal{P}\subset \R^2$, since this is the situation in which the global method described above would be applicable due to the smaller problem sizes.  The example we focus on uses the Cassini oval as the target level curve to recover, which can be described as the zero-level set of the implicit function
\begin{equation}
f(\vx) = f(x_1,x_2)=(x_1^2+x_2^2)^2-2a^2(x_1^2-x_2^2)+a^4-b^4,
\label{eq:cassini}
\end{equation} 
where we take $a=1$ and $b=1.1$; see the Figure~\ref{cassini}(b) for the resulting curve (dashed-line). We sample this curve at $N$ nodes $\Pts=\{\pts_j\}_{j=1}^N$ and compute the corresponding (unit) normal vectors using $\vn_j=\nabla f(\pts_j)$; see Figure~\ref{cassini}(a) for a plot of the exact data used for the case of $N=30$.
As mentioned in the introduction, the key to our method is the fact that the normal vectors to a level curve (or surface) are curl-free.  We thus fit the data $(\pts_j,\vn_j)$, $j=1,\ldots,N$, using the curl-free RBF interpolant~\eqref{CF_PHS_interpolant} and from this extract out the potential $\cfpot$ as in~\eqref{eq:PHSCurlInterp_Rd}. Since the potential for a curl-free field is only unique up to a constant (which is a consequence of the Helmholtz decomposition theorem~\cite{FuselierWright2017}), the zero-level curve of $\cfpot$ will not necessarily approximate the zero-level curve of $f$.  To fix this we set $\tilde{\cfpot}(\vx) = \cfpot(\vx)-\mu$, where $\mu$ is the discrete mean of $\cfpot$ at the nodes $\Pts$.  The result from this experiment is shown in Figure~\ref{cassini}(b), where we see excellent agreement between the zero-level curve of $\tilde{\cfpot}$ and $f$.  

While this global method is reasonable to use for reconstructing level curves and surfaces when $N$ is small, the cost of solving the linear system \eqref{rbf_phs_ls} becomes too expensive for large $N$, which will be the case for any sampling of a complex surface. We address this issue next.

%%%%%%%%%%%%%%%%%%%%%%%%%%%%%%%%%%%%%%%%%%%%%%%%%%%%%%%%%%%%
%%%%%%%%%%%%%%%%%%%%%%%%%%%%%%%%%%%%%%%%%%%%%%%%%%%%%%%%%%%%%
\section{The CFPU method}
\label{sec:PUM}
Partition of unity (PU) methods offer a way to split up a global approximation problem on a domain $\Omega$ into local approximation problems on overlapping patches covering $\Omega$.  These local approximations are then blended together to form a global approximant using weight functions that form a partition of unity~\cite{Babuska_PU}.  This procedure can drastically reduce the computational cost of the original approximation problem. In order to explain the CFPU method, we first give a brief description of the ideas behind PU methods as they pertain to our problem, and introduce some necessary notation for what follows.

\subsection{PU methods}\label{sec:pum}
Let $\Omega\subset\R^d$ be an open, bounded domain of interest for approximating some function $f:\Omega\longrightarrow\R$. Let $\Omega_1,\dots,\Omega_M$ be a collection of distinct overlapping patches that form an open cover of $\Omega$, i.e., $\cup_{m=1}^M \Omega_m \supseteq \Omega$, and let the overlap between patches be limited such that at most \refone{$K \ll M$} patches overlap at any given point $\vx\in\Omega$.  For each $m=1,\ldots,M$, let $w_m:\Omega_m\longrightarrow[0,1]$ be a weight function such that $w_m$ is compactly supported on $\Omega_m$, and let the set of weight functions $\{w_m\}$ have the property that 
$\sum_{m=1}^{M} w_m \equiv 1$.
Suppose $s_m$ is some approximation to $f$ on each patch $\Omega_m$.  Then the PU approach of Babu\v{s}ka and Melenk~\cite{Babuska_PU} is to form an approximant $s$ to $f$ over the whole domain $\Omega$ by ``blending'' the local approximants $s_m$ with $w_m$ as follows:
$s = \sum_{m=1}^M w_m s_m$.  

Balls are common choices for the patches, since then the weight functions $w_m$ can be easily constructed using Shepard's method~\cite{Shepard1968} as follows.  Let $\kappa:\R^{+}\rightarrow\R$ have compact support over the interval $[0,1)$.  For each patch $\Omega_m$, let $\vom_m\in\R^d$ denote its center and $\prad_m>0$ denote its radius, and define $\kappa_m(\vx) := \kappa\left(\|\vx-\vom_m\|/\prad_m\right)$.  The weight function for patch $\Omega_m$ is then given by
\begin{equation}
w_m(\vx) = \kappa_m(\vx)/\sum_{j=1}^{M} \kappa_j(\vx),\; m=1,\ldots,M.
\label{eq:wght_func}
\end{equation}
Note that each $w_m$ is only supported over $\Omega_m$ and that the summation in the denominator only involves terms that are non-zero over patch $\Omega_m$, which is bounded by $K$.  In this study, we choose the patches as balls, and use the $C^1$ quadratic $B$-spline
%Figure~\ref{fig:RBF_PU} (b) illustrates one of these weights functions for the example domain in part (a), where $\kappa$ is chosen as the $C^1$ quadratic $B$-spline
\begin{align}
\kappa(r) = \begin{cases} 1 - 3r^2, & 0\leq r \leq \frac{1}{3}, \\ \frac{3}{2}(1 - r)^2, & \frac{1}{3}\leq r \leq 1 \end{cases}
\end{align}
to define the weight functions. \reftwo{With this weight function the global approximant $s$ is at most $C^1$.  However, weight functions with more smoothness can be easily created from smoother compactly supported univariate functions $\kappa$ if global approximants with higher smoothness are desired.}  

RBFs are commonly used with the PU approach to reduce the computational cost.  They have been used for approximating a function from scattered samples (e.g.~\cite{wendland_2002,cavoretto_2010,CaDeRoPe:2015}) and solving differential equations (e.g.~\cite{Safdari-Vaighani2015,Larsson2017,shankar_2018}).  Recently the present authors presented a PU method for interpolation of divergence-free and curl-free vector fields~\cite{DFW2020a}, which is what the current approach is based on.  \reftwo{For RBF-PU methods, it should be noted that the smoothness of the weight functions \eqref{eq:wght_func} does not govern the rates of convergence of the global PU approximants~\cite{Larsson2017,DFW2020a}.}

%\begin{figure}[h]
%	\centering
%	\begin{tabular}{cc}
%	\includegraphics[width=0.4\textwidth]{disk_patches.pdf} & \includegraphics[width=0.45\textwidth]{pu_patach_disk.pdf} \\
%	(a) & (b)
%	\end{tabular}
%\caption{(a) Illustration of partition of unity patches (outlined in blue lines) for a node set $\Pts$ (marked with black dots) contained  in a domain $\Omega$. (b) Illustration of one of the PU weight functions for the patches from part (a), where the color transition from white to yellow to red to black correspond to weight function values from $0$ to $1$.}
%\label{fig:RBF_PU}
%\end{figure}

\subsection{Description of CFPU}\label{sec:pu_description}
For brevity, we describe the CFPU method for reconstructing a zero-level surface $\mathcal{P}$ in $\mathbb{R}^3$ defined by $f(\vx) = 0$ using curl-free PHS of order $\ell$.  Let $\Pts=\{\pts_j\}_{j=1}^N$ be a given set points on $\mathcal{P}$ and let $\{\vn_j\}_{j=1}^N$ denote the unit normals (or approximations to the normals) of $\mathcal{P}$ at $\Pts$.  Let $\Omega_1,\dots,\Omega_M$ be a set of overlapping patches that form an open cover of $\mathcal{P}$.  Finally, let $\Pts_m$ denote the nodes contained in $\Omega_m$ and $n_{m}$ denote the cardinality of $\Pts_m$.  For each $\Omega_m$, we fit a curl-free RBF interpolant $\vs_m$  of the form~\eqref{CF_PHS_interpolant} to the normals at the points in $\Pts_m$ and then extract from this its scalar potential $\cfpot_m$ using~\eqref{eq:PHSCurlInterp_Rd}.  

A natural first approach to constructing a global potential for approximating $\mathcal{P}$ would be to blend these local potentials $\cfpot_m$ using the PU weight functions $w_m$ as $\cfpot = \sum_{m=1}^M w_m \cfpot_m$.  However, this will lead to an issue since each $\cfpot_m$ is only unique up to a constant.  This means that for two patches $\Omega_{k}$ and $\Omega_m$ that overlap, $\cfpot_{k}$ and $\cfpot_m$ may be shifted from one another in the overlap region, which would then lead to an inaccurate global PU approximant in the overlap.  To fix this issue we can shift each potential by a different constant so that they approximately agree in the overlap region.  To determine these constants, we use the fact that the points $\Pts_m$ reside on the zero-level surface $\mathcal{P}$, and we want each potential $\cfpot_{m}$ to be approximately zero on $\Pts_m$.  One way to achieve this result is to enforce that the discrete mean of the local potentials is zero over the patch nodes.  To this end, let \revision{$\mu_m$ equal the mean of the values of $\cfpot_m$ at $\Pts_m$}, and then define the shifted potential 
\begin{align*}
\mcfpot_m := \cfpot_m-\mu_m.
\end{align*}
The global CFPU approximant for the underlying implicit function $f$ is given by
\begin{align}
\cfpot(\vx):=\sum_{m=1}^M w_m(\vx)\mcfpot_m(\vx).
\label{eq:curlfree_mcf_pum}
\end{align}
We can then approximate $\mathcal{P}$ as the surface defined by the set of all $\vx$ in $\cup_{m=1}^M \Omega_m$ such that $\cfpot(\vx) = 0$.

The CFPU method requires solving $M$ linear systems of size $(3n_m + L)$-by-$(3n_m+L)$ rather than one large $(3N+L)$-by-$(3N+L)$  system for the global curl-free method described in Section \ref{sec:globalCF}.  We select \refone{$n_m \ll N$}, for all $m$ so that the computational cost is significantly reduced, i.e.\ the complexity is $\bigO(\sum_{m=1}^M (3n_m + L)^3 N)$ rather than $\bigO( (3N + L)^3)$, when using a direct solver.  Furthermore, each of these smaller systems can be solved independently, making the CFPU method pleasingly parallel compared to the global method.  Finally we note that the computational complexity for evaluating the CFPU approximant \eqref{eq:curlfree_mcf_pum} is also significantly less than the global method.  For each evaluation point, only a small subset (equal to the number of patches that contain the evaluation point) of the local potentials $\mcfpot_m$ need to be evaluated.  The cost of evaluating each $\mcfpot_m$ is $\bigO(3n_m + L)$ rather than $O(3N + L)$ for the global method. These potentials can also be evaluated independently.

%%%%%%%%%%%%%%%%%%%%%%%%%%%%%%%%%%%%%%%%%%%%%%%%%%%%%%%%%%%%
%%%%%%%%%%%%%%%%%%%%%%%%%%%%%%%%%%%%%%%%%%%%%%%%%%%%%%%%%%%%%
\subsection{Exact interpolation}\label{sec:exact_interp}
The CFPU method as described above will not in general exactly interpolate the zero-level surface $\mathcal{P}$ at the points in $\Pts$, i.e.\ $\cfpot(\pts_j) \neq 0$, $j=1,\ldots,N$, which is often desirable when the points are assumed to be exactly on the surface~\cite{Macedo_HRBF}.  We can, however, enforce this condition by simply subtracting an interpolant of the residual from each patch potential $\cfpot_m$.  We describe the details of this procedure below.

Let the points in patch $\Omega_m$ from $\Pts$ be denoted by $\Pts_m = \{\pts_j^m\}_{j=1}^{n_m}$ and let $\revision{\sigma_m}$ be a scalar PHS interpolant to the values of $\cfpot_m$ at $\Pts_m$.  Using the same notation from Section \ref{sec:CFPHS}, this interpolant can be written as
\begin{align}
\sigma_m(\vx) = \sum_{j=1}^{n_m} c_j^m \phi_{\ell}(\|\vx - \pts_j^m\|) + \sum_{k=0}^{L} b_k^m p_k(\vx),
\label{eq:tps_interp}
\end{align}
where the coefficients are determined from the interpolation conditions $\sigma_m\bigr|_{\Pts_m} = \cfpot_m\bigr|_{\Pts_m}$ and the moment conditions $\sum_{j=1}^{n_m} c_j^m p_k(\pts_j^m) = 0$, $k=0,\ldots,L$.  The potential $\cfpot_m$ on patch $\Omega_m$ can be shifted by $\sigma_m$ to obtain
\begin{align*}
\icfpot_m := \cfpot_m-\sigma_m.
\end{align*}
By construction, this shifted potential satisfies $\icfpot_m(\pts_j^m) = 0$, so that an interpolatory global CFPU approximant for the underlying implicit function $f$ can then be obtained as
\begin{align}
\cfpot(\vx):=\sum_{m=1}^M w_m(\vx)\icfpot_m(\vx).
\label{eq:curlfree_icf_pum}
\end{align}
An approximation to the surface $\mathcal{P}$ is again given as the set of all $\vx$ in $\cup_{m=1}^M \Omega_m$ such that $\cfpot(\vx) = 0$.

We note that this exact interpolation technique is more expensive than just shifting the potentials by the mean as in \eqref{eq:curlfree_mcf_pum}, but only by a constant factor.  Additionally, just as with $\cfpot_m$, each $\sigma_m$ can be determined independently of the others. 

The residual of the patch potentials $\cfpot_m$ can be highly oscillatory, which could lead to spurious oscillations in the interpolants $\sigma_m$ and hence also $\cfpot$.  For 3D reconstructions, we thus recommend using the PHS kernel $\phi_0$ (i.e.\ $\phi_0(r) = r$) for the residual since the function will have minimal bending energy among all interpolants~\cite{Du77}.  This is the choice we use in all of our examples.

%
% the approximate curves or surfaces, but a modification must be made in order to enforce interpolation of the implicit function. An exact interpolant of the curve or surface means that the recovered potential on each patch evaluates to zero at the sample nodes on the patch, $\Pts_m$. We achieve this by evaluating the recovered potential at the nodes on the patch, interpolating this with a scalar polyharmonic RBF, and then subtracting this interpolant from the recovered potential. We then blend the combination of these two approximants with the partition of unity weight functions. In practice we have found that it is best to use a low order $k=1$ polyharmonic spline.

%%%%%%%%%%%%%%%%%%%%%%%%%%%%%%%%%%%%%%%%%%%%%%%%%%%%%%%%%%%%
%%%%%%%%%%%%%%%%%%%%%%%%%%%%%%%%%%%%%%%%%%%%%%%%%%%%%%%%%%%%%
\subsection{Regularization}\label{sec:regularization} There are two regularization techniques we can apply to deal with noise.
\subsubsection{Normals}\label{sec:nrml_regularization}
If the samples of the normal vectors of $\mathcal{P}$ are corrupted with noise, then interpolating them exactly on each patch to recover the potentials may cause issues, such as producing spurious sheets in the reconstructed surface.  In this case, it may make sense to instead introduce some regularization in the vector approximants on the patches.  Regularized kernel approximation, such as smoothing splines or ridge regression~\cite{splines_wahba}, offers one effective way to do this.

For a curl-free PHS approximant given in \eqref{CF_PHS_interpolant}, the smoothing spline regularization approach amounts to solving the following minimization problem:
\begin{equation}
\min_{\vc \in \R^{3n}} \left[\frac{1}{3n}\sum_{j=1}^n\|\vs(\pts_j)-\vn_j\|^2 + \lambda \vc^T A \vc\right],\;\text{subject to}\; P^T\vc = 0,
\label{eq:functional}
\end{equation}
where $A$ and $P$ are the matrices from \eqref{rbf_phs_ls} for $n$ points. The first term in the quadratic functional measures the goodness of fit of the approximant while the second term measures its smoothness.\footnote{This term arises from the minimization of a Hilbert space semi-norm related to the function space of the vector approximants---the so called native space norm associated with the $\Phicurl$ kernel~\cite{Poly_Vec}.}  The regularization parameter $\lambda\geq 0$ controls the tradeoff between these terms, with larger $\lambda$ resulting in smoother approximants. For a given $\lambda$, we can obtain the minimizer of the constrained quadratic functional \eqref{eq:functional} by solving the following modified version of the system \eqref{rbf_phs_ls}:
%that balances the fit of the approximant  and $J_m(f)$ is a smoothness penalty defined by the partial derivatives of total order $m$:
%\begin{equation*}
%J_m(f)=\int\limits_{-\infty}^{\infty} \int\limits_{-\infty}^{\infty} \sum_{i=0}^m \binom{m}{i} \left[\frac{\partial^m f(x,y)}{\partial x^i \partial y^{m-i}}\right]^2 dx dy.
%\end{equation*}
\begin{equation}
    \begin{bmatrix}
        A+3n\lambda I &P\\
        P^T &\mathbf{0}
    \end{bmatrix}
    \begin{bmatrix}
        \mathbf{c}\\
        \mathbf{b}
    \end{bmatrix}
    =
    \begin{bmatrix}
        \mathbf{u}\\
        \mathbf{0}
    \end{bmatrix},
    \label{eq:regsys}
\end{equation}
where $I$ is the $3n$-by-$3n$ identity matrix and $\vu$ contains the normals.  

For the CFPU method, we use this regularization approach on each patch $\Omega_m$ to obtain regularized potentials $\cfpot_m$.  This opens up the option of using a different regularization parameter $\lambda_m$ on each patch, and thus controlling the regularization of the approximants spatially.

In Section \ref{sec:results}, we present some results for selecting the regularization parameters $\lambda_m$ using generalized cross validation (GCV)~\cite{splines_wahba} and use the method from~\cite{GCV87} to compute these parameters.  The GCV approach is more computationally expensive than an ad hoc approach of choosing $\lambda_m$ and then solving \eqref{eq:regsys} directly, but it does provide a means for automating the choice of parameters.

\subsubsection{Residual}\label{sec:res_regularization}
If the point samples in the point cloud are also noisy, as often occurs in range scans of real 3D objects, then enforcing exact interpolation on the patch potentials by interpolating the residual may again cause issues in the reconstruction.  We can also introduce regularization in this process using smoothing splines.  In this case, one uses a similar minimization problem as \eqref{eq:functional}, but for the scalar TPS approximant \eqref{eq:tps_interp}.  In fact, smoothing splines were developed for this scalar problem~\cite{splines_wahba}.  Since different regularization parameters can be used for fitting the normals vs.\ fitting the residual of the potential, we let $\alpha$ denote the residual regularization parameter to avoid confusion.  We also experiment with picking this parameter using GCV in Section \ref{sec:results}.

\subsection{Choosing the PU patches}
As discussed above, each PU patch $\Omega_m$, $m=1,\ldots,M$, is chosen as a ball and can thus be determined by its center $\vom_m$ and radius $\rho_m$.  In this work, we choose the set of patch centers $\Xi = \{\vom_m\}_{m=1}^M$ directly from the point cloud using the \revision{point cloud simplification technique implemented in MeshLab~\cite{meshlab}, which is based on Poisson disk sampling ~\cite{poissondisk}.}  The input to the function is an estimate for the number of subsamples of the point cloud that should be produced.  The result is a quasi-uniform point cloud with respect to the (approximated) geodesic distance for the surface represented by the original point cloud.  We use this reduced point cloud for the patch centers.  The \revision{point cloud simplification} technique is computationally efficient, with an estimated complexity of $\mathcal{O}(M\log M)$. Another option for generating the patch centers is to use the weighted sample elimination method from~\cite{Yuksel2015}, which is also a Poisson disk sampling technique.  This method allows one to produce an exact number for the reduced point cloud, but is not as computationally efficient as the method in MeshLab.

\revision{The radii for the patches are determine using the following procedure.  For each patch center $\vom_m \in \Xi$, we determine the distance to its nearest neighbor (using the Euclidean distance) in the set $\Xi \setminus \vom_m$ using a $k$-d tree.  We then set $\tau$ as the maximum of these distances.  An initial radius $\rho$ for all the patches is then chosen as $\rho=(1+\delta)\tau/2$, where $\delta > 0$.}  The parameter $\delta$ controls the overlap of the patches and we have found that using $\delta=1$ gives consistently good results.  Setting $\rho_m = \rho$ initially for each patch $\Omega_m$, we then check if the number of samples, $n_m$, it contains satisfies $n_m \geq n_{\min}$.  If it does not, then we increase $\rho_m$ until this inequality is satisfied.  In this work, we choose $n_{\min} = 2L = 2(\binom{\ell+3}{3}-1)$, where $\ell$ is the order of the curl-free PHS kernel used.  It is necessary for the curl-free PHS interpolant to use $L$ points to be well-posed, and we have found that doubling this value works well in practice.  \reftwo{The final step involves adjusting the radii of the patches so that every sample point is included in at least one patch.  We do this by finding the closest patch center to each outlying point and increasing its radius so that point is enclosed.}

In our implementation, we build a $k$-d tree on the patch centers and use this for range queries in determining which sample points belong to which patches.  

\section{Results}\label{sec:results}
In this section, we test the CFPU method on several different 3D surface reconstruction problems.  We start with a comparison of the CFPU method with the global curl-free method for a small point cloud.  We then focus solely on the CFPU method.  We use point clouds generated from a known surface to test the accuracy of the method.  Using this same known surface, we then add noise to the normals and show how regularization can help reconstruct smooth surfaces with noisy normals.  Next, we consider a problem with raw range data that contains misalignments of the points and noise to show how both regularization of the normals and the residuals can result in a smooth reconstruction.  \revision{We then show how the method performs for reconstructing various common surfaces found in the literature.  Lastly, we present results on computational performance of the CFPU method.  For many of the examples, we also include comparisons to the surface reconstructions from two popular methods, the Screened Poisson Surface Reconstruction (SPR) method~\cite{Pois_SR,Screen_Pois} and the Smooth Signed Distance (SSD) method~\cite{SSR_2011}.  Unlike many other methods (e.g.\ HRBF Implicits), open-source versions of SPR and SSD are freely available~\cite{Kazhdan:CGF:2019}.}
Table \ref{tbl:parameters} lists the names of the surfaces we use in the experiments together with the size of the point cloud and the number patches $M$ used for the CFPU reconstructions.  The choices for $M$ were chosen in a somewhat ad hoc manner and were not optimized.
%A goal of these tests is to show how the method parameters, such as the order $\ell$ for the curl-free PHS kernels $\Phi_{\ell}$ and the regularization parameters $\lambda$ and $\alpha$ affect the reconstructions.  See Table \ref{tbl:parameters} for a listing of  various properties and parameters used in the tests.

\begin{table}[htb]
\centering
\begin{tabular}{|c||c|c|c|}
\hline
Surface & Num.\ points $N$ & Num.\ patches $M$ & Normals \\
\hline \hline
Homer  & 5103 & 135 & MeshLab \\
Trefoil  & [6114,32856]  & 864 & Exact \\
Stanford Bunny  & 362271 & 14472 & Raw \\
Stanford Dragon & 434856 & 14400 & Raw \\
Happy Buddha & 583079 & 42861 & MeshLab \\
Armadillo & 172974 & 14349 & MeshLab \\
Raptor & 135740 & 10337 & Raw \\
Filigree & 514300 & 35130 & Raw \\
Pump Carter & 238141 & 14351 & Raw \\
Dancing children & 724588 & 17296 & Raw \\
Gargoyle & 499548 & 21436 & Raw\\
\hline
\end{tabular}
\caption{Various properties and parameters used in the numerical examples. Entries in the Normals column indicate the technique for generating the normals: Exact=computed exactly from the surface; MeshLab=computed using the ``Compute normals for point set'' in MeshLab using 10 neighbors and no smoothing; Raw=normals are used directly from the values given in the downloaded files from either the Stanford University 3D Scanning Repository or the AIM@SHAPE-VISIONAIR repository.\label{tbl:parameters}}
\end{table}

While the problems we consider have normals that can be either computed analytically or that are supplied with the point cloud data, this may generally not be the case.  There are many algorithms available that can approximate the normals directly from a point cloud data (e.g.,~\cite{Jin_05,Feng_05,Klas_09}).  In some of the results that follow, we have experimented with the point cloud normal estimation function from MeshLab.  The last column of Table \ref{tbl:parameters} lists the surfaces for which we used this function.

In all of the tests, we use either exact interpolation or its regularized version to approximate the potentials $s$ (cf.\ \eqref{eq:curlfree_icf_pum}).  We found that both of these procedures result in better reconstructions of the zero-level surface than the mean-shift version (cf.\ \eqref{eq:curlfree_mcf_pum}). Additionally, we use the \refone{curl-free} PHS kernels generated from the $\phi_{\ell}(r) = (-1)^{\ell+1}r^{2\ell+1}$ since the point clouds are in $\mathbb{R}^3$.  As discussed previously, the order $\ell$ controls the smoothness of the kernels. 

All results were obtained from a MATLAB implementation of the CFPU method.  Additionally, we used the MATLAB isosurface function to extract the zero-level set for the surfaces based on a sufficiently dense sampling on a uniform 3D grid surrounding the point cloud.  The code is freely available from the CFPU GitHub repository~\cite{cfpusoftware}, including MEX interfaces to the MeshLab functions (VCGLib library) used determining the patch centers and normals.  All the CFPU results presented below can be reproduced from the code and data in the repository.

\subsection{Comparison with the global method}
\begin{figure}[htb]
\centering
\begin{tabular}{c}
Point cloud \\
\includegraphics[width=0.18\textwidth]{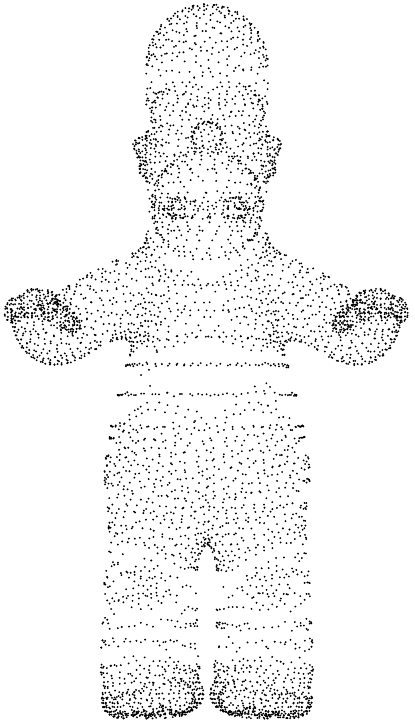}
\end{tabular}
\begin{tabular}{ccc}
\multicolumn{2}{c}{Global CF method} &  \revision{SPR} \\
$\ell=1$ & $\ell=2$ & \\ 
\includegraphics[width=0.18\textwidth]{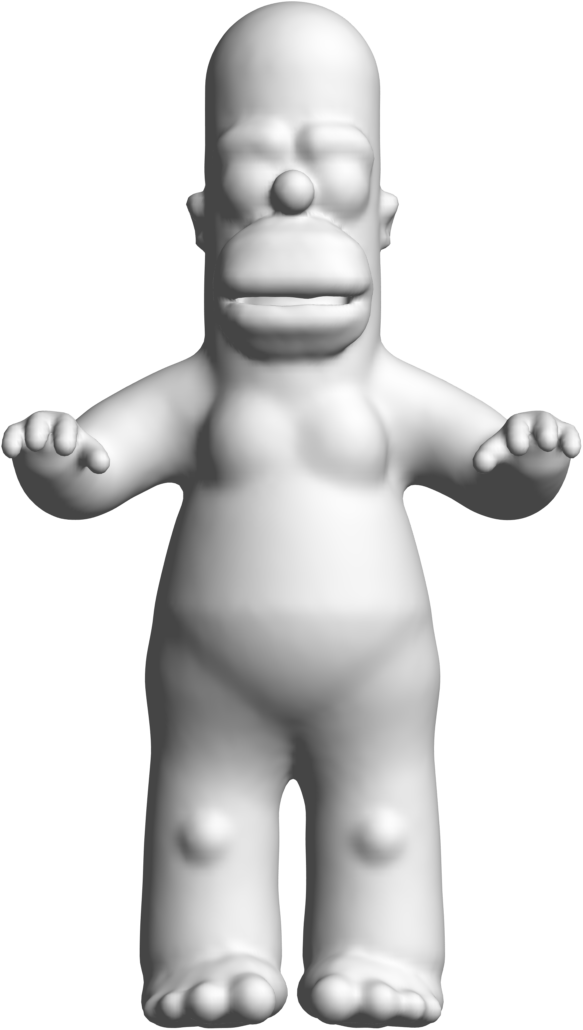} & 
\includegraphics[width=0.18\textwidth]{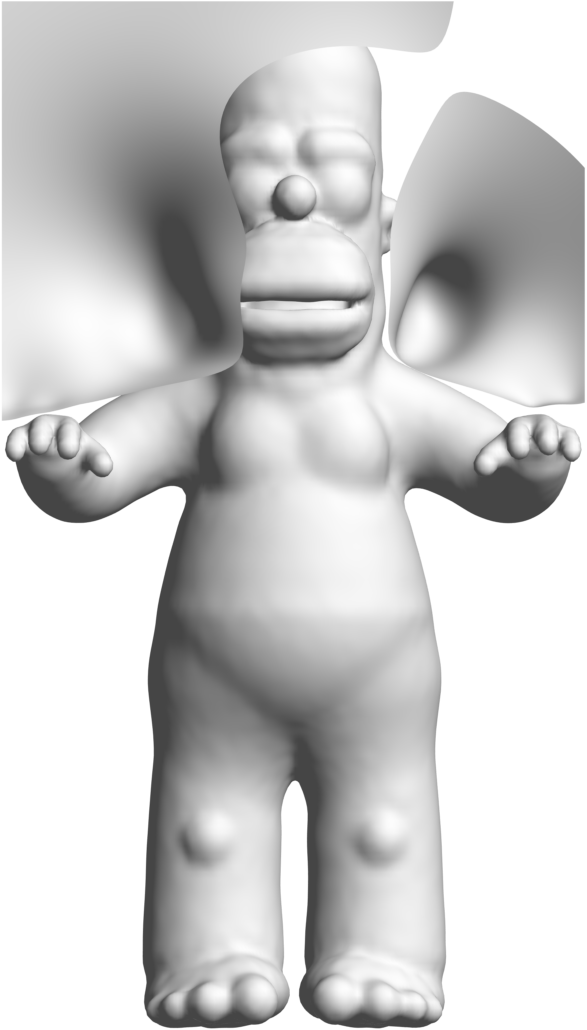} & 
\includegraphics[width=0.18\textwidth]{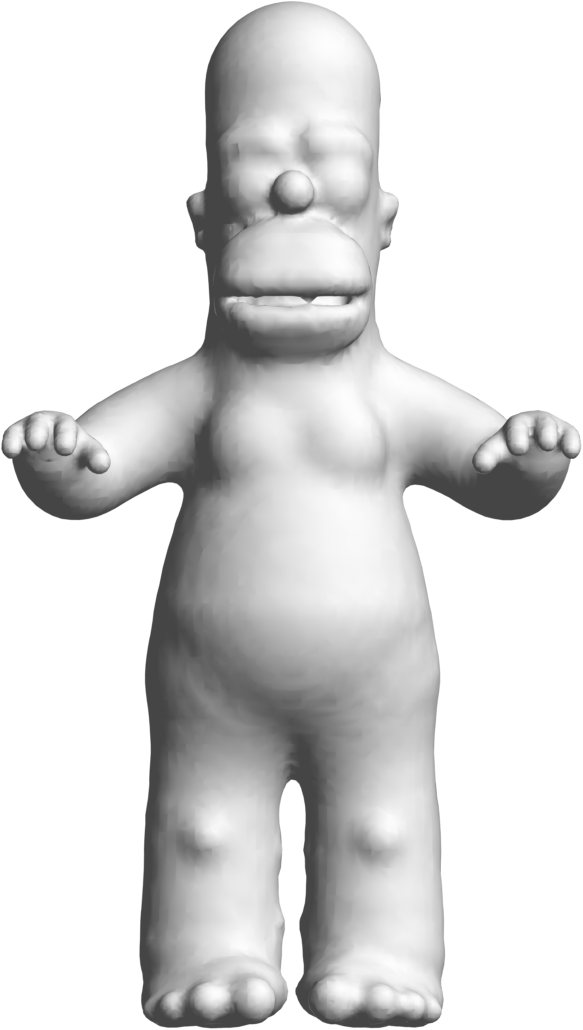} \\
\multicolumn{2}{c}{CFPU method}  & \revision{SSD} \\
$\ell=1$ & $\ell=2$  \\ 
\includegraphics[width=0.18\textwidth]{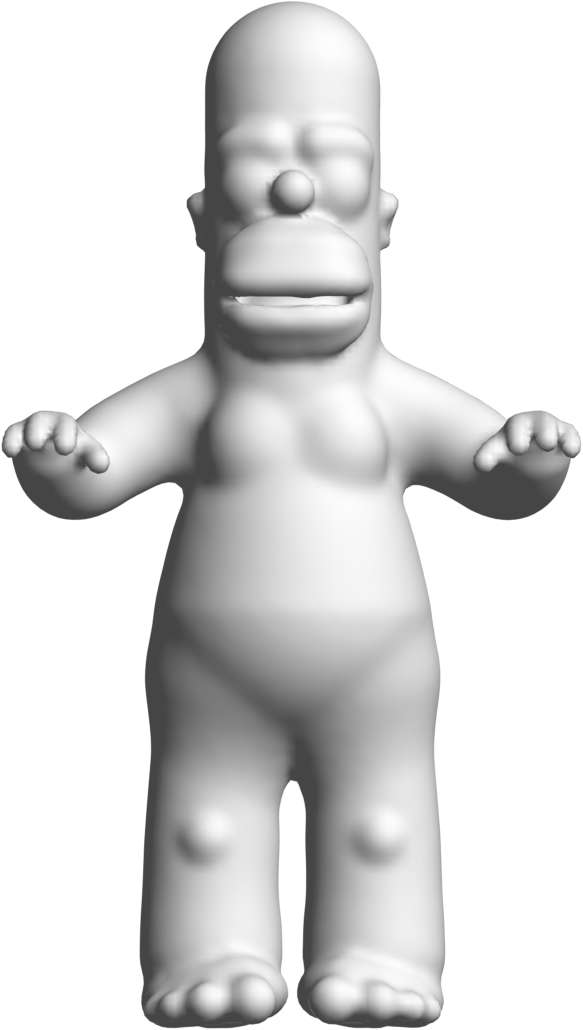} & 
\includegraphics[width=0.191\textwidth]{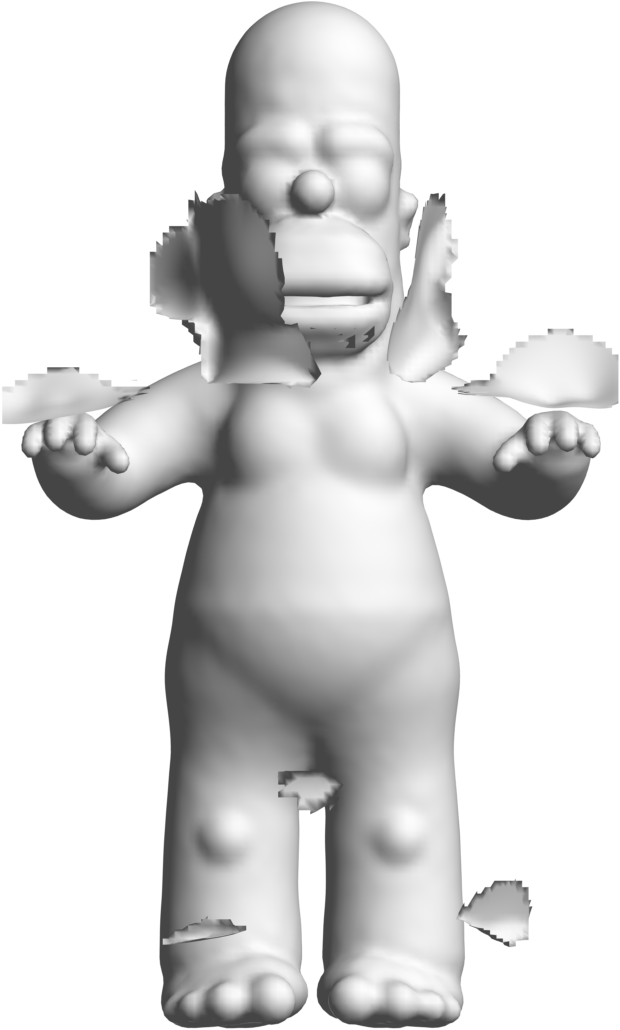} & 
\includegraphics[width=0.18\textwidth]{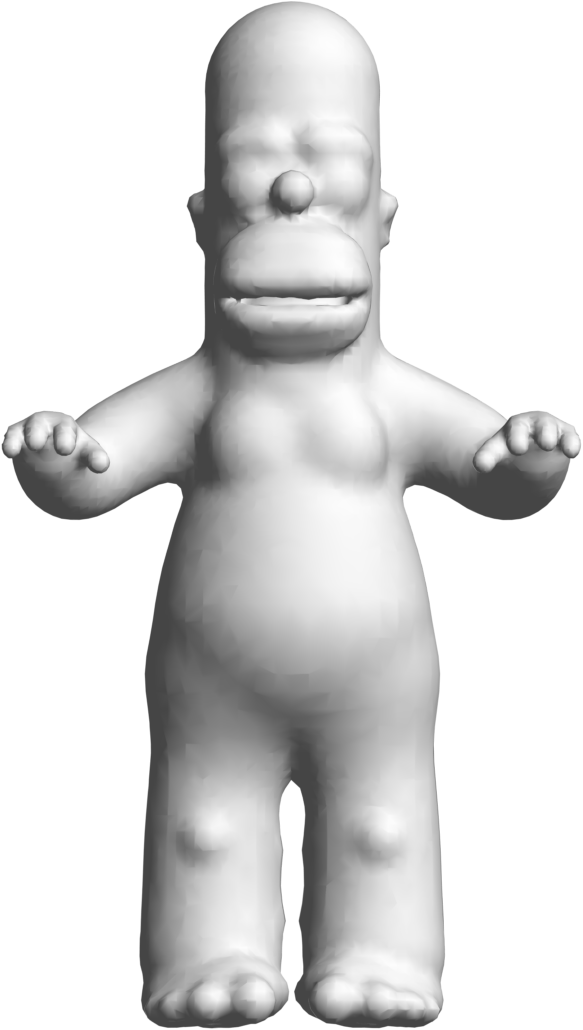} \\
\end{tabular}
\caption{\revision{Comparisons of different reconstructions of Homer.\label{fig:homer}}}
\end{figure}

We first compare the CFPU method with the global curl-free method discussed in Section \ref{sec:globalCF}, but now for reconstructing a 3D surface.  \revision{As mentioned previously, the global method is prohibitively expensive for large data sets, so we restrict our example to the small Homer point cloud shown on the left of Figure \ref{fig:homer}, which consists of $N=5103$ points.  The resulting reconstructions using $\ell=1$ and $2$ for the two curl-free methods are shown in the second and third columns of Figure \ref{fig:homer}, respectively.  We can see from the figure that for $\ell=1$, both methods produce visually pleasing, essentially indistinguishable level-surfaces from the point cloud.  However, when using the smoother $\ell=2$ kernel, both methods result in spurious sheets.  While these can be resolved by including some regularization of the normals, we do not present the results here, as later tests will experiment with this procedure.}

\revision{Reconstructions of Homer using the SPR and SSD methods are displayed in the fourth column of Figure \ref{fig:homer}.  We see from the figure that these methods result in rougher, less sharp surfaces than the curl-free methods with $\ell=1$.  This may be due to the point cloud having sparse sampling in certain regions~\cite{Screen_Pois}.  The curl-free methods seem to be less sensitive to these issues.}

%\begin{figure}[htb]
%\centering
%\begin{tabular}{ccc}
%Point  & \multicolumn{2}{c}{Global CF method} \\
%cloud & $\ell=1$ & $\ell=2$  \\ 
%\includegraphics[width=0.18\textwidth]{HomerPointCloudN5103.pdf} &
%\includegraphics[width=0.18\textwidth]{homer_global_order1_interp.png} & 
%\includegraphics[width=0.18\textwidth]{homer_global_order2_interp.png} \\
%Screened  & \multicolumn{2}{c}{CFPU method}  \\
%Poisson & $\ell=1$ & $\ell=2$  \\ 
%\includegraphics[width=0.18\textwidth]{homer_poissonrecon} & \includegraphics[width=0.18\textwidth]{homer_cfpu_order1_interp.png} & 
%\includegraphics[width=0.191\textwidth]{homer_cfpu_order2_interp.png} \\
%\end{tabular}
%\caption{\revision{Comparisons of the global curl-free and CFPU methods for different orders $\ell$ and the screened Poisson method for reconstructing Homer.\label{fig:homer}}}
%\end{figure}

%\begin{figure}[tbh]
%\centering
%\begin{tabular}{cc}
%\includegraphics[width=0.41\textwidth]{tooth_nodes_N3136.png} & 
%\includegraphics[width=0.37\textwidth]{tooth_cfpu_N3136_M468.png} \\
%(a) Point cloud \& oriented normals  & (b) CFPU reconstruction
%\end{tabular}
%\caption{(a) $N=3136$ point cloud and corresponding normals for the tooth.  (b) CFPU reconstruction of the tooth from the data in part (a).\label{fig:tooth}}
%\end{figure}
\subsection{Accuracy of CFPU}
We consider a pipe surface of radius 0.7 generated from a (2,5) torus knot.  The 3D curve defining the knot can be written parametrically as
\begin{align*}
(x(t),y(t),z(t)) = (\cos(2t)(\cos(5t)+3),\sin(2t)(\cos(5t)+3),\sin(5t)),
\end{align*}
\begin{figure}[htb]
\centering
\begin{tabular}{cc}
\includegraphics[width=0.22\textwidth]{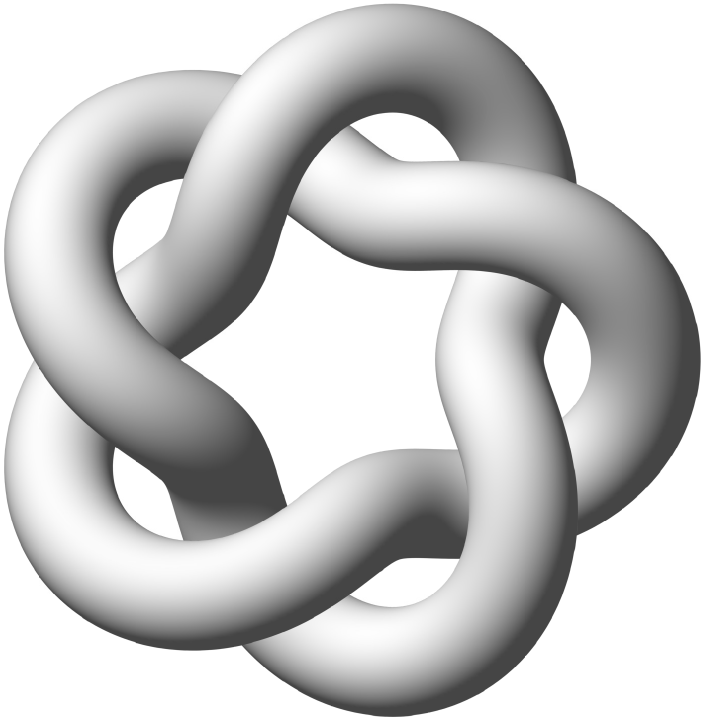} & \includegraphics[width=0.22\textwidth]{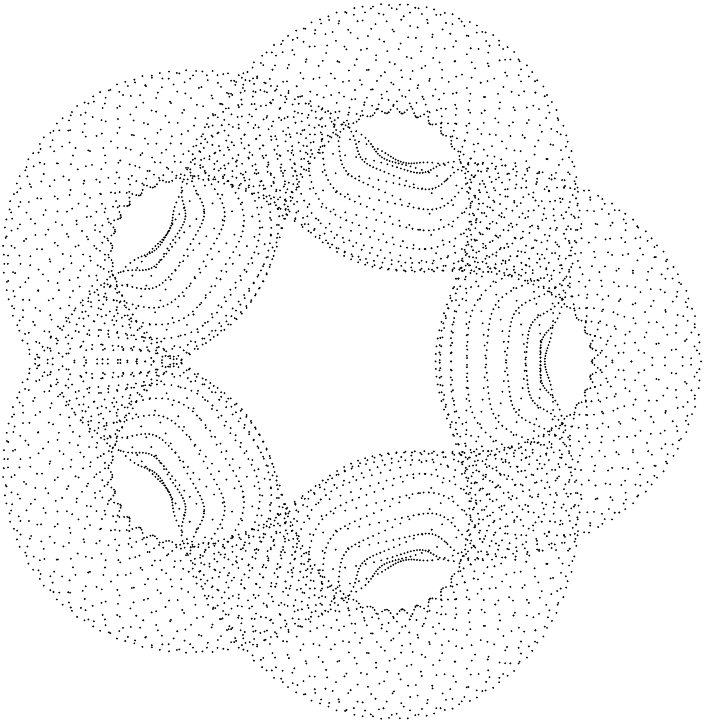}  \\
(a) & (b)
\end{tabular}
\caption{(a) Illustration of the knot used in the accuracy and noisy normals tests. (b) Example of a synthetic point cloud for the knot with $N=6144$ samples used in the tests.\label{fig:knot}}
\end{figure}
where $0\leq t \leq 2\pi$.  From this parametric representation, we generate a point cloud $\Pts=\{\pts_1,\ldots,\pts_N\}$ and the (unit) normals for the knot (see Figure \ref{fig:knot} for an illustration), and then use these in the CFPU method to reconstruct a potential whose zero-level surface approximates the exact knot surface. To test the accuracy of these reconstructions, we generate a dense set of $131424$ points exactly on the knot and compute the \refone{root mean square (RMS)} of the differences between the potential evaluated at these points and the exact solution (which is zero).

Table \ref{tbl:knot_errors} displays the results for this test for $\ell=1$ and $\ell=2$ and for increasing numbers of sample points.  In these tests, the synthetic point clouds were generated so that the average spacing between points decreases like $N^{-1/2}$.  \refone{So, we also include estimates of the algebraic convergence rates $\nu$ of the reconstructions assuming the RMS error is $\mathcal{O}(N^{-\nu / 2})$.}  We see from the table that for both $\ell=1$ and $\ell=2$, the CFPU reconstructions appear to be converging to \revision{the} true zero-level surfaces as the density of the samples increases.  \refone{Furthermore, the reconstructions using $\ell=2$ appear to be converging at an approximate rate of $\nu=5$, while they appear to be converging at the approximate rate of $\nu=3$ for $\ell=1$.  The faster rate is expected for $\ell=2$ since the knot is a smooth surface and higher values of $\ell$ produce smoother approximants.}  

\reftwo{More detailed results on the accuracy in approximating general curl-free fields (and their potentials) in two and three dimensions with curl-free RBF interpolants can be found in~\cite{DFW2020a}.  For example, it is shown that for certain piecewise smooth radial kernels, the approximation error typically decreases when more points are enclosed in each patch.  However, the rates of convergence remain the same.}

\begin{table}[t]
\centering
\begin{tabular}{|c|cc|cc|}
\hline
    & \multicolumn{2}{c|}{$\ell=1$} & \multicolumn{2}{c|}{$\ell=2$} \\
$N$ & RMS error & Rate & RMS error & Rate \\
\hline
\hline
$6144$ & $2.92\times 10^{-4}$ & -- & $1.88\times 10^{-5}$ & -- \\
$8664$ & $1.67\times 10^{-4}$ & $3.26$ & $8.60\times 10^{-6}$ & $4.56$ \\
$11616$ & $1.09\times 10^{-4}$ & $2.88$ & $4.21\times 10^{-6}$ & $4.87$ \\
$18816$ & $5.05\times 10^{-5}$ & $3.20$ & $1.23\times 10^{-6}$ & $5.12$ \\
$23064$ & $3.80\times 10^{-5}$ & $2.81$ & $7.46\times 10^{-7}$ & $4.89$ \\
$27744$ & $2.88\times 10^{-5}$ & $2.99$ & $4.73\times 10^{-7}$ & $4.92$ \\
$32856$ & $2.19\times 10^{-5}$ & $3.25$ & $3.08\times 10^{-7}$ & $5.07$ \\
\hline
\end{tabular}
\caption{\refone{Comparison of the errors and estimated algebraic convergence rates in the CFPU reconstruction of the knot for increasing numbers of samples $N$ using the curl-free PHS kernel $\Phicurl_{\ell}$, for $\ell=1,2$.  All results use a fixed number of $M=864$ PU patches.\label{tbl:knot_errors}}}
\end{table}

%\begin{table}[h]
%\centering
%\begin{tabular}{|c|cc|cc|}
%\hline
%    & \multicolumn{2}{c|}{$\ell=1$} & \multicolumn{2}{c|}{$\ell=2$} \\
%$N$ & RMS error & Rate & RMS error & Rate \\
%\hline
%\hline
%$6144$ & $9.90\times 10^{-5}$ & -- & $8.08\times 10^{-6}$ & -- \\
%$8664$ & $4.12\times 10^{-5}$ & $5.10$ & $3.45\times 10^{-6}$ & $4.95$ \\
%$11616$ & $2.27\times 10^{-5}$ & $4.07$ & $1.69\times 10^{-6}$ & $4.85$ \\
%$18816$ & $7.46\times 10^{-6}$ & $4.61$ & $4.53\times 10^{-7}$ & $5.47$ \\
%$23064$ & $5.21\times 10^{-6}$ & $3.52$ & $2.67\times 10^{-7}$ & $5.22$ \\
%$27744$ & $4.10\times 10^{-6}$ & $2.61$ & $1.68\times 10^{-7}$ & $4.97$ \\
%$32856$ & $2.87\times 10^{-6}$ & $4.21$ & $9.69\times 10^{-8}$ & $6.54$ \\
%\hline
%\end{tabular}
%\caption{\refone{Comparison of the errors in the CFPU reconstruction of the knot for increasing numbers of samples $N$ using the curl-free PHS kernel $\Phicurl_{\ell}$, for $\ell=1,2$.  Also included are the estimated algebraic convergence rates $\nu$ of the reconstructions assuming $\text{err} \sim N^{-\nu/2}$.  All results use a fixed number of $M=864$ PU patches.\label{tbl:adfad}}}
%\end{table}

 %%%%%%%%%%%%%%%%%%%%%%%%%%%%%%%%%%%%%%%%%%%%%%%%%%%%%%%%%%%%
%%%%%%%%%%%%%%%%%%%%%%%%%%%%%%%%%%%%%%%%%%%%%%%%%%%%%%%%%%%%%
\subsection{Noisy normals}
\begin{figure}[h!]
\centering
\begin{tabular}{ccccc||c}
& \multicolumn{4}{c||}{\small CFPU} & \\ 
& {\tiny No reg.\ ($\lambda=0$)} & {\tiny  $\lambda=10^{-4}$} & {\tiny  $\lambda=10^{-2}$} & {\tiny $\lambda$ GCV} & {\small \revision{SPR}} \\
\rotatebox{90}{\hspace{0.05\textwidth} \small $\ell=1$} &
\includegraphics[width=0.163\textwidth]{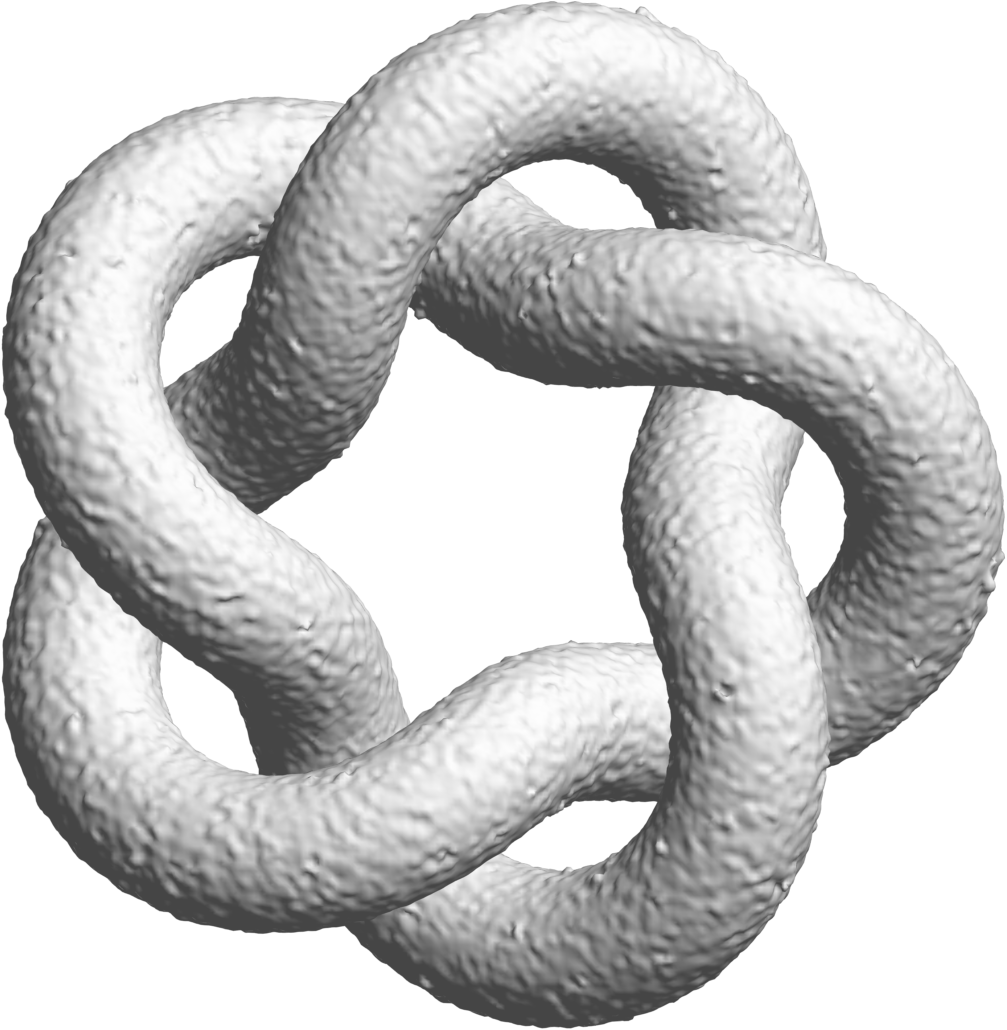} & 
\includegraphics[width=0.163\textwidth]{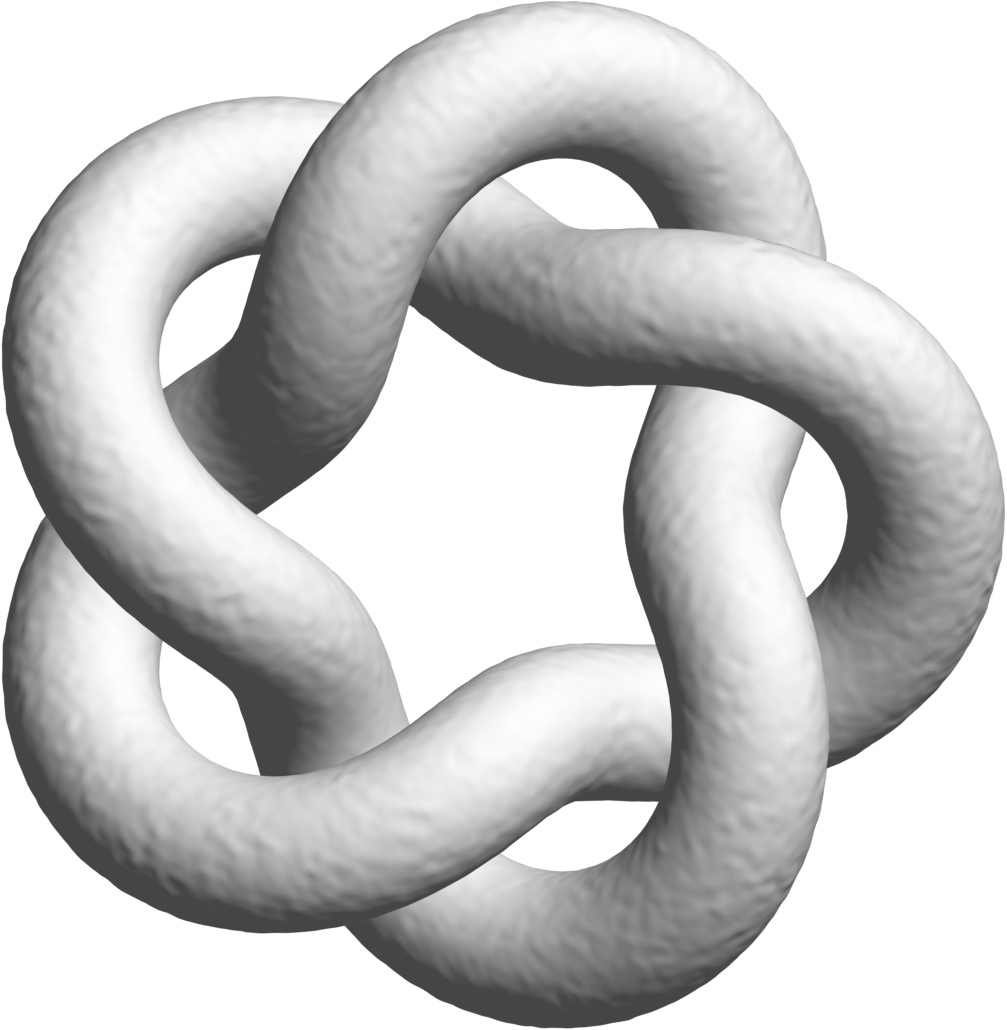} &
\includegraphics[width=0.163\textwidth]{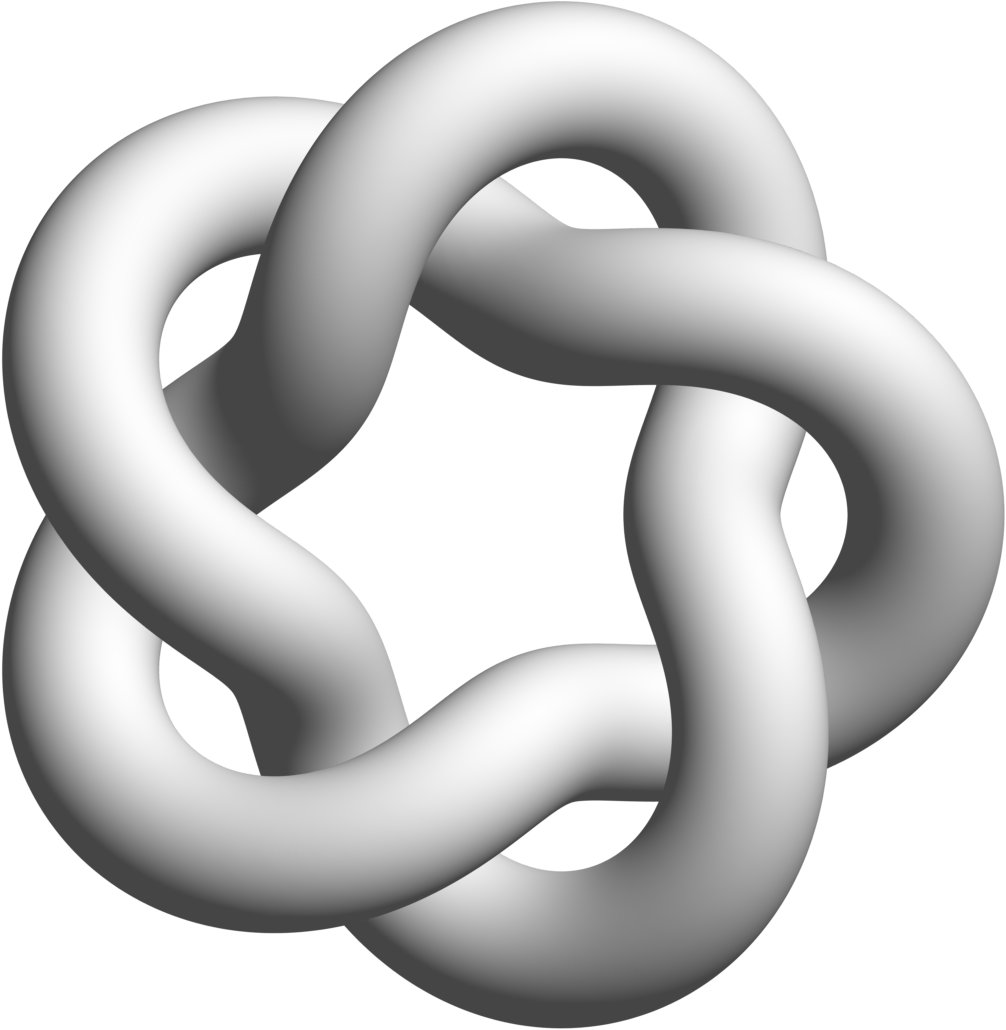} &
\includegraphics[width=0.163\textwidth]{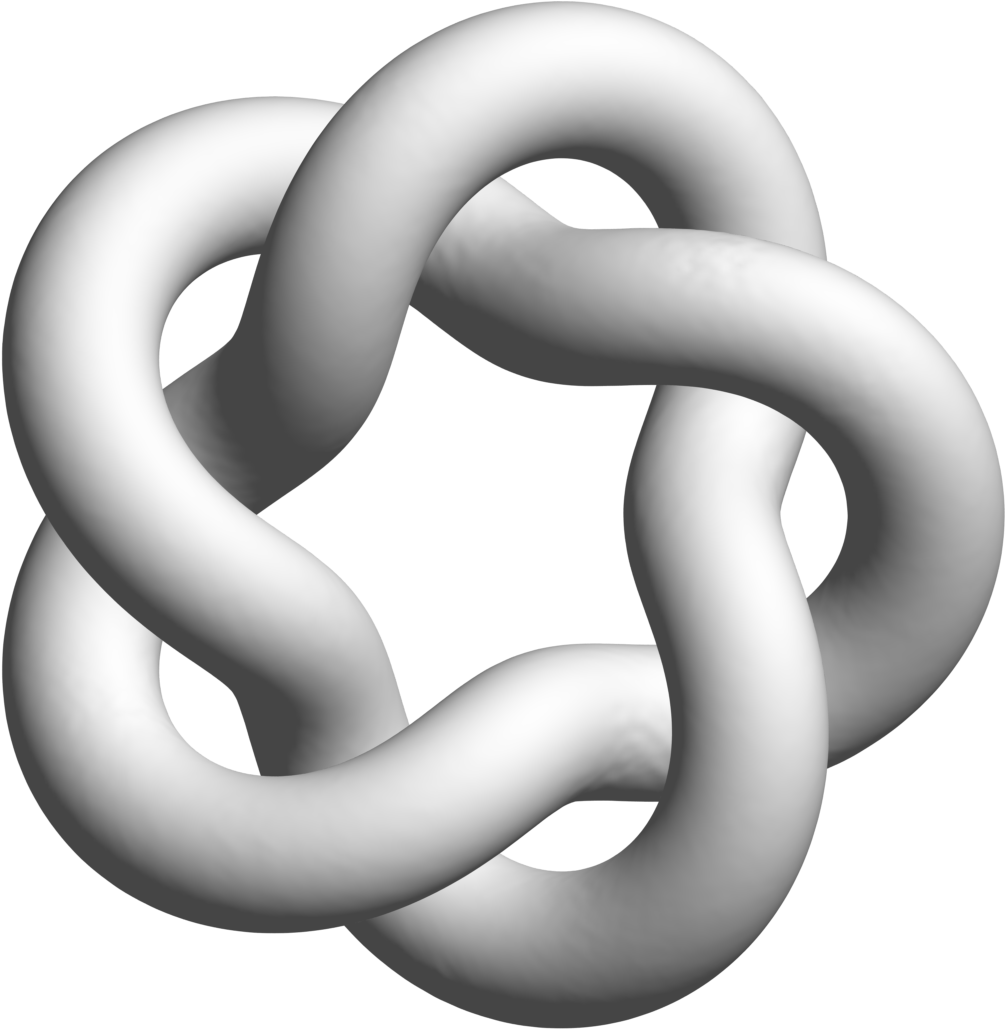}  & \includegraphics[width=0.163\textwidth]{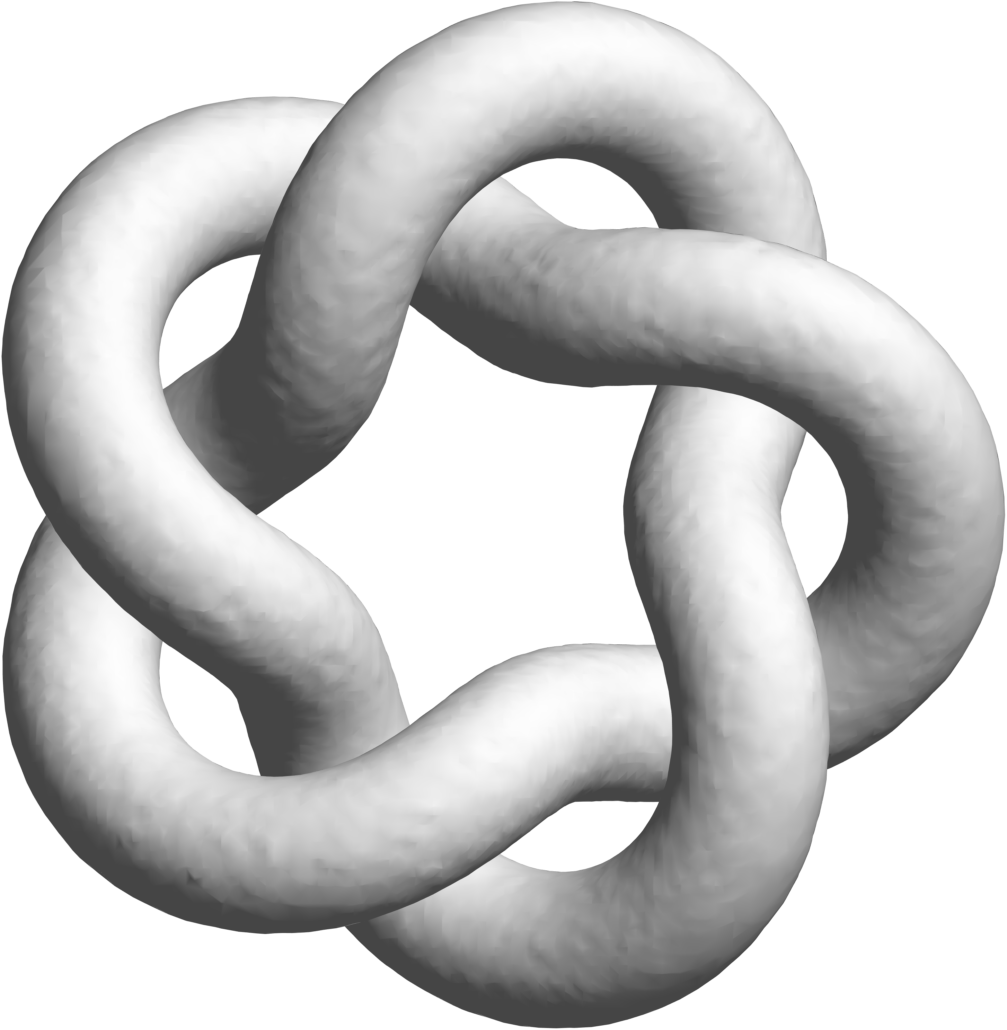} \\ 
& {\tiny No reg.\ ($\lambda=0$)} & {\tiny  $\lambda=10^{-7}$} & {\tiny  $\lambda=10^{-5}$} & {\tiny $\lambda$ GCV}  & {\small \revision{SSD}} \\
\rotatebox{90}{\hspace{0.05\textwidth} \small $\ell=2$ }&
\includegraphics[width=0.163\textwidth]{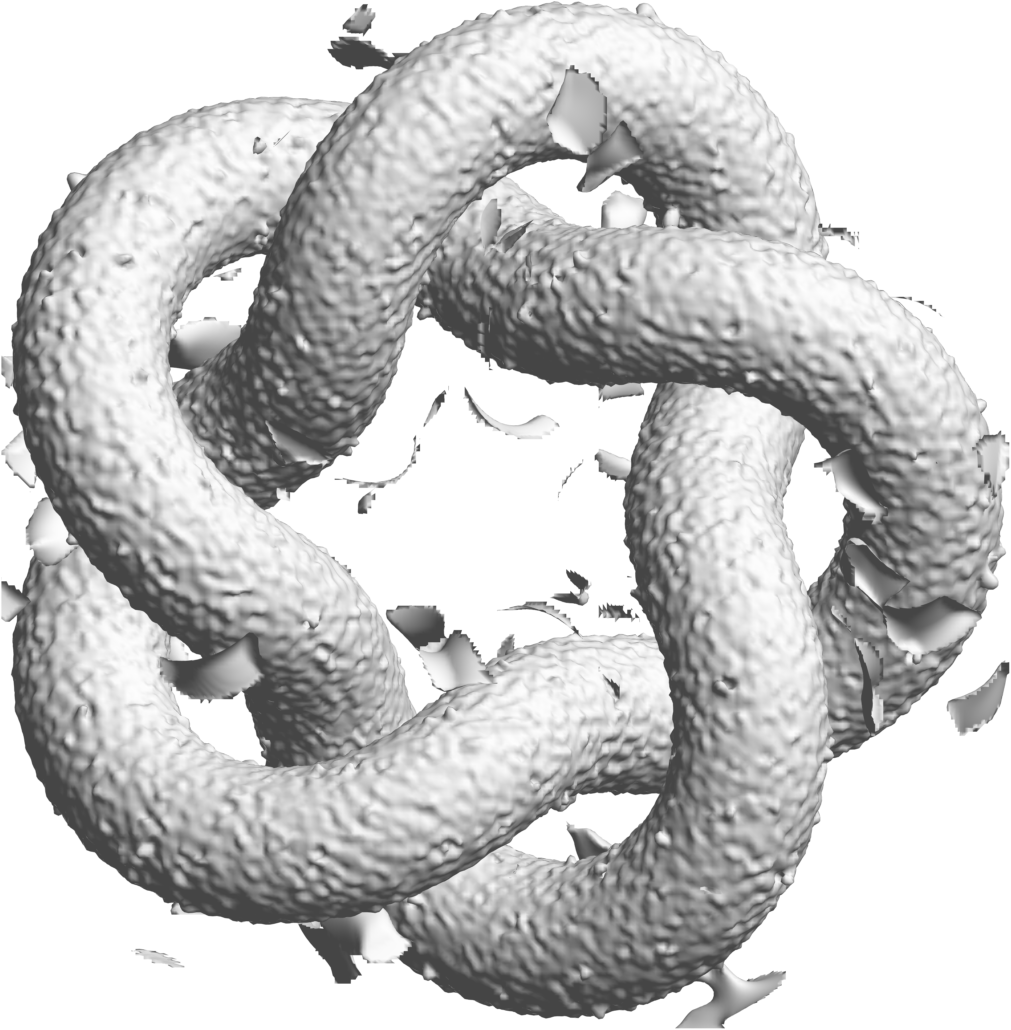} & 
\includegraphics[width=0.163\textwidth]{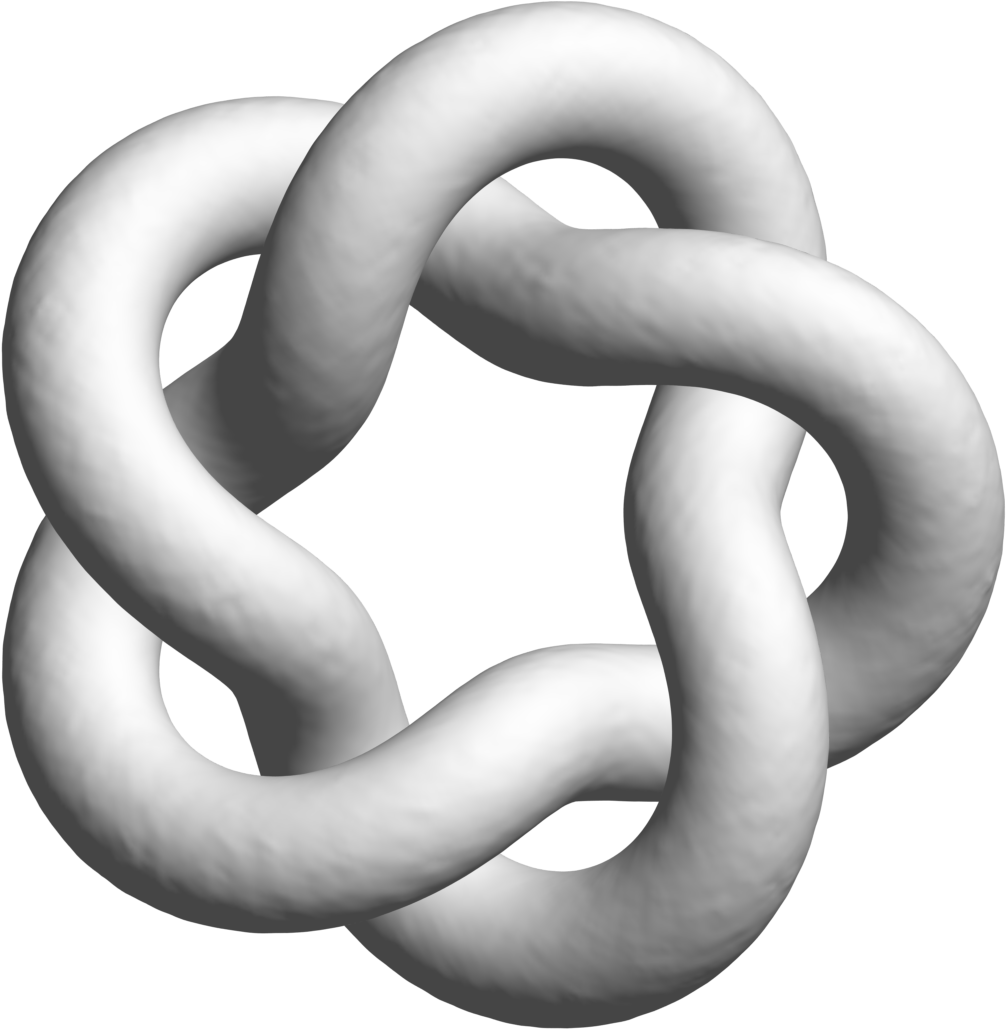} &
\includegraphics[width=0.163\textwidth]{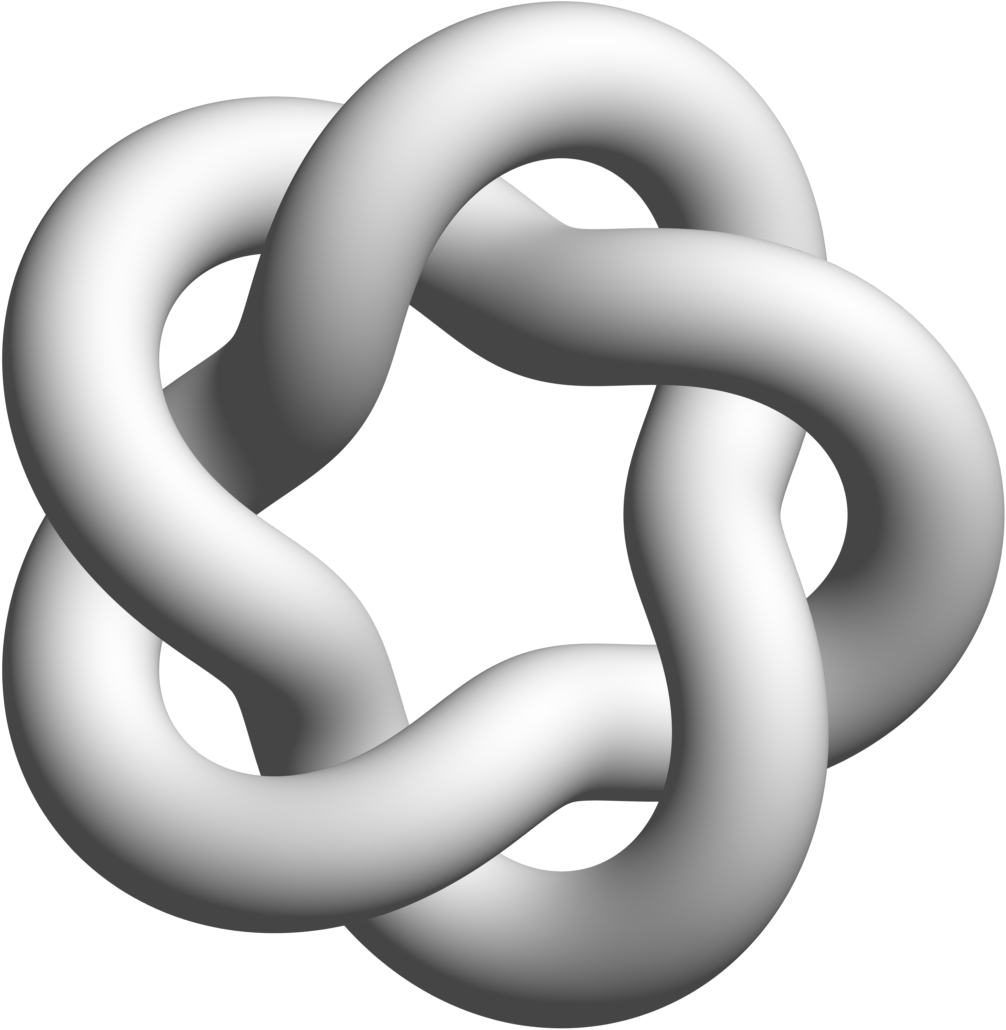} &
\includegraphics[width=0.163\textwidth]{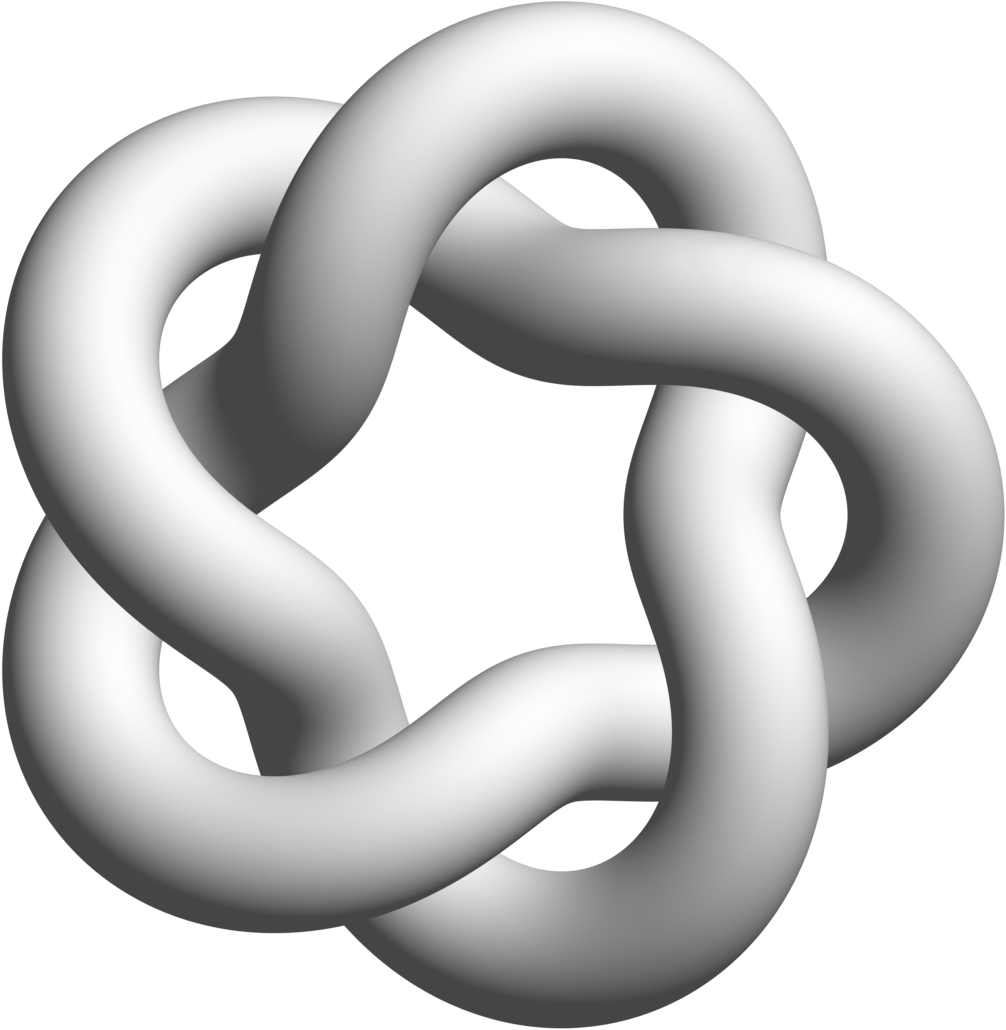} & \includegraphics[width=0.163\textwidth]{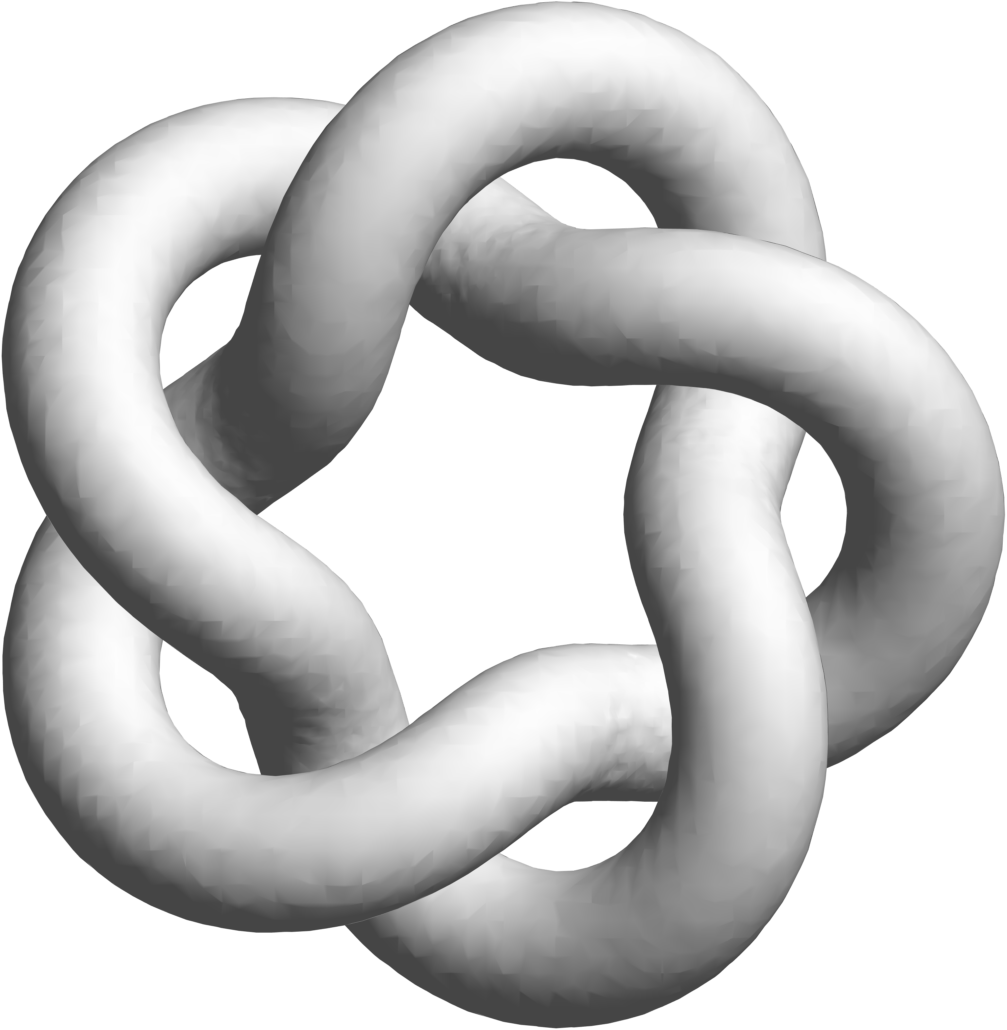}
\end{tabular}
\caption{\revision{Comparison of reconstructions of the knot with zero mean Gaussian white noise added to the normals. First column shows the CFPU reconstructions without any regularization of the normals.  Second and third columns show the CFPU reconstructions using regularization with a fixed parameter $\lambda$ chosen for all the patches.  Fourth column shows the CFPU reconstructions with the regularization parameter chosen using GCV on each patch.  Fifth column shows the surfaces from the SPR and SSD methods for comparison purposes.  All results use $N=23064$.\label{fig:knot_noise}}}
\end{figure}

In this test we demonstrate how the regularization method from Section \ref{sec:nrml_regularization} can help with noise in the normals; we do not test the residual regularization method here.  We use the knot example from the previous section and add noise to the exact (unit) normals $\{\vn_1,\ldots,\vn_N\}$ according to
%\begin{align*}
$\vn_j^* = \vn_j + \boldsymbol{\epsilon}_j$,
%\end{align*}
where $\boldsymbol{\epsilon}_j\in\R^3$ and each of its components is an independent normally distributed random variable with mean zero and standard deviation 0.3.  The results are displayed in Figure \ref{fig:knot_noise}. 

The first column of this figure shows the reconstructed surfaces for both $\ell=1$ and $\ell=2$ when no regularization is applied.  Both results give a noisy surface, but only the $\ell=2$ results give spurious sheets.  The second and third columns show the results when a fixed regularization parameter $\lambda$ is used for all patches.  We see from the figure that including regularization results in a smoother surfaces and increasing $\lambda$ increases the smoothness.  These columns also show that larger values of $\lambda$ are needed for $\ell=1$ to obtain the same level of smoothness.  \revision{The fourth column} of the figure shows the results when automating the selection of the regularization parameters using GCV on each patch.  We see from the figure that this technique appears to be quite effective at reducing the noise for both $\ell=1$ and $\ell=2$.  However, in the $\ell=1$ case, there are still some remnants of oscillations where the knot crosses over itself.

\revision{For comparison purposes, the last column of Figure \ref{fig:knot_noise} shows reconstructions of the ``noisy'' knot using the SPR and SSD methods.  To filter out the noise and produce these surfaces, we had to modify the default parameters to introduce more smoothing, similar to the CFPU method.  For the SPR method we set the \texttt{pointWeight} parameter to 10 and the \texttt{samplesPerNode} parameter to 2, while for the SSD method we set the \texttt{samplesPerNode} parameter to 4 and the \texttt{biLapWeight}  parameter to 4.  We found these parameters produced the best results. We can see from the figure that both the SPR and SSD reconstructions produce more noisy surfaces than the regularized CFPU method, with the noise being more prevalent in the former.}

\subsection{Raw range data}
\begin{figure}[t]
\centering
\begin{tabular}{ccccc||c}
& \multicolumn{4}{c||}{\small CFPU} & \\ 
& {\tiny No reg.\ ($\lambda=0$, $\alpha=0$)} & {\tiny  $\lambda=10^{-2}$, $\alpha=0$} & {\tiny  $\lambda=10^{-2}$, $\alpha=10^{-4}$} & {\tiny $\lambda$, $\alpha$ GCV} & {\small \revision{SPR}} \\
\rotatebox{90}{\hspace{0.05\textwidth} \small $\ell=1$} &
\includegraphics[width=0.163\textwidth]{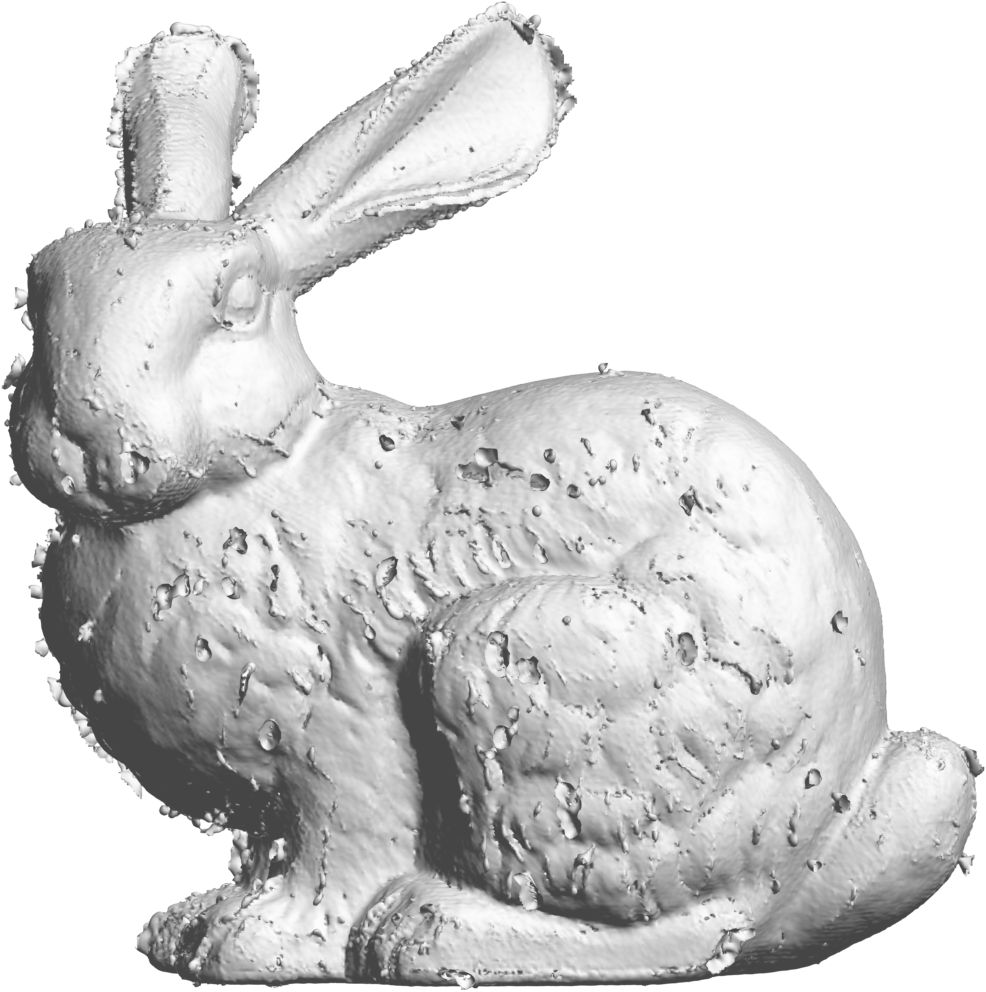} & 
\includegraphics[width=0.163\textwidth]{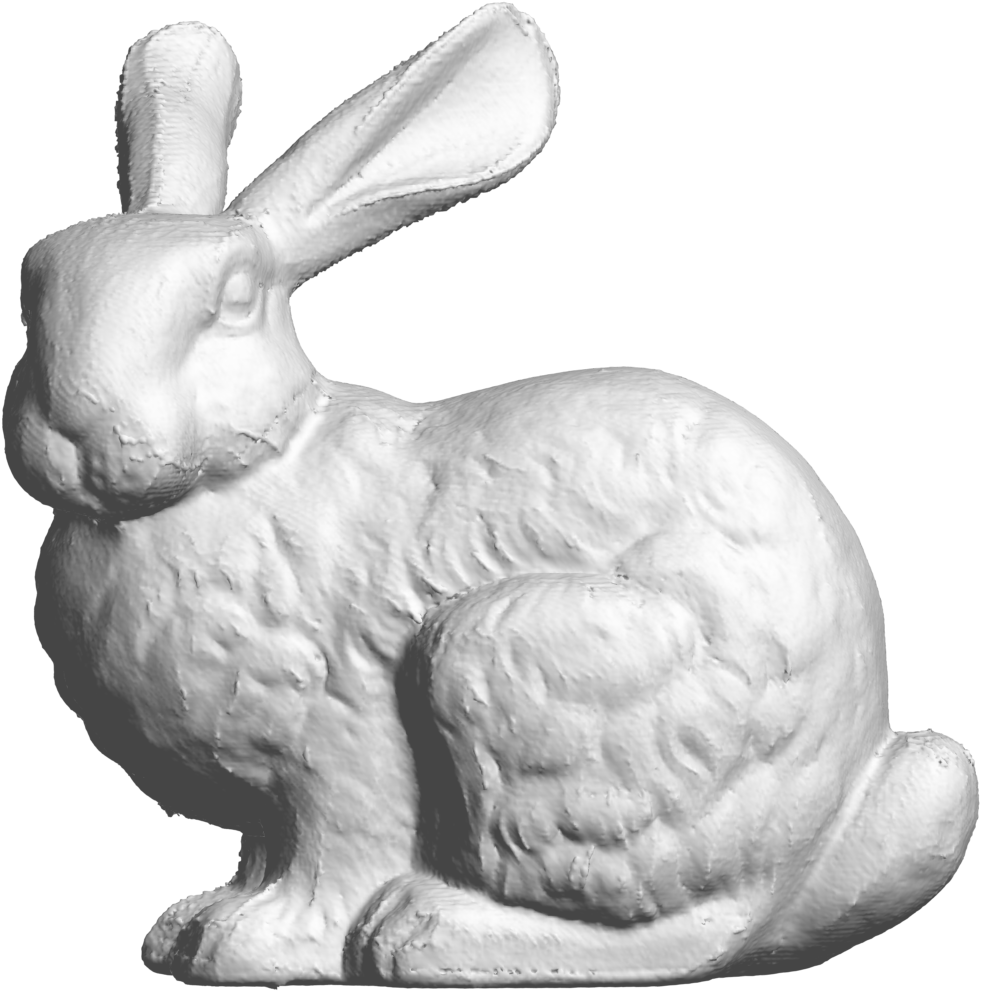} &
\includegraphics[width=0.163\textwidth]{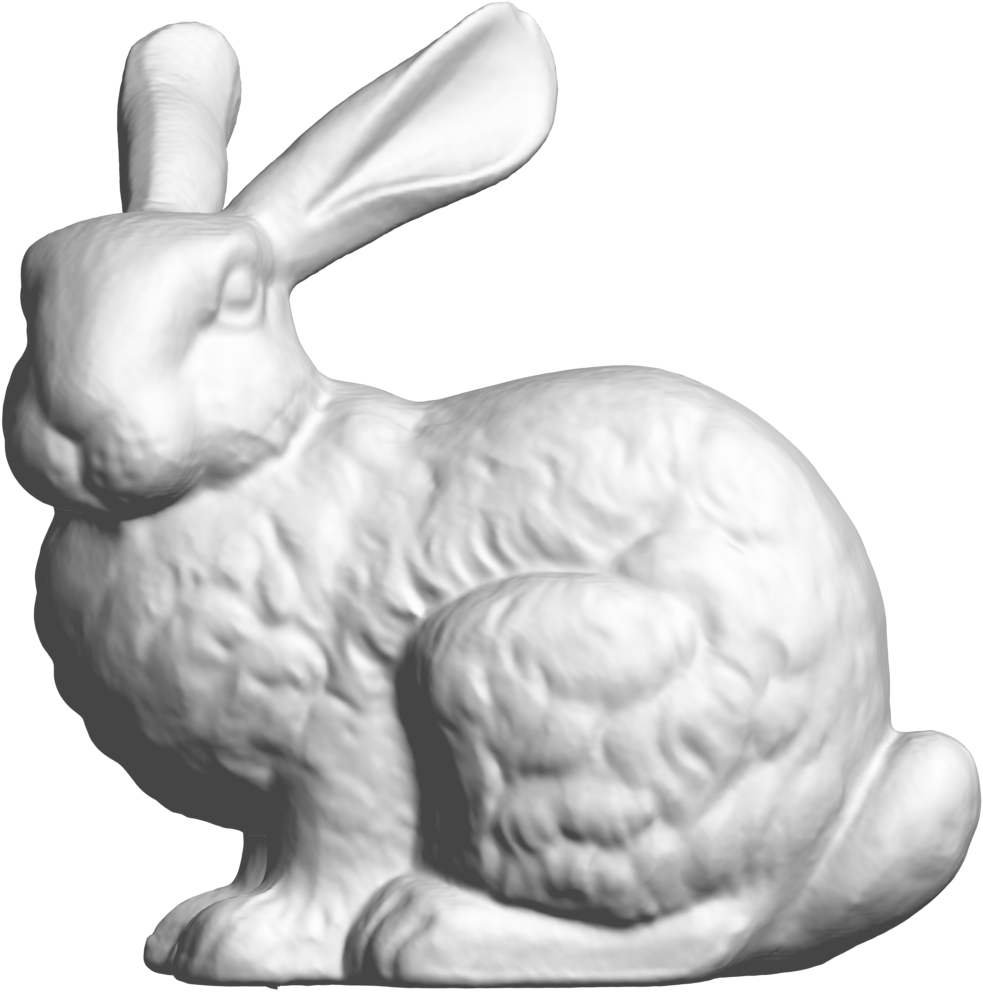} &
\includegraphics[width=0.163\textwidth]{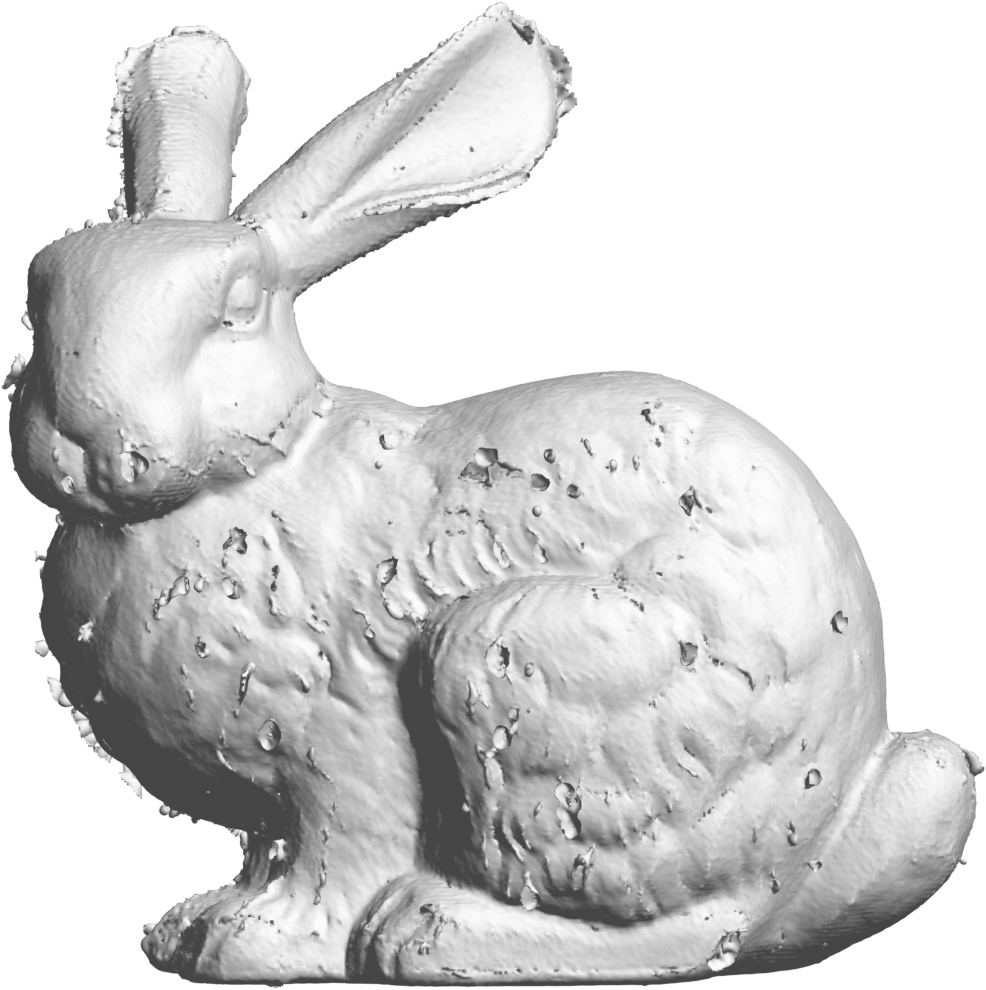} &
\includegraphics[width=0.163\textwidth]{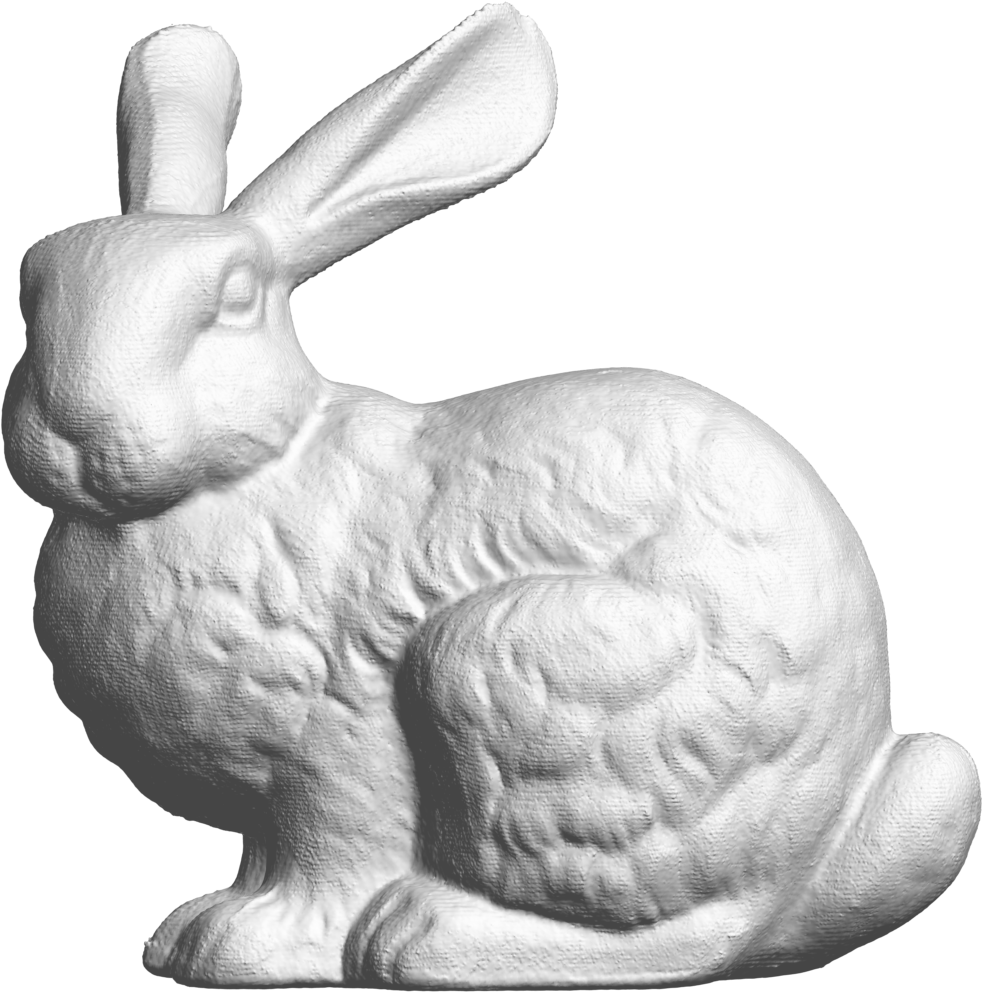} \\

& {\tiny No reg.\ ($\lambda=0$, $\alpha=0$)} & {\tiny  $\lambda=10^{-2}$, $\alpha=0$} & {\tiny  $\lambda=10^{-2}$, $\alpha=10^{-4}$} & {\tiny $\lambda$, $\alpha$ GCV} & {\small \revision{SSD}} \\
\rotatebox{90}{\hspace{0.05\textwidth} \small $\ell=2$ }&
\includegraphics[width=0.163\textwidth]{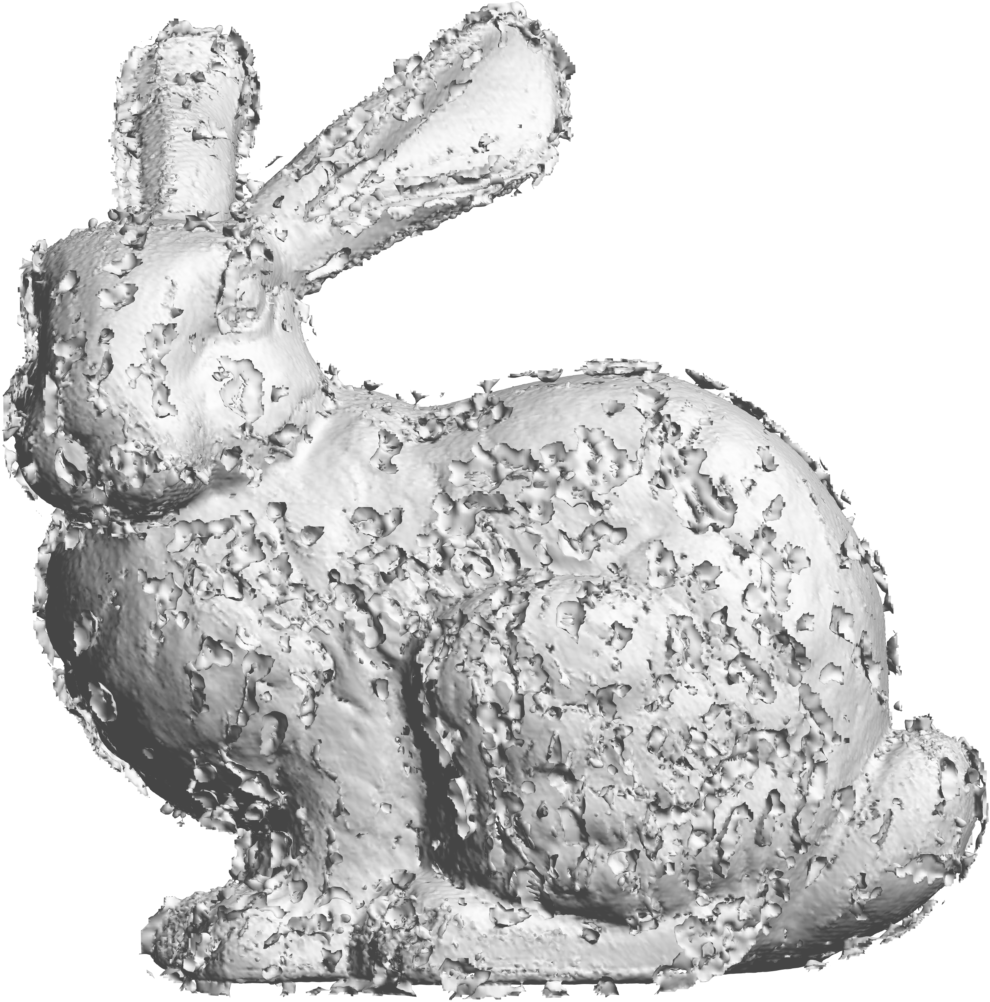} & 
\includegraphics[width=0.163\textwidth]{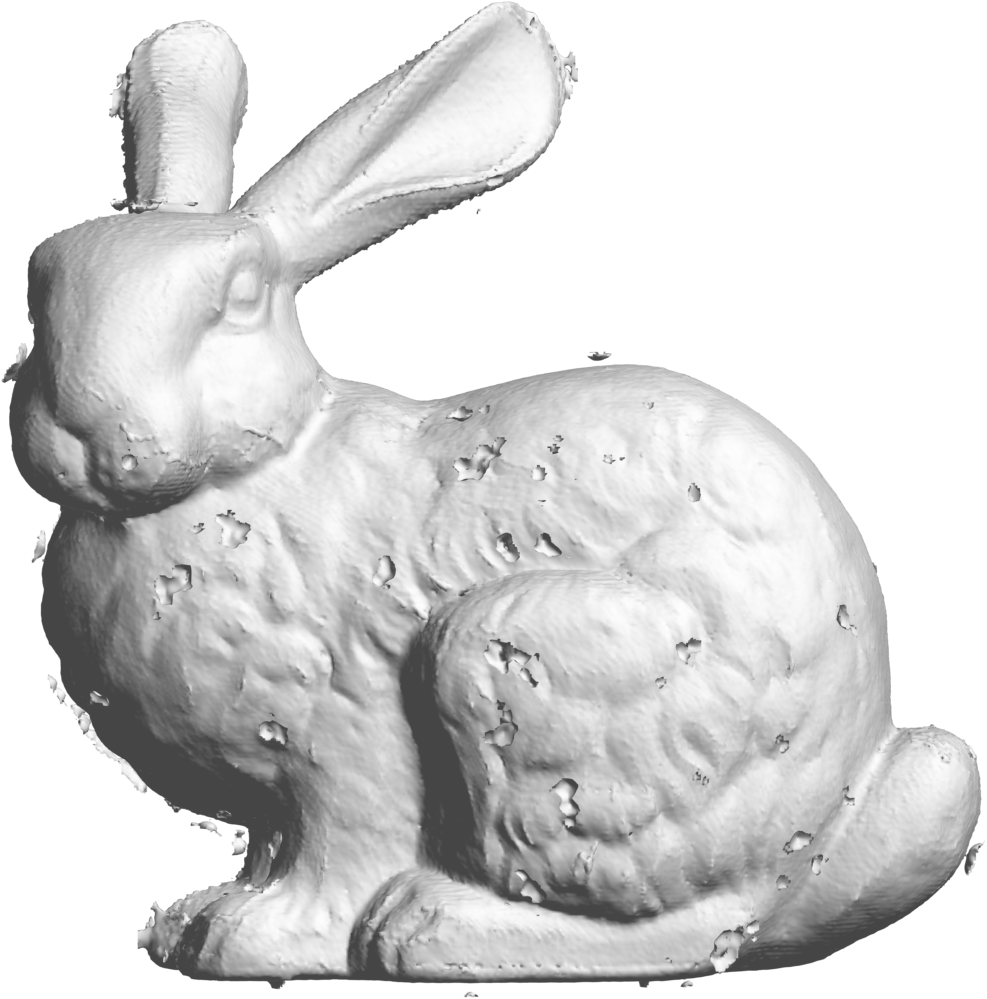} &
\includegraphics[width=0.163\textwidth]{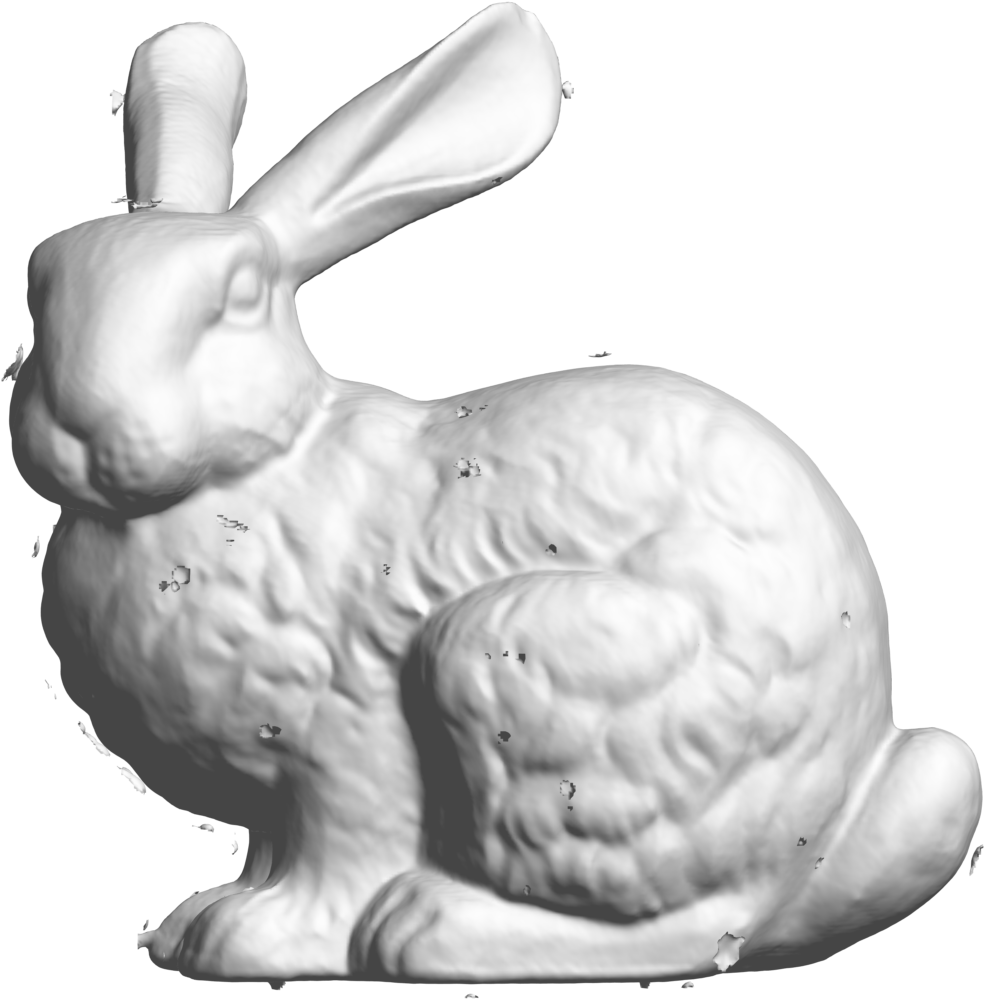} &
\includegraphics[width=0.163\textwidth]{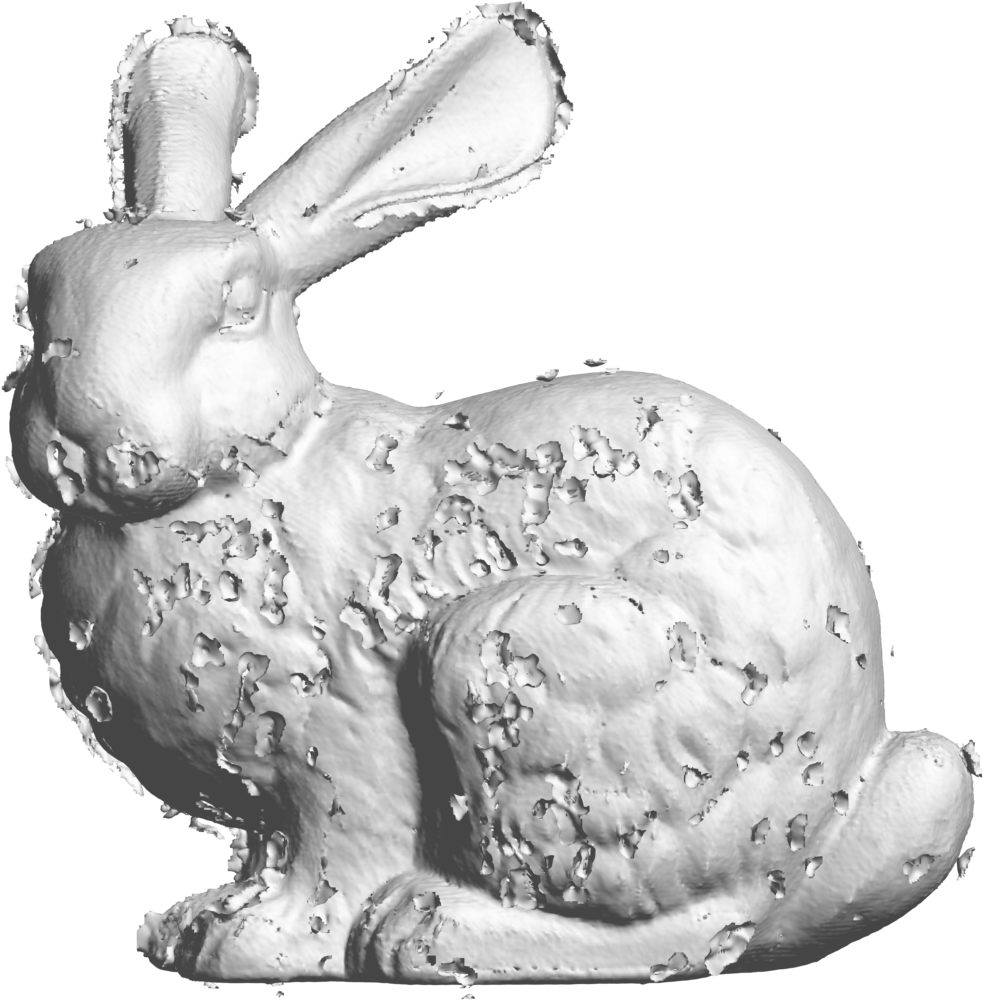} &
\includegraphics[width=0.163\textwidth]{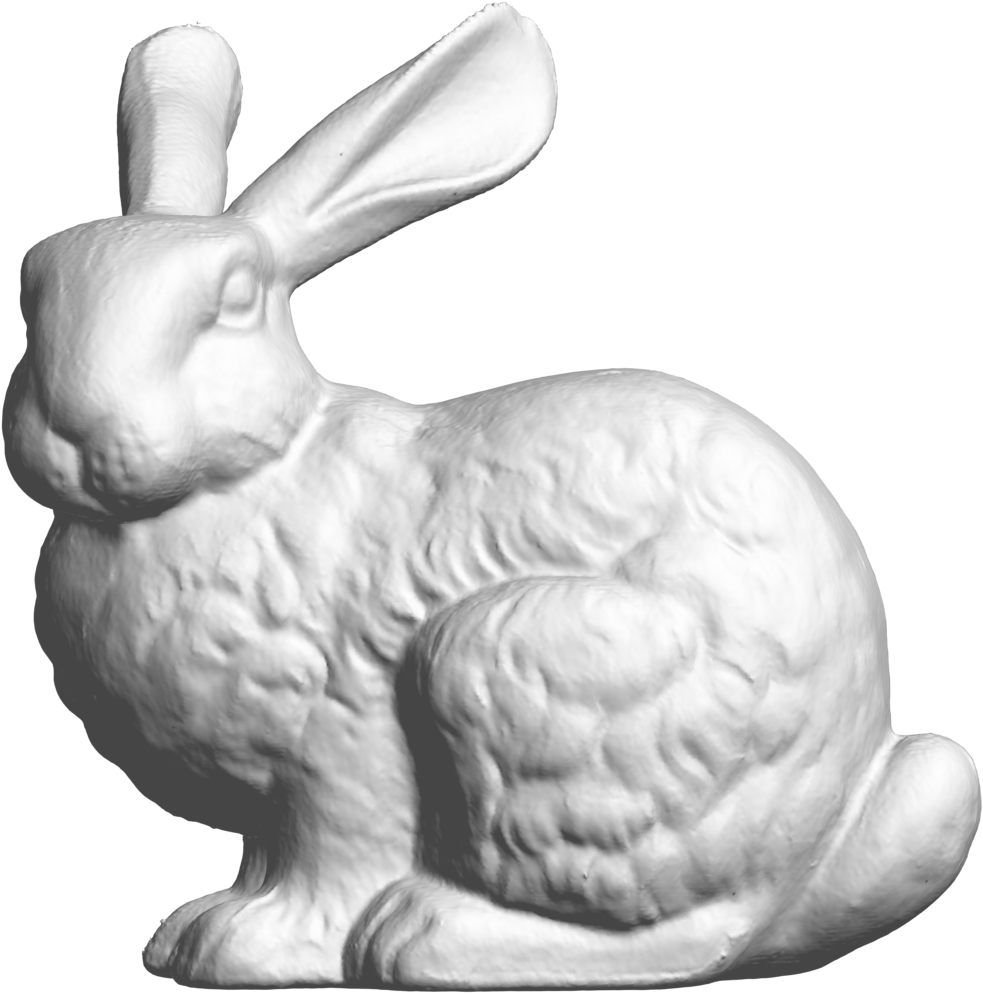} \\
\end{tabular}
\caption{Comparison of the CFPU reconstructions of the Stanford Bunny using raw range data.  First column shows the reconstructions without any regularization.  The remaining columns show the effects of regularizing the fit of the normal vectors ($\lambda > 0$) and the fit of the residual of the potential ($\alpha > 0$).  For the second and third columns, the parameters were fixed for all patches, while for the fourth column, the parameters were chosen on a per-patch basis using GCV.  \revision{The last column shows the surfaces from the SPR and SSD methods for comparison purposes. \label{fig:bunny}}}
\end{figure}

%\begin{figure}[h!]
%\centering
%\begin{tabular}{c|cccc}
%& {\tiny No reg.\ ($\lambda=0$, $\alpha=0$)} & {\tiny  $\lambda=10^{-2}$, $\alpha=0$} & {\tiny  $\lambda=10^{-2}$, $\alpha=10^{-4}$} & {\tiny $\lambda$, $\alpha$ GCV} \\
%\hline
%\rotatebox{90}{\hspace{0.05\textwidth} \small $\ell=1$} &
%\includegraphics[width=0.2\textwidth]{bunny_cfpu_order1_noreg.png} & 
%\includegraphics[width=0.2\textwidth]{bunny_cfpu_order1_reg2_1em2_regi20.png} &
%\includegraphics[width=0.2\textwidth]{bunny_cfpu_order1_reg2_1em2_regi2_1em4.png} &
%\includegraphics[width=0.2\textwidth]{bunny_cfpu_order1_reg1_gcv_regi1_gcv.png} \\
%\hline
%\rotatebox{90}{\hspace{0.05\textwidth} \small $\ell=2$ }&
%\includegraphics[width=0.2\textwidth]{bunny_cfpu_order2_noreg.png} & 
%\includegraphics[width=0.2\textwidth]{bunny_cfpu_order2_reg2_1em2_regi20.png} &
%\includegraphics[width=0.2\textwidth]{bunny_cfpu_order2_reg2_1em2_regi2_1em4.png} &
%\includegraphics[width=0.2\textwidth]{bunny_cfpu_order2_reg1_gcv_regi1_gcv.png} \\
%\end{tabular}
%\caption{Comparison of the CFPU reconstructions of the Stanford Bunny using raw range data.  First column shows the reconstructions without any regularization.  While the remaining columns show the effects of regularizing the fit of the normal vectors ($\lambda > 0$) and the fit of the residual of the potential ($\alpha > 0$).  For the second and third columns, the parameters were fixed for all patches, while for the fourth column, the parameters were chosen on a per-patch basis using GCV.\label{fig:bunny}}
%\end{figure}
We now consider reconstructing a surface from range data, which typically consists of optical scans of an object from different orientations and produces data sets that are misaligned and contain noise.  For this test, we use the Stanford Bunny problem and generate the point cloud and normals by combining the 10 available range scans.  Figure \ref{fig:bunny} shows the resulting CFPU reconstructions of the surface for $\ell=1$ and $\ell=2$ and different regularization parameters.  

We see from the first column of the figure that without any regularization, the CFPU method produces a zero-level surface with several spurious sheets and that these are more prevalent in the $\ell=2$ reconstruction than the $\ell=1$.  The second column shows the reconstruction from regularizing the fit of the normal vectors using the same parameter for each patch.  We see that this regularization reduces the spurious sheets for the $\ell=2$ case and completely eliminates them for the $\ell=1$ case.  The third column shows the results of using both regularization of the normals and the fit of the residual of the potential, again using the same parameter for all the patches.  For both values of $\ell$, we see that this added regularization results in a smoother surface with less fine-scale detail.  For the $\ell=2$ case, there are still spurious sheets present in the reconstruction and further tuning of the parameters could not eliminate these.  Finally, the \revision{fourth} column of Figure \ref{fig:bunny} shows the reconstructions when using GCV to determine both regularization parameters on a per-patch basis. Unlike the knot problem, we see from the figure that GCV is not effective at producing parameters that smooth the results and eliminate the spurious sheets.  This may be due to the noise in the data being non-Gaussian and from inconsistent noise levels across the different scans.

\revision{The final column of Figure \ref{fig:bunny} shows the reconstructions of the bunny using the SPR and SSD using the default parameters.  These surfaces are similar in terms of fine-scale detail and lack of spurious sheets to the CFPU surface with $\ell=1$ in the second column.}

\subsection{Common 3D models}
We finally test the CFPU method for reconstructing some common 3D models found in the literature.  We only present results for the method using $\ell=1$ since this gave the best reconstructions and always produced zero-level surfaces without spurious sheets.  \revision{We also only include comparisons to the SPR method, since the SSD method gave similar results.}

\begin{figure}[htb]
\setlength{\tabcolsep}{2pt}
\centering
\begin{tabular}{ccccc}
\rotatebox{90}{\hspace{0.08\textwidth} \small \revision{CFPU}} & \phantom{H} &
\includegraphics[width=0.3\textwidth]{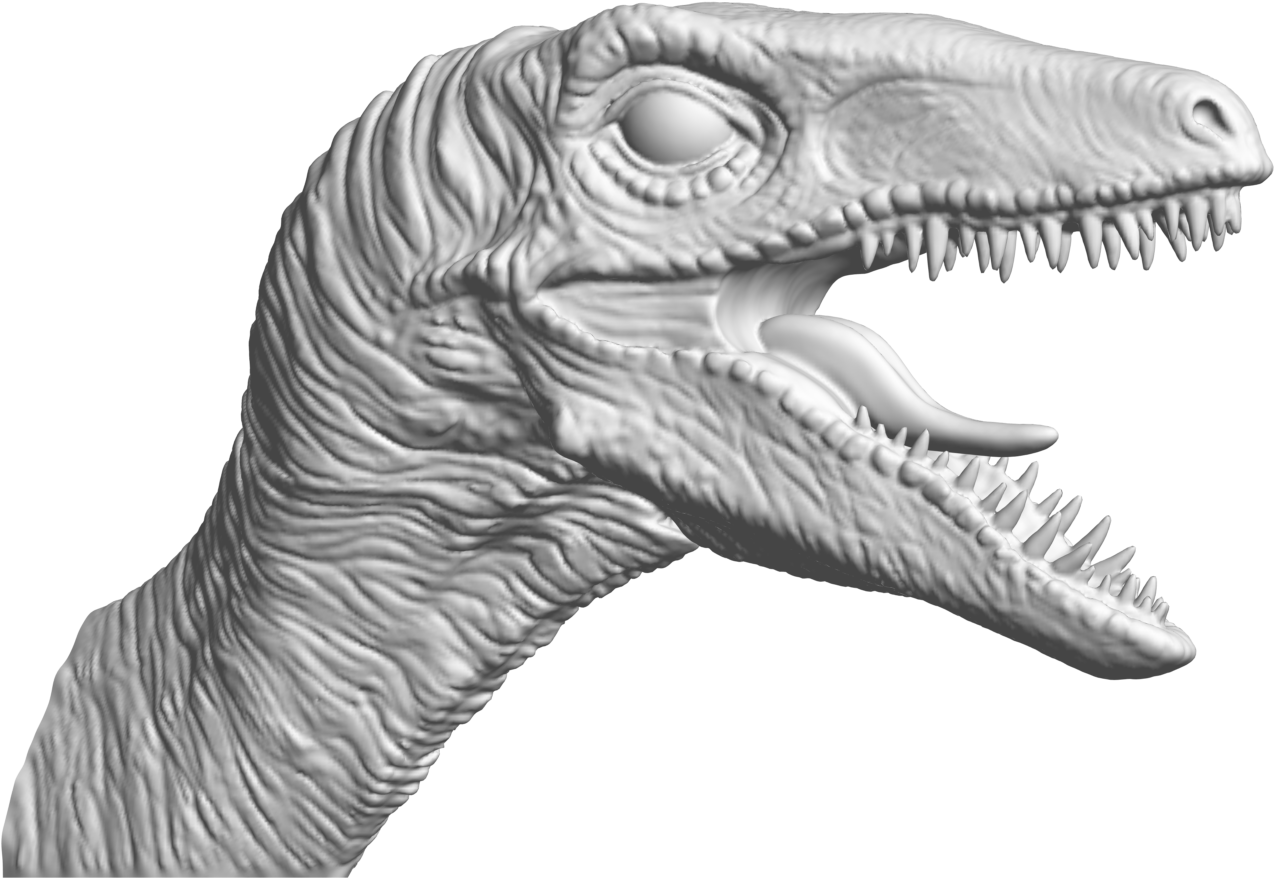} & \includegraphics[width=0.25\textwidth]{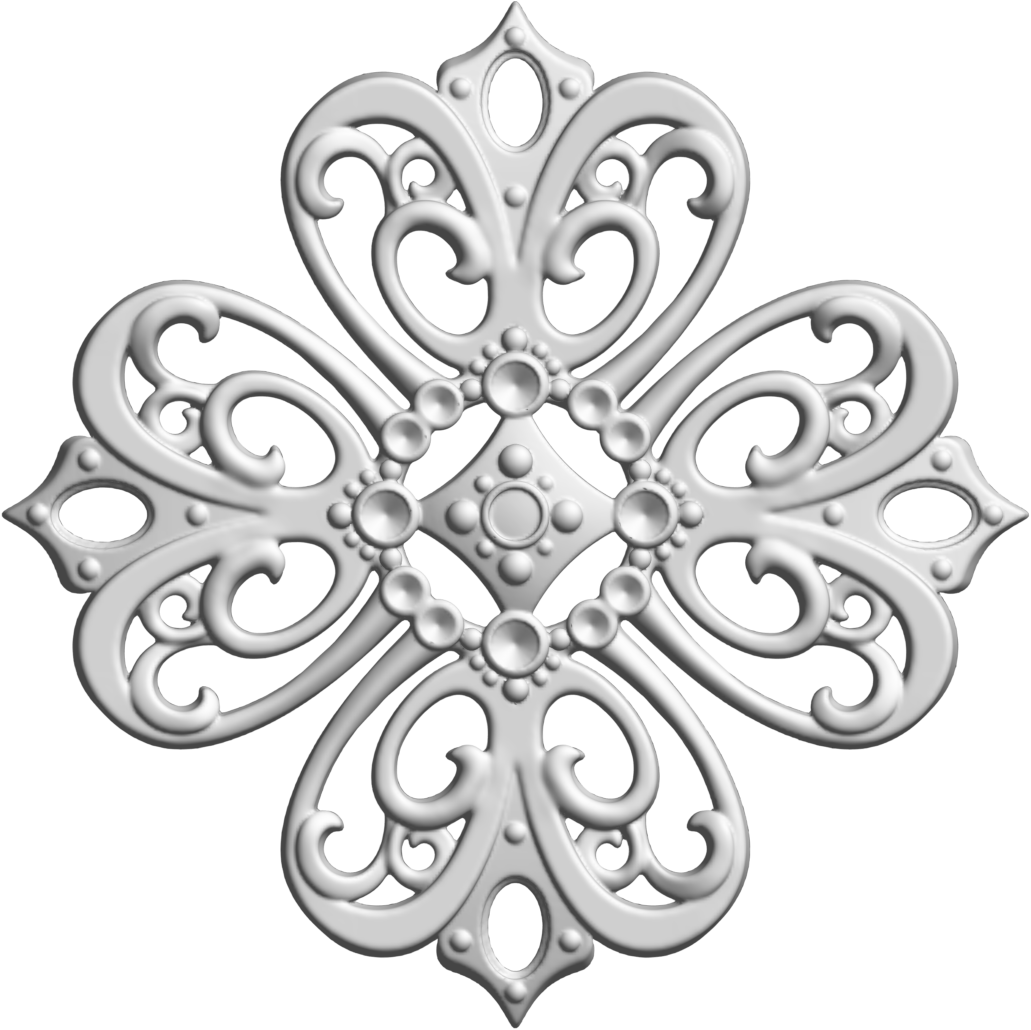} & 
\includegraphics[width=0.3\textwidth]{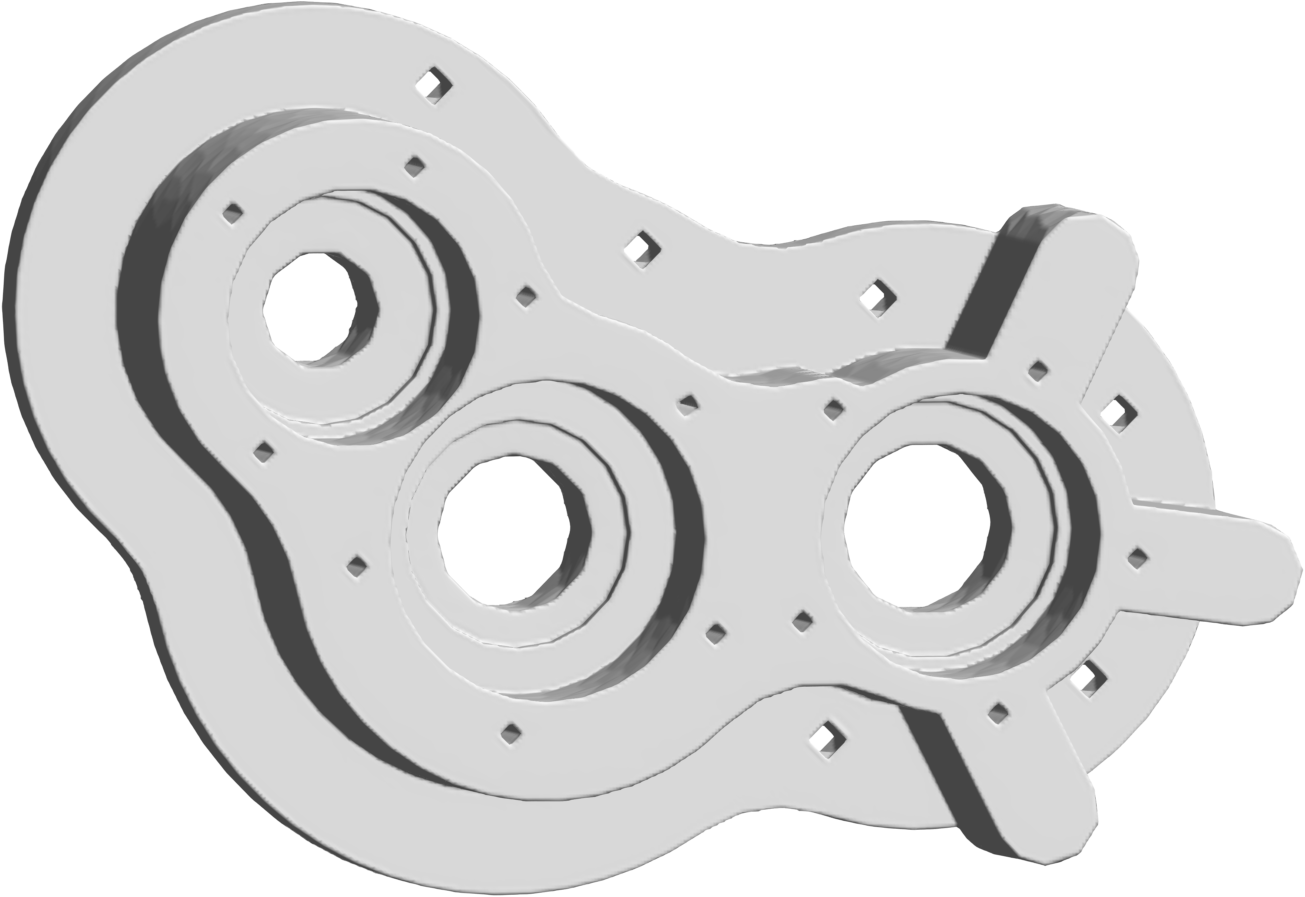} \\
\rotatebox{90}{\hspace{0.08\textwidth} \small \revision{SPR}}  & \phantom{H} &
\includegraphics[width=0.3\textwidth]{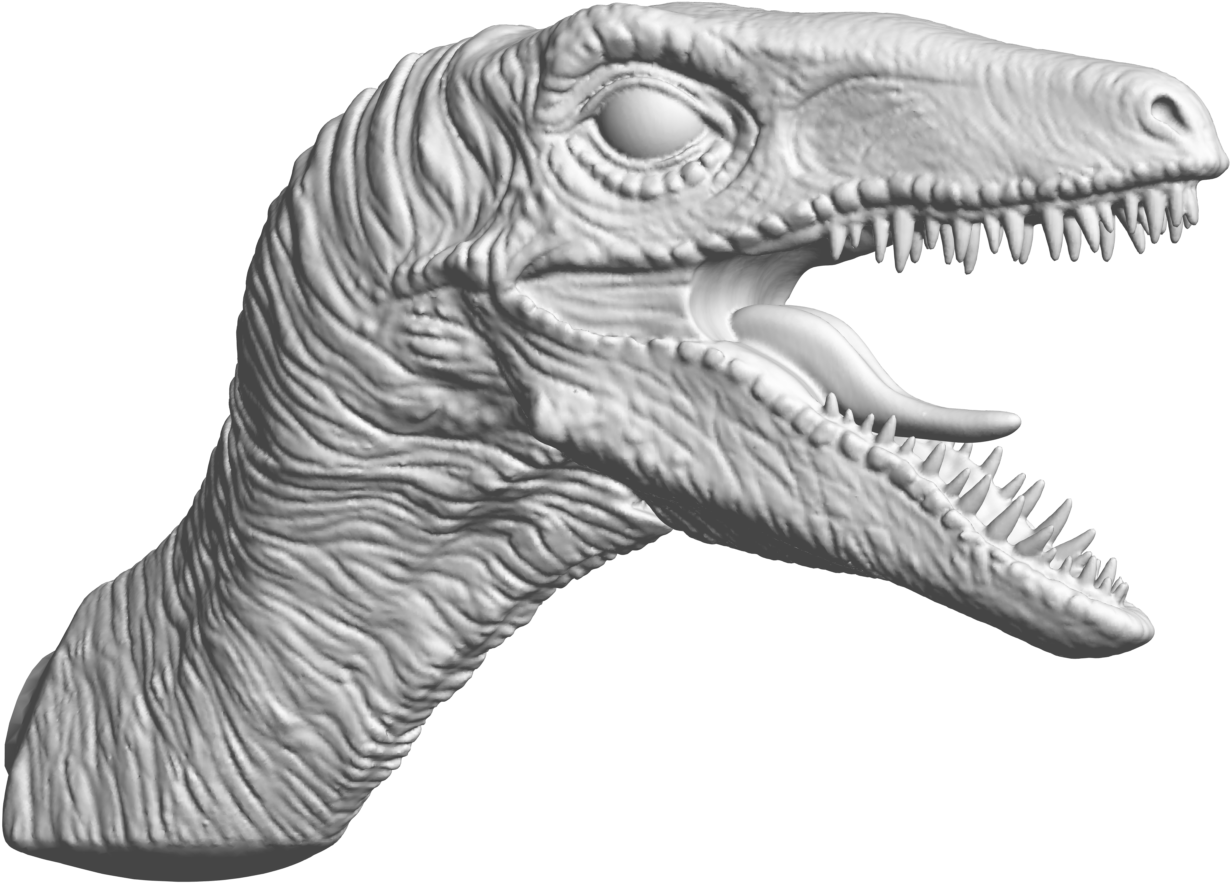} & \includegraphics[width=0.25\textwidth]{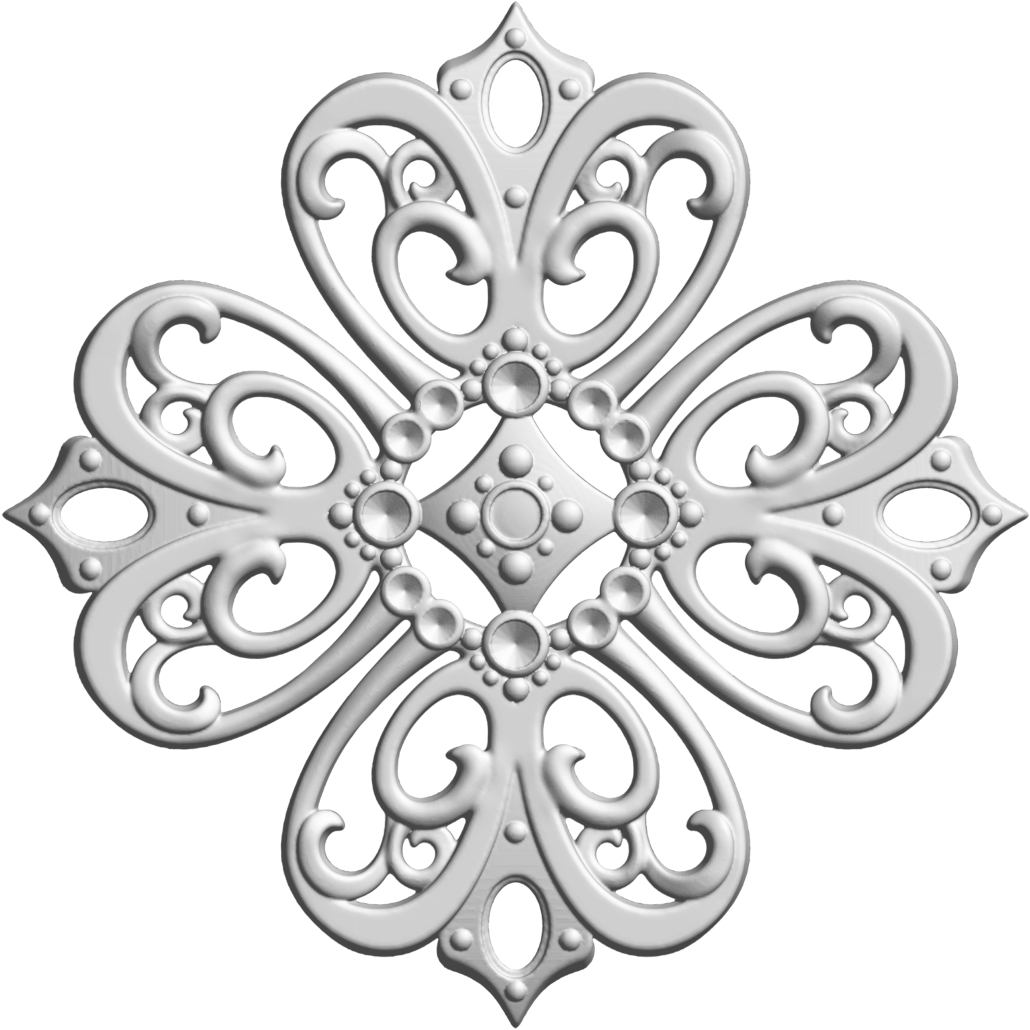} & 
\includegraphics[width=0.3\textwidth]{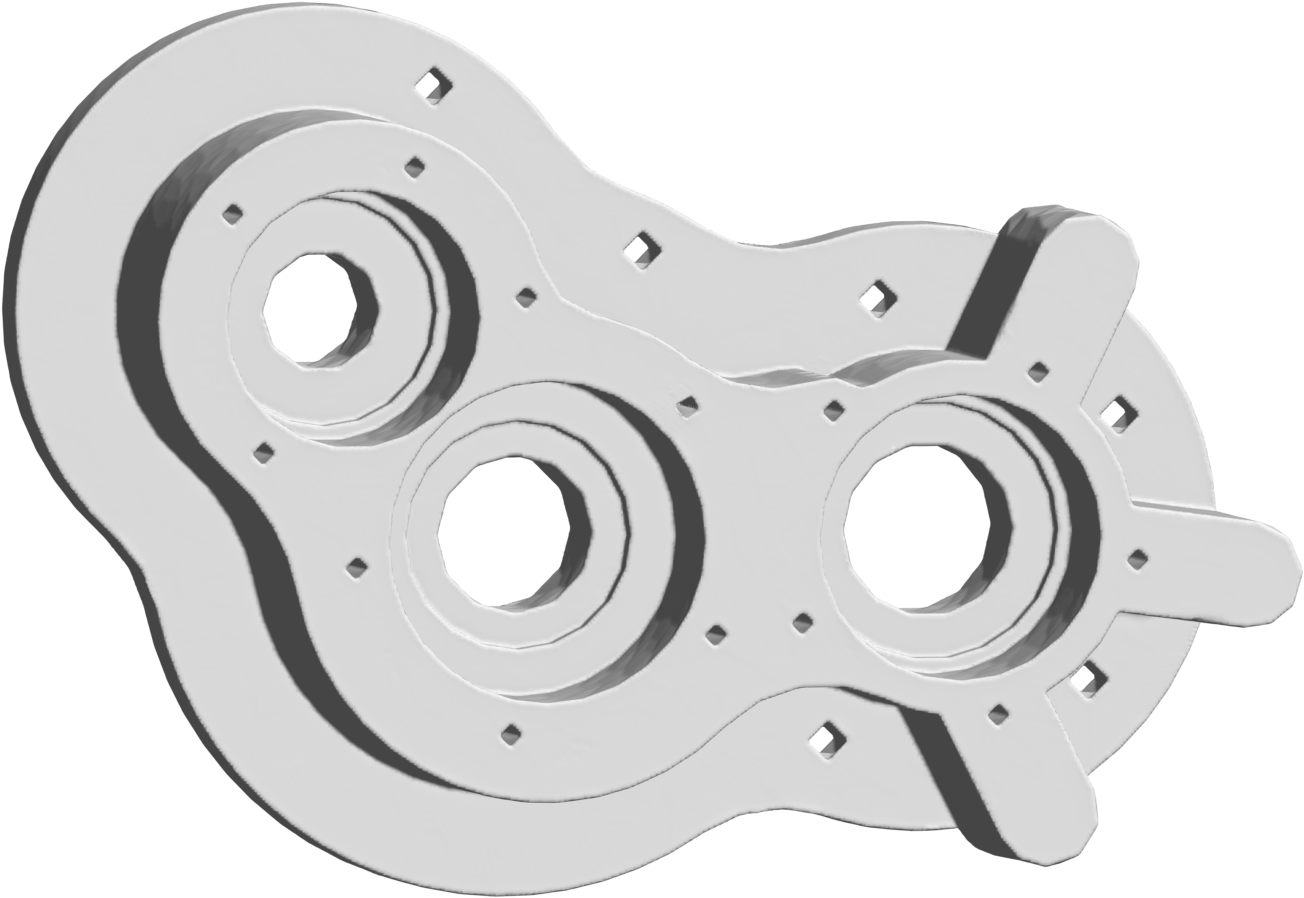}
\end{tabular}
\caption{\revision{Surface reconstructions for the CFPU (top row) and SPR (bottom row) methods for three different surfaces produced from CAD models (Raptor, Filigree, and Pump Carter). \label{fig:cadmodels}}}
\setlength{\tabcolsep}{6pt}
\end{figure}

In the first set of experiments, we consider reconstructing three different CAD models, the Raptor, Filigree, and Pump Carter, using only the vertices and normals from the triangulations that define the models.  In the case of the Pump Carter, the point cloud is obtained from a uniform subdivision of the original CAD model.  The results of the CFPU reconstructions are displayed in \revision{the top row of Figure \ref{fig:cadmodels}, while the bottom row displays the SPR reconstructions.  No regularization was used for the CFPU results and the default parameters were used for the SPR results. We can see from the figure that CFPU and SPR produces very similar results.  The fine-scale detail of these models is captured and there are no spurious sheets.  Furthermore, the sharp features present in the models are well captured without over smoothing.}

\begin{figure}[htb]
\setlength{\tabcolsep}{2pt}
\centering
\begin{tabular}{ccccc}
\rotatebox{90}{\hspace{0.15\textwidth} \small \revision{CFPU}} & \phantom{H} &
\includegraphics[width=0.3\textwidth]{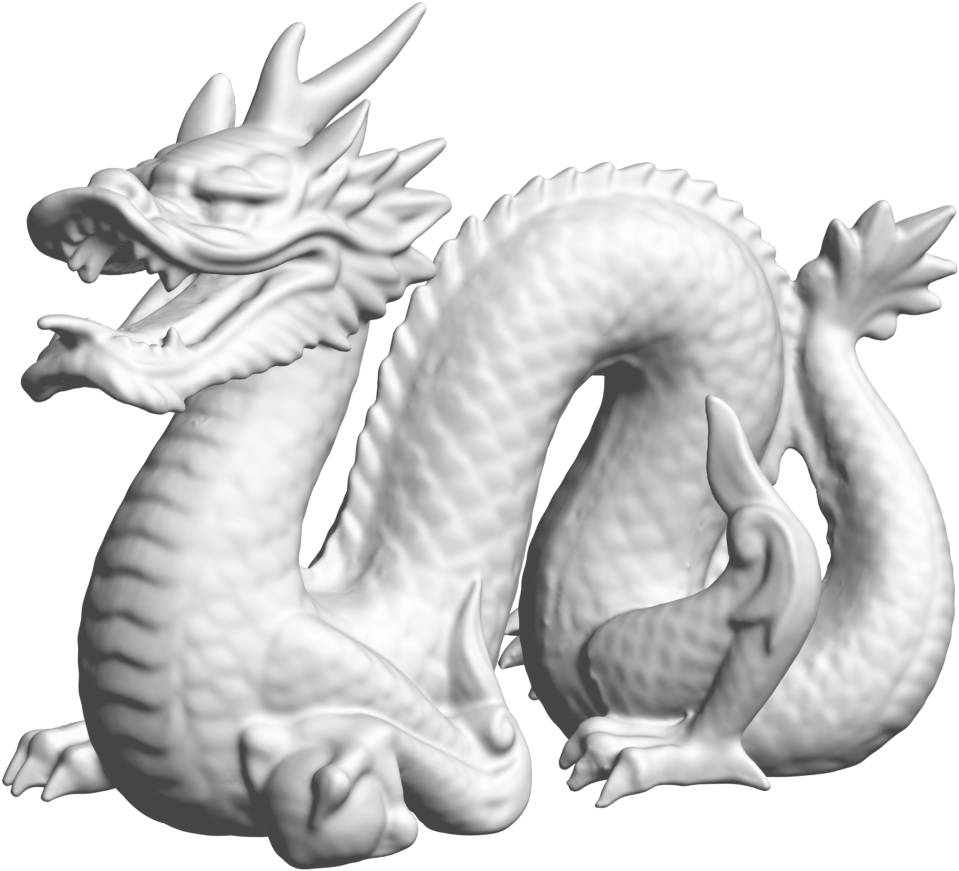} & \includegraphics[width=0.16\textwidth]{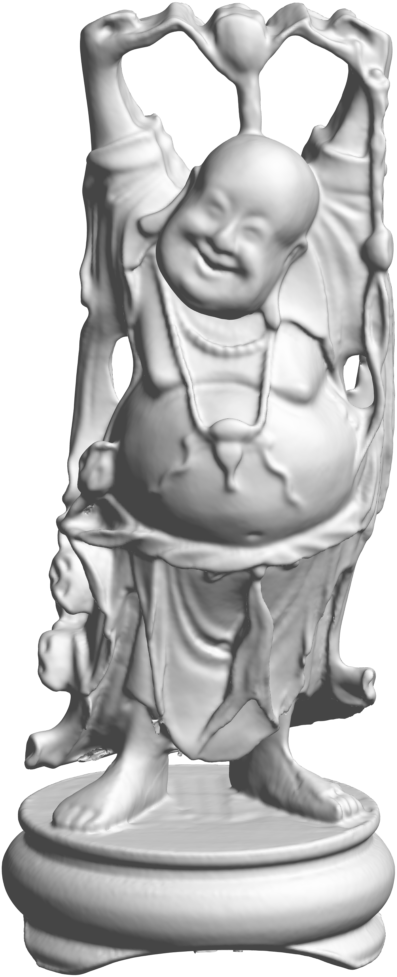} & 
\includegraphics[width=0.3\textwidth]{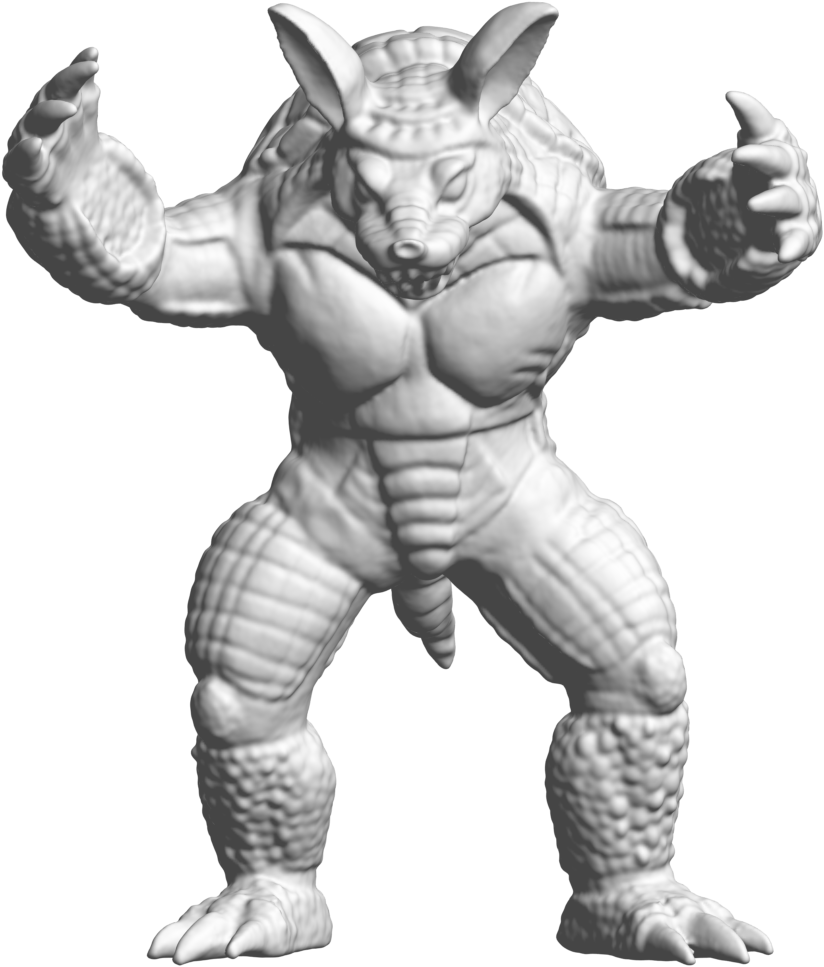}  \\
\rotatebox{90}{\hspace{0.15\textwidth} \small \revision{SPR}}  & \phantom{H} &
\includegraphics[width=0.3\textwidth]{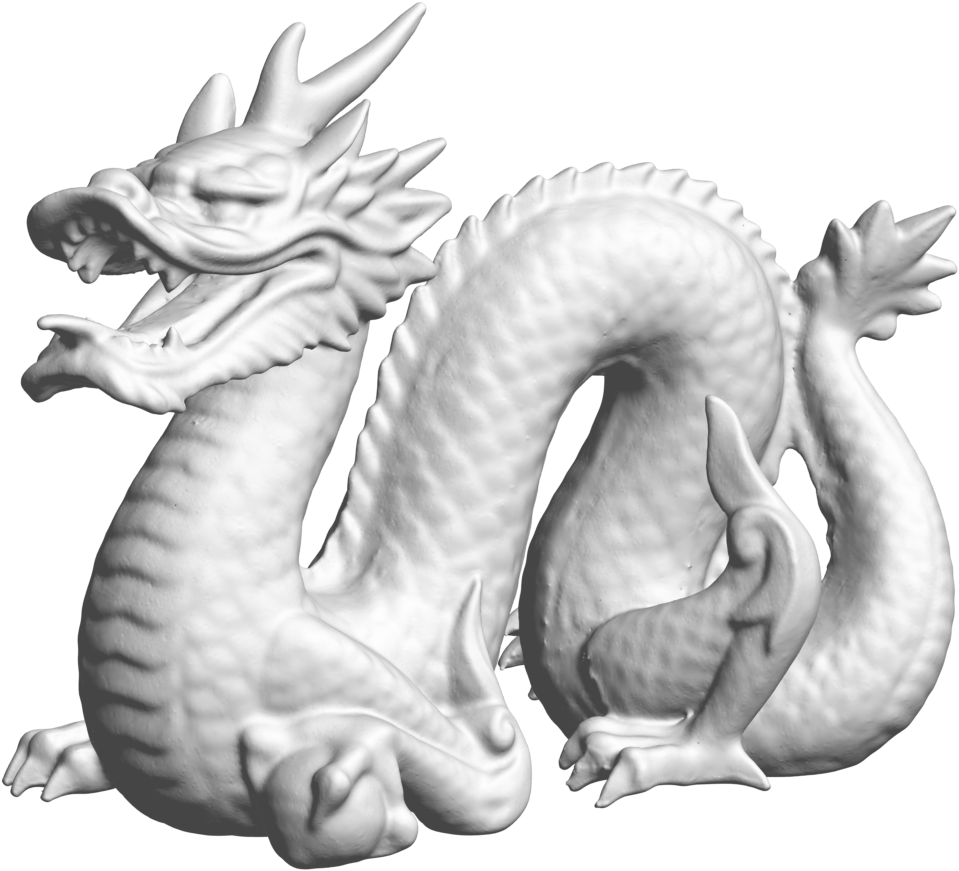} & \includegraphics[width=0.16\textwidth]{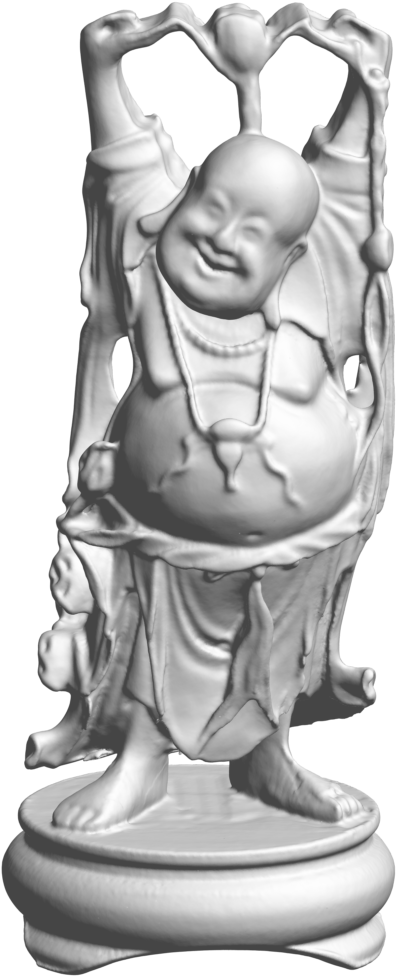} & 
\includegraphics[width=0.3\textwidth]{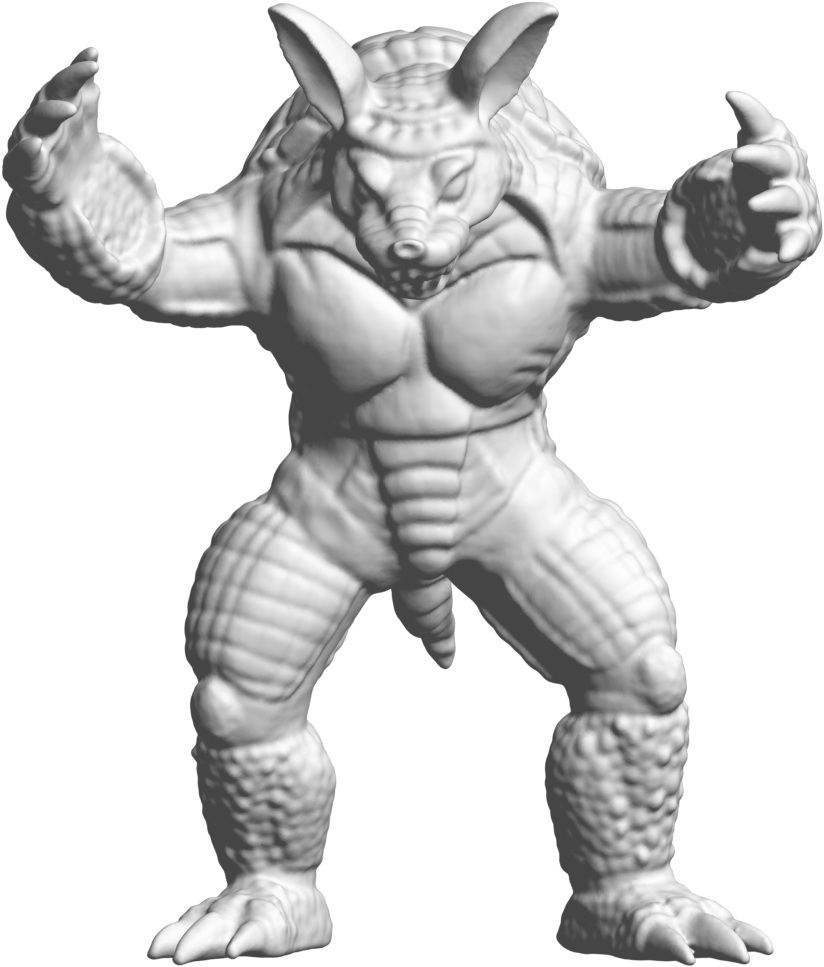}
\end{tabular}
\caption{\revision{Surface reconstructions for the CFPU (top row) and SPR (bottom row) methods for three different surfaces produced from cleaned-up scans of real models (Stanford Dragon, Happy Buddha, and Armadillo).\label{fig:realmodels}}}
\setlength{\tabcolsep}{6pt}
\end{figure}

For the next set of experiments, we consider three models produced from scans of real objects, the Stanford Dragon, Happy Buddha, and Armadillo.  In these tests, we use point clouds generated from post-processed versions of the scans where any outliers in the range data have been removed and the point clouds have been simplified.  \revision{The top row of Figure \ref{fig:realmodels} shows the resulting CFPU reconstructions of these models, while the bottom row shows the SPR reconstructions.  As with the CAD models, no regularization was used for the CFPU results and the default parameters were used for the SPR results.  We can see from the figure that the CFPU and SPR methods again produce similar results. The fine-scale detail is well-captured, without any noticeable smoothing, and there are no spurious sheets.} 
%Furthermore, the resulting reconstructions are visually very similar to other reconstructions from state-of-the-art methods found in the literature (e.g.,~\cite{

\begin{table}[htb]
\centering
{\color{black}
\begin{tabular}{|c|c||c|c|c|c|}
\hline
   \multirow{3}{*}{Surface}  & \multirow{3}{*}{Grid size} & \multicolumn{4}{c|}{Time in seconds} \\
\cline{3-6}
   & & \multicolumn{4}{c|}{Number of cores} \\
 &  & 1 & 2 & 4 & 8 \\
\hline\hline
Raptor & 384\texttimes143\texttimes 258 & 18.2 (16.5) & 6.67 (8.84) &  4.14 (7.22) & 3.35 (5.69) \\
Filigree & 383\texttimes63\texttimes384 & 30.4 (7.70) & 10.9 (4.10) & 6.63 (3.39) & 4.99 (3.13)   \\
Pump Carter & 384\texttimes85\texttimes271 & 18.1 (12.4) & 6.79 (6.01) & 4.22 (4.68) & 3.55 (4.11)  \\
Stanford Dragon & 384\texttimes175\texttimes272 & 42.73 (11.83) & 17.24 (6.46) & 10.7 (5.46) & 9.20 (4.82)  \\
Happy Buddha & 161\texttimes162\texttimes384 & 37.34 (7.71) & 14.0 (4.20) & 8.32 (3.43) & 5.99 (3.02)  \\
Armadillo & 323\texttimes294\texttimes384 & 9.42 (10.0) & 3.62 (5.54) & 2.25 (4.85) & 1.61 (4.55) \\
\hline
\end{tabular}}
\caption{\revision{Runtime performance of the CFPU method for reconstructing various surfaces.  The first number in each of the timing columns is the time it takes to solve all the local linear systems for determining the expansion coefficients of every $\icfpot_m$ in \eqref{eq:curlfree_icf_pum}.  The number in parenthesis is the time it takes to evaluate the potential $s$ in \eqref{eq:curlfree_icf_pum} over the grid points to extract the level surface.  Refer to Table \ref{tbl:parameters} for further details on the parameters used.\label{tbl:runtime}}}
\end{table}

%\begin{table}[htb]
%\centering
%{\color{blue}
%\begin{tabular}{|c||c|c|}
%\hline
%Surface & Fitting time (sec.) & Evaluation time (sec.) \\
%\hline\hline
%Raptor & 3.35 & 5.69 \\
%Filigree & 4.99 & 3.13 \\
%Pump Carter & 3.55 & 4.11 \\
%Stanford Dragon & 9.20 & 4.82 \\
%Happy Buddha & 5.99 & 3.02 \\
%Armadillo & 1.61 & 4.55 \\
%\hline
%\end{tabular}}
%\caption{\revision{Runtimes for the CFPU method for reconstructing various surfaces.\label{tbl:runtime}}}
%\end{table}

%\setlength{\tabcolsep}{0pt}
\begin{figure}[htb]
\centering
\begin{tabular}{ccc}
& $\lambda=0$, $\alpha=0$ & $\lambda=0$, $\alpha=0$ \\
\rotatebox{90}{\hspace{0.05\textwidth} \small No regularization} & 
\includegraphics[width=0.44\textwidth]{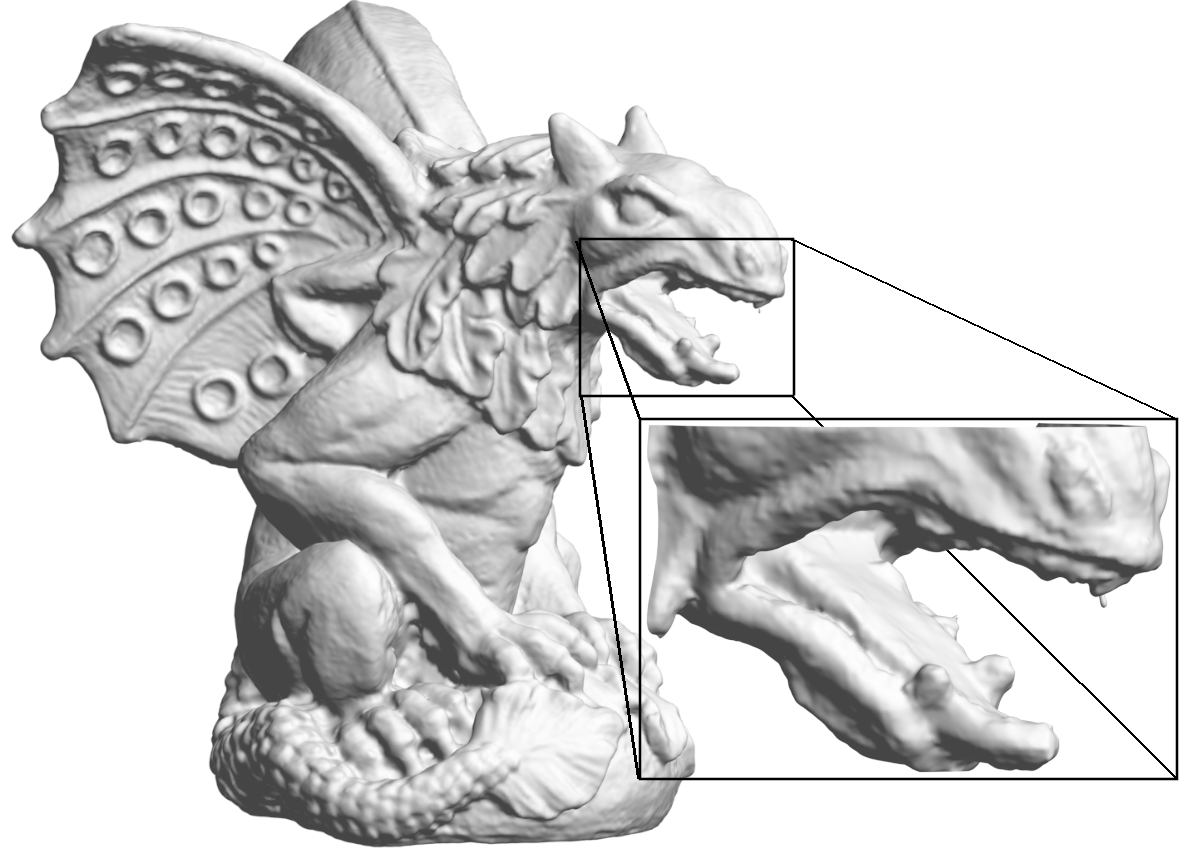} & \includegraphics[width=0.44\textwidth]{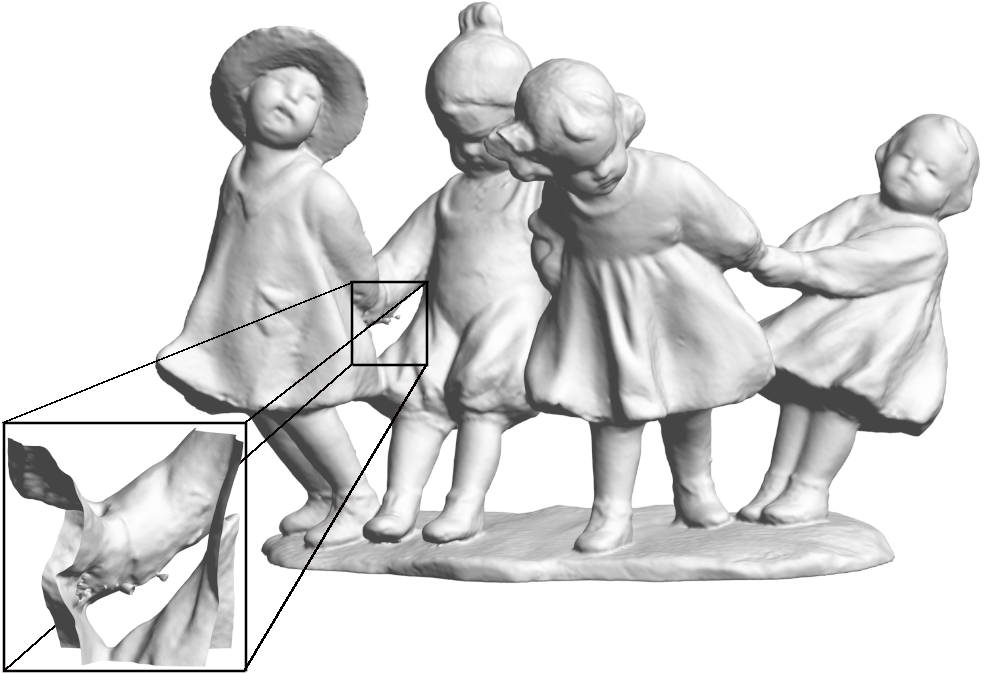} \\ 
& $\lambda=10^{-3}$, $\alpha=0$ & $\lambda=10^{-3}$, $\alpha=10^{-4}$ \\
\rotatebox{90}{\hspace{0.05\textwidth} \small Uniform regularization} & 
\includegraphics[width=0.44\textwidth]{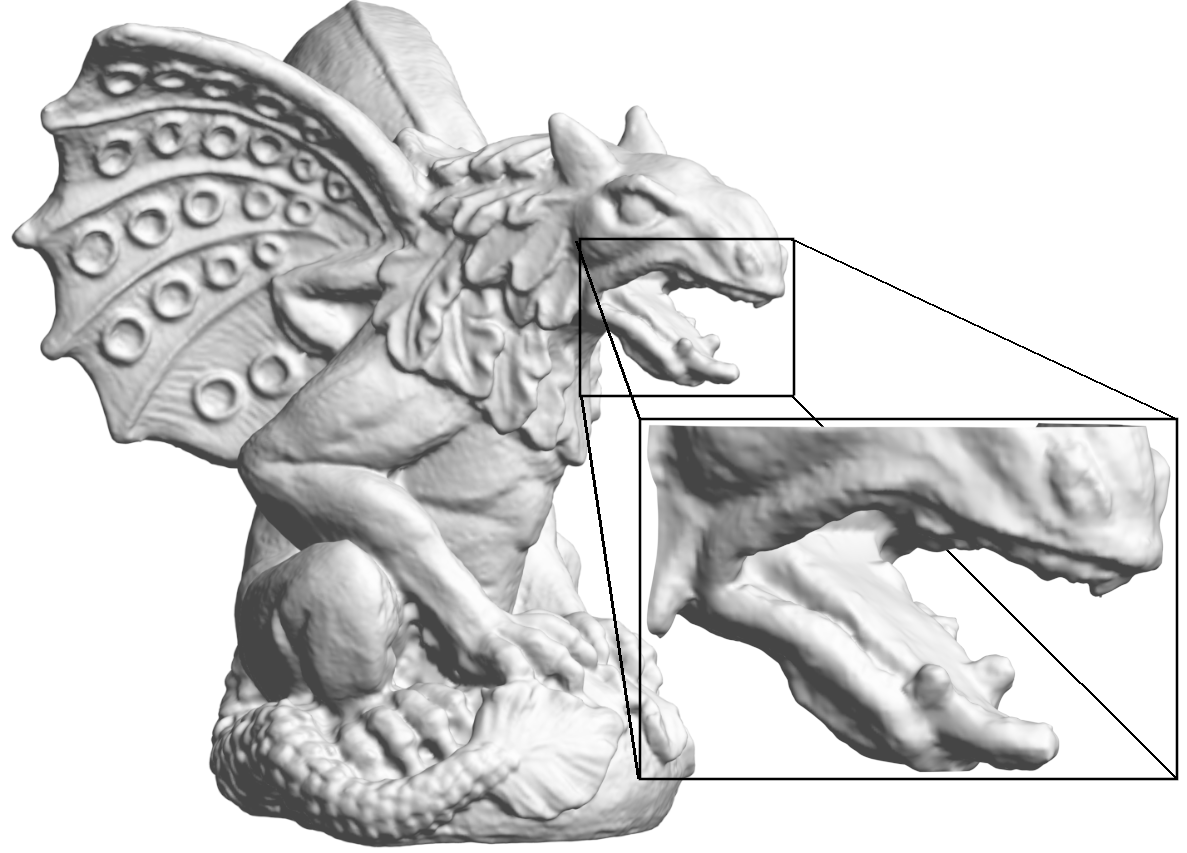} & \includegraphics[width=0.44\textwidth]{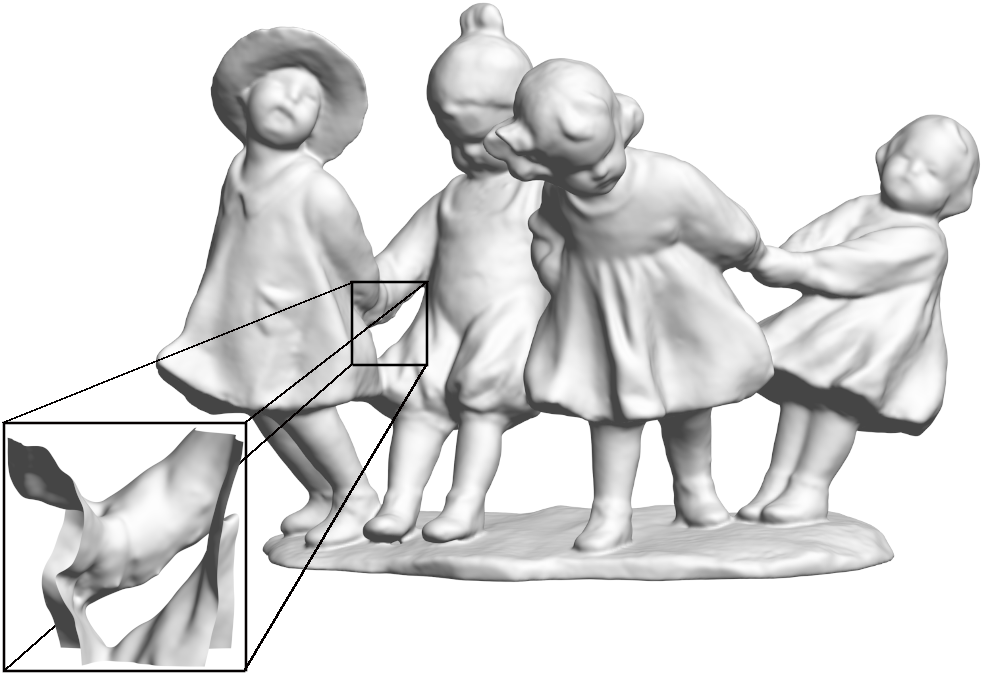} \\
& $\lambda=10^{-3}$, $\alpha=0$ & $\lambda=10^{-3}$, $\alpha=10^{-4}$ \\
\rotatebox{90}{\hspace{0.05\textwidth} \small Localized regularization} & 
\includegraphics[width=0.44\textwidth]{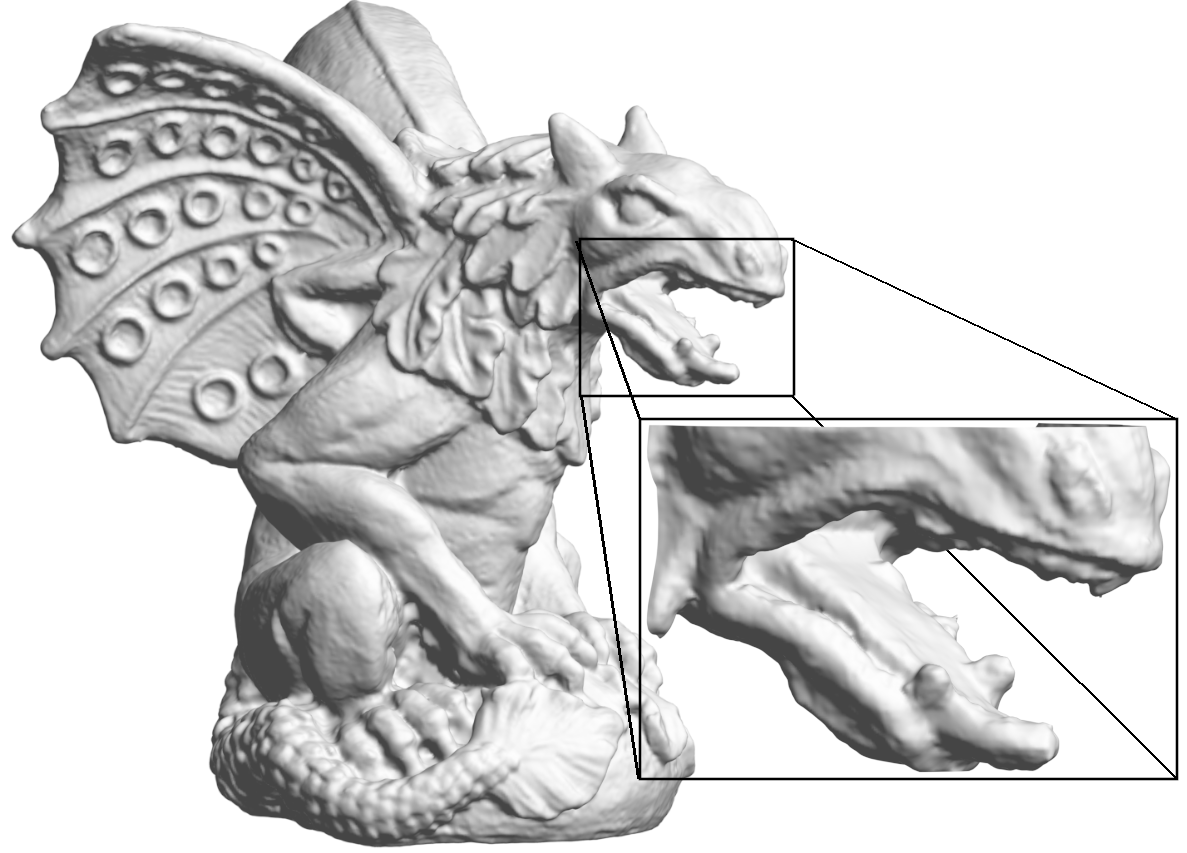}  & \includegraphics[width=0.44\textwidth]{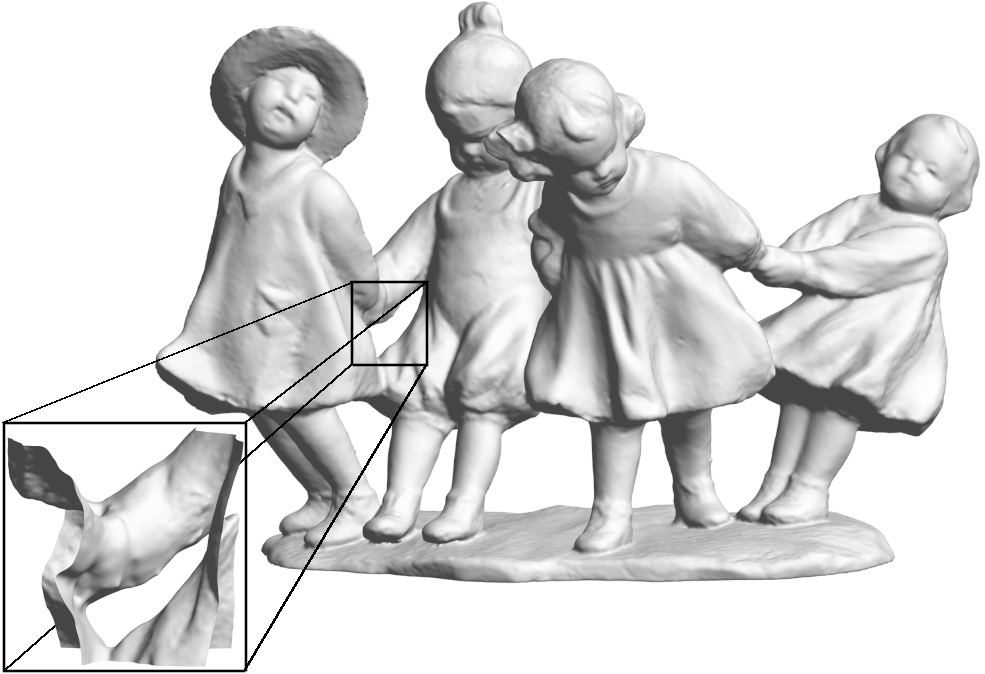} \\
\end{tabular}
\caption{CFPU reconstructions of the Gargoyle and Dancing Children models for $\ell=1$.  The top row shows the
results when no regularization is used, the second row shows the results when the same regularization is used across all patches, and the third row shows the results when using regularization only locally where issues arise in the original reconstructions.  The insets of the figures show a more detailed view of the issues.\label{fig:gargoylechildren}}
\end{figure}

\begin{figure}[htb]
\centering
\begin{tabular}{cc}
\includegraphics[width=0.44\textwidth]{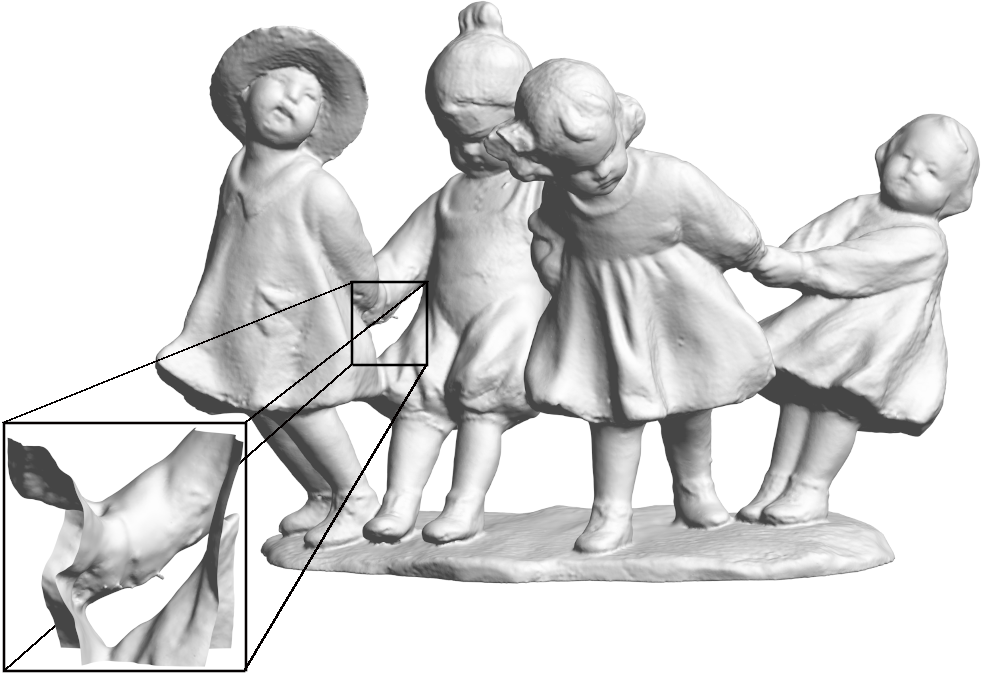} & 
\includegraphics[width=0.44\textwidth]{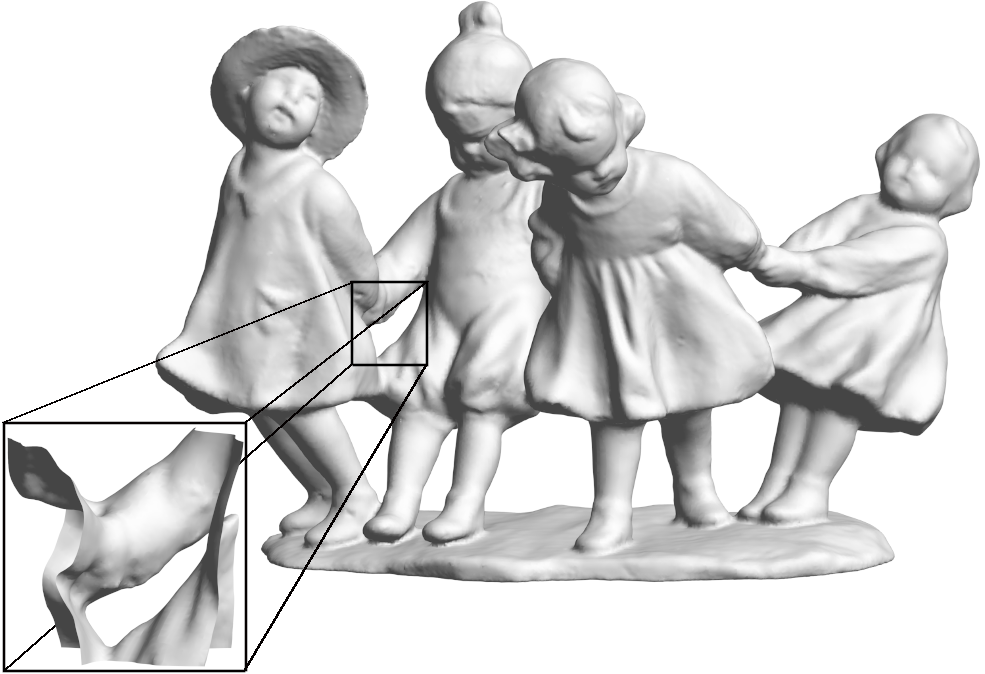} \\
\end{tabular}
\caption{\revision{SPR results for the Dancing Children model using the default parameters (left) and tuned parameters to resolve the issues around the joined hands.\label{fig:poissonchildren}}}
\end{figure}

In the final set of experiments, we again consider post-processed models produced from scans of real objects.  However, unlike the previous experiments, the CFPU reconstructions produce unwanted effects in some very localized areas of these models.  Figure \ref{fig:gargoylechildren} displays the results for the Gargoyle and Dancing Children models.  The first row of the figure shows the CFPU reconstructions without any regularization.  The insets of these plots show a more detailed view of where some spurious sheets are produced in the reconstructions.  These are exactly in the regions where the original polygonal models also have issues. The next row of the figure shows the results when including regularization in the reconstructions across all the patches.  For the Gargoyle model, we found that the spurious sheets around the front teeth could be eliminated by simply regularizing the fit of the normal vectors.  However, for the Dancing Children model it was necessary include some regularization of both the normals and the residual to eliminate the spurious sheets near the joined hands.  As such, we see that more smoothing has been introduced in Dancing Children results from the non-regularized version.  

Since the issues in these models is localized to a very small area, it is possible to only use regularization in the patches that surround these areas.  The last row of Figure \ref{fig:gargoylechildren} shows the results of such a localized regularization strategy.  We can see from these plots that the spurious sheets have been eliminated with this approach while not over-smoothing the entire model.  This is especially prevalent in the Dancing Children results.

\revision{Similar issues to the non-regularized CFPU reconstruction of the Dancing Children model occur for the SPR method, as shown in the left display of Figure \ref{fig:poissonchildren}, which was obtained with the default parameters.  These issues can also be resolved by introducing some smoothing, as shown in the right display of this same figure.  This latter result was obtained by setting the \texttt{samplesPerNode} parameter to 7.  As with the CFPU method with global regularization, this results in over smoothing for the whole surface.  The localized regularization approach of the CFPU method offers an effective solution to this problem.}

%\setlength{\tabcolsep}{6pt}
%\begin{figure}[htb]
%\centering
%%\begin{tabular}{p{0.47\textwidth}p{0.47\textwidth}}
%\begin{tabular}{llll}
%\def\stackalignment{r}
%\bottominset{\includegraphics[width=0.18\textwidth,frame]{gargoyle_zoom_cfpu_order1_noreg_noregi.png}}{\includegraphics[width=0.25\textwidth,]{gargoyle_cfpu_order1_noreg_noregi.png}}{0pt}{-60pt} \hspace{0.1in}& & & 
%\def\stackalignment{r}
%\bottominset{\includegraphics[width=0.18\textwidth]{gargoyle_zoom_cfpu_order1_reg2_1em4_noregi.png}}{\includegraphics[width=0.25\textwidth]{gargoyle_cfpu_order1_reg2_1em4_noregi.png}}{0pt}{-60pt} \\
%(a) No regularization  & & & (b) Regularization
%\end{tabular}
%\caption{CFPU reconstructions of the Happy Buddha with (a) no regularization and (b) with regularization.  In (b) GCV was used to determine the regularization parameter on each patch.  Both experiments used the highest resolution zippered model of the dragon consisting of  $N=583079$ points and normals vectors and $M=14226$ patches.\label{fig:gargoyle}}
%\end{figure}

\subsection{\revision{Computational performance}}
\revision{
%We conclude with some results on the computational performance of the method.  
Timing results  (wall clock times) for the CFPU method are displayed in Table \ref{tbl:runtime}.  Here we have reported the timings for solving the linear systems associated with the potential coefficients in each patch and evaluating these potentials on a uniform grid for the extraction of the level surface using MATLAB's isosurface function (which is based on the Marching Cubes algorithm).  As discussed at the end of Section \ref{sec:pu_description}, the CFPU method is amenable to parallel implementations, so we also report timings with an increasing number of processor cores being used.  These results were obtained using MATLAB's Parallel Computing Toolbox.  Specifically, we used parfor loops over the patches to solve the linear systems and evaluate the potential.  All results were obtained using MATLAB R2020b on a MacBook Pro with 2.4 GHz 8-Core Intel Core i9 processor and 32 GB RAM.  

We see from the table that increasing the number of cores from 1 to 2 reduces the computational time considerably, but that these reductions are diminished as the number of cores continues to increase.  We believe that this due to non-optimal memory usage associated with parfor and anticipate that a C/C++ implementation using OpenMP will show better scaling.  Additionally, we anticipate that the evaluation times (the values in parenthesis) can be further improved by using an octree-based isosurface extraction algorithm, such as~\cite{Kazhdan:SGP:2007}.  \reftwo{A more detailed study on the computational performance of the method as the number of patches vary with increasing sample points for general curl-free RBF interpolants can be found in~\cite{DFW2020a}.}
}
 %%%%%%%%%%%%%%%%%%%%%%%%%%%%%%%%%%%%%%%%%%%%%%%%%%%%%%%%%%%%
%%%%%%%%%%%%%%%%%%%%%%%%%%%%%%%%%%%%%%%%%%%%%%%%%%%%%%%%%%%%%
\section{Discussion and concluding remarks}
In this work, we introduced the CFPU method, a novel method that combines curl-free RBFs and a PU framework, for reconstructing an implicit surface from oriented point cloud data.  
By using PHS kernels for the RBFs, the method avoids the (sometimes) tedious task of selecting good shape or support parameters.  However, the method does feature an (integer) smoothness parameter ($\ell$).   For point clouds with noise and/or sharp features, using $\ell=1$ appears to give good reconstructions of the underlying surface across a range of examples.  However, for smooth surfaces without noise in the point clouds, increasing this value appears to result in higher accuracy and convergence rates of the reconstructions. The PU framework allows the method to be computationally efficient for large point sets.  \revision{Additionally, numerical experiments show that it can represent sharp features of the surfaces.}

The CFPU method can provide exact interpolation of a point cloud.  Regularization can also be introduced to deal with noise in the normals and the point positions.  We found that these regularization methods were effective at eliminating spurious sheets in the zero-level surfaces for point clouds with these issues.  We tested the method on a synthetic model and demonstrated its convergence properties.  We also tested the method on raw range data and demonstrated the robustness of the method when using PHS of order $\ell=1$ with some regularization. \revision{Finally, we presented results for several common 3D models found in the literature.  Comparisons of CFPU with the popular SSR and SSD methods show that it can provide similar surface reconstructions of the common 3D models.  In the case of the Homer model, which features regions with sparse sampling,  CFPU produced smoother reconstructions than SSR and SSD.  We additionally found this to be the case for the knot model with synthetically generated noise added to the normals.}

\revision{The CFPU method is similar to HRBF Implicits~\cite{Macedo_HRBF}  and SPR~\cite{Pois_SR,Screen_Pois} in that it uses the oriented normals to determine a potential for the level surface approximating the point cloud.  However, there are a few key differences. First, unlike HRBF Implicits, CFPU does not fit the potential directly together with the normal vectors, but exploits properties of matrix-valued curl-free RBFs to extract a potential from a fit of just the normals.  It then corrects this potential with a scalar RBF approximant.  This two-stage process allows the method to introduce regularizations that can separately deal with noisy normals and noisy point samples. In the case where just the normals are noisy, these regularization can still produce a potential that interpolates the point cloud exactly.  Furthermore, the PU framework allows for spatially localized regularization to be \refone{applied}.  Second,  CFPU solves for the potential in a strong form instead of a weak form like SPR.  This means the normals for the surface do not need to be spread to a 3D background grid and integrated against test functions. The strong form has the added benefit that it can produce surfaces that exactly interpolate the point cloud.  Third, HRBF Implicits and SPR solve large (sparse) linear systems involving all the data in the point cloud to determine an approximating potential.  Alternatively, CFPU solves small, decoupled linear systems on each PU patch to determine a local approximating potential that is then blended through a sparse process to determine a global potential. This should make the method more amenable to parallel implementations.}

There are several aspects of CFPU that we plan to investigate in the future.  One of these is an automated approach for generating the PU patches for non-uniformly sampled point clouds.  Another, along these same lines, is to look at ways of adapting the shapes of the PU patches from balls to ellipsoids that better conform to the underlying surface.  We also plan to investigate automated approaches for selecting the regularization parameters.  While GCV worked well for synthetic noise, we found that it did not produce acceptable results for raw range data.  We plan to look into other techniques along the lines of~\cite{Lee2003} for selecting these parameters.  Another regularization approach we will investigate is based on least squares, where fewer RBF centers are used than data sites.  This approach has been effective in other RBF-PU methods~\cite{Larsson2017} and has the additional benefit of reducing the computational cost.  \revision{Finally, we plan to implement a parallel version of the  method in C/C++ for even further improvements in computational efficiency.}

 %%%%%%%%%%%%%%%%%%%%%%%%%%%%%%%%%%%%%%%%%%%%%%%%%%%%%%%%%%%%
%%%%%%%%%%%%%%%%%%%%%%%%%%%%%%%%%%%%%%%%%%%%%%%%%%%%%%%%%%%%%
\section*{Acknowledgements}
KPD's work was partially supported by a SMART Scholarship, which is funded by The Under Secretary of Defense-Research and Engineering, National Defense Education Program / BA-1, Basic Research. GBW's work was partially supported by National Science Foundation grants 1717556 and 1952674.  The Stanford Bunny, Happy Buddha, Stanford Dragon, and Armadillo data were obtained from the Stanford University 3D Scanning Repository.
The Homer, Raptor, Filigree, Pump Carter, Dancing Children, and Gargoyle data were obtained from the AIM@SHAPE-VISIONAIR Shape Repository.

\bibliographystyle{siamplain}
\bibliography{PUM_SR_Refs,references}

\end{document}